\def\version{\today}

\documentclass[reqno,twoside,11pt]{amsart}
\usepackage{cite}
\usepackage{amsmath}
\usepackage{amsfonts}
\usepackage{amssymb}
\usepackage{epsfig}
\usepackage{verbatim}
\usepackage{color}

\usepackage{ulem}

\usepackage{mathpazo}

\DeclareFontFamily{OT1}{cmss}{} \DeclareFontShape{OT1}{cmss}{m}{n} {<5> <6> <7> <8> <9> <10> <11> <12> <13> <14.4> cmss10}{}
\DeclareMathAlphabet{\cmss}{OT1}{cmss}{m}{n}

\DeclareFontFamily{OT1}{fraktura}{}
\DeclareFontShape{OT1}{fraktura}{m}{n} {<5> <6> <7> <8> <9> <10> <11> <12> <13> <14.4> [1.1] eufm10}{}
\DeclareMathAlphabet{\fraktura}{OT1}{fraktura}{m}{n}

\newtheoremstyle{thm}{1.5ex}{1.5ex}{\itshape\rmfamily}{} {\bfseries\rmfamily}{}{2ex}{}

\newtheoremstyle{def}{1.5ex}{1.5ex}{\rmfamily\sl}{} {\bfseries\rmfamily}{}{2ex}{}

\newtheoremstyle{rem}{1.3ex}{1.3ex}{\rmfamily}{} {\bfseries\rmfamily}{}{2ex}{}

\newtheoremstyle{ass}{1.5ex}{1.5ex}{\rmfamily\sl}{} {\bfseries\rmfamily}{}{2ex}{}

\newenvironment{proofsect}[1] {\vskip0.1cm\noindent{\rmfamily\itshape#1.}}{\qed\vspace{0.15cm}}

\theoremstyle{thm}
\newtheorem{theorem}{Theorem}[section]
\newtheorem{lemma}[theorem]{Lemma}
\newtheorem{proposition}[theorem]{Proposition}

\newtheorem*{Main Theorem}{Main Theorem.}
\newtheorem{corollary}[theorem]{Corollary}

\newtheorem{definition}[theorem]{Definition}

\theoremstyle{rem}

\numberwithin{equation}{section}

\renewcommand{\section}{\secdef\sct\sect}
\newcommand{\sct}[2][default]{\refstepcounter{section}
\addcontentsline{toc}{section}
{{\tocsection {}{\thesection}{\!\!\!\!#1\dotfill}}{}}
\vspace{0.7cm}
\centerline{ 
\scshape\arabic{section}.\ #1} \nopagebreak \vspace{0.2cm}}
\newcommand{\sect}[1]{
\vspace{0.4cm} \centerline{\large\scshape\rmfamily #1}
\vspace{0.2cm}}

\renewcommand{\subsection}{\secdef\subsct\sbsect}
\newcommand{\subsct}[2][default]{\refstepcounter{subsection}
\addcontentsline{toc}{subsection}
{{\tocsection{\!\!}{\hspace{1.2em}\thesubsection}{\!\!\!\!#1\dotfill}}{}}
\nopagebreak\vspace{0.45\baselineskip} {\flushleft\bf
\thesection.\arabic{subsection}~\bf #1.~}
\\*[3mm]\noindent
\nopagebreak}
\newcommand{\sbsect}[1]{
\vspace{0.1cm}\noindent
\textbf{#1.~}\vspace{0.1cm}}

\renewcommand{\subsubsection}{%
\secdef \subsubsect\sbsbsect}
\newcommand{\subsubsect}[2][default]{%
\refstepcounter{subsubsection} 
\addcontentsline{toc}{subsubsection}{{\tocsection{\!\!}
{\hspace{3.05em}\thesubsubsection}{\!\!\!\!#1\dotfill}}{}}
\nopagebreak
\vspace{0.15\baselineskip} \nopagebreak {\flushleft\rmfamily
\itshape\arabic{section}.\arabic{subsection}.\arabic{subsubsection}
\ \rmfamily #1\/.}\ }
\newcommand{\sbsbsect}[1]{\vspace{0.1cm}\noindent
\rmfamily \itshape
\arabic{section}.\arabic{subsection}.\arabic{subsubsection} \
\sffamily #1\/.\ }

\renewcommand{\caption}[1]{%
\vglue0.5cm
\refstepcounter{figure}
\begin{center}
\begin{minipage}[c]{0.8\textwidth}\small {\sc Fig.~\thefigure\ }#1\end{minipage}
\end{center}
}

\newcommand{\supp}{\operatorname{supp}}

\newcommand{\textd}{\text{\rm d}\mkern0.5mu}
\newcommand{\texti}{\text{\rm  i}\mkern0.7mu}
\newcommand{\texte}{\text{\rm  e}\mkern0.7mu}

\newcommand{\Var}{\text{\rm Var}}

\newcommand{\EE}{\mathcal E}
\newcommand{\FF}{\mathcal F}
\newcommand{\GG}{\mathcal G}
\newcommand{\HH}{\mathcal H}

\newcommand{\NN}{\mathcal N}
\newcommand{\OO}{\mathcal O}

\newcommand{\E}{\mathbb E}

\newcommand{\N}{\mathbb N}

\newcommand{\R}{\mathbb R}

\newcommand{\Z}{\mathbb Z}

\newcommand{\twoeqref}[2]{(\ref{#1}--\ref{#2})}

\newcommand{\cc}{{\text{\rm c}}}

\newcommand{\fraka}{\fraktura a}

\newcommand{\frake}{\fraktura e}
\newcommand{\frakd}{\fraktura d}

\def\myffrac#1#2 in #3{\raise 2.6pt\hbox{$#3 #1$}\mkern-1.5mu\raise 0.8pt\hbox{$#3/$}\mkern-1.1mu\lower 1.5pt\hbox{$#3 #2$}}
\newcommand{\ffrac}[2]{\mathchoice%
	{\myffrac{#1}{#2} in \scriptstyle}
	{\myffrac{#1}{#2} in \scriptstyle}
	{\myffrac{#1}{#2} in \scriptscriptstyle}
	{\myffrac{#1}{#2} in \scriptscriptstyle}
}

\newcommand{\wh}{\widehat}
\newcommand{\wt}{\widetilde}

\newcommand{\laweq}{\,\overset{\text{\rm law}}=\,}

\newcommand{\leb}{{\rm Leb}}
\newcommand{\Cov}{\text{\rm Cov}}
\newcommand{\dd}{\textd}

\newcommand{\Lawarrow}{{\,\overset{\text{\rm law}}\longrightarrow\,}}

\newcommand\independent{\protect\mathpalette{\protect\independenT}{\perp}}
\def\independenT#1#2{\mathrel{\rlap{$#1#2$}\mkern3mu{#1#2}}}

\newcommand{\MB}{}
\newcommand{\eMB}{\normalcolor}

\newcommand{\cspecial}{\fraktura c}

\newcommand{\myemph}[1]{\textit{#1}}

\begin{document}

\title[Random walk local time\hfill]{Exceptional points of discrete-time\\random walks in planar domains}
\author[\hfill Y.~Abe, M.~Biskup, S.~Lee]
{Yoshihiro Abe$^1$, Marek~Biskup$^{2}$ and Sangchul Lee$^{2}$}
\thanks{\hglue-4.5mm\fontsize{9.6}{9.6}\selectfont\copyright\,\textrm{2019}\ \ \textrm{Y. Abe, M.~Biskup and S.~Lee.
Reproduction, by any means, of the entire
article for non-commercial purposes is permitted without charge.\vspace{2mm}}}
\maketitle

\vspace{-5mm}
\centerline{\textit{$^1$
Department of Mathematics and Informatics, Chiba University, Chiba, Japan}}
\centerline{\textit{$^2$
Department of Mathematics, UCLA, Los Angeles, California, USA}}
\smallskip

\smallskip
\centerline{\version}

\vskip0.5cm
\begin{quote}
\footnotesize \textbf{Abstract:}
Given a sequence of lattice approximations $D_N\subset\mathbb Z^2$ of a bounded continuum domain~$D\subset\mathbb R^2$ with the vertices outside $D_N$ fused together into one boundary vertex~$\varrho$, we consider discrete-time simple random walks  in $D_N\cup\{\varrho\}$ run for a time proportional to the expected cover time and describe the scaling limit of the exceptional level sets of the thick, thin, light and avoided points.  We show that these are distributed, up a spatially-dependent log-normal factor, as the zero-average Liouville Quantum Gravity measures in~$D$. The limit law of the local time configuration at, and nearby, the exceptional points is determined as well. The results extend earlier work by the first two authors who analyzed the continuous-time problem in the parametrization by the local time at~$\varrho$. A novel uniqueness result concerning divisible random measures and, in particular, Gaussian Multiplicative Chaos, is derived as part of the proofs. \end{quote}


\section{Introduction}
\begin{figure}[t]
\vglue0.2cm
\centerline{\includegraphics[height = 2.0in]{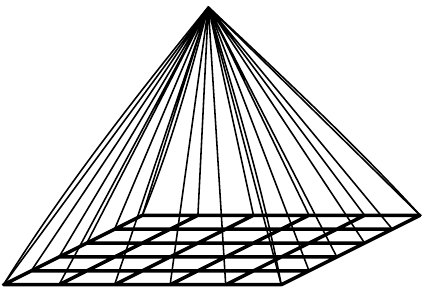}
}
\begin{quote}
\small 
\vglue-0.2cm
\caption{
\label{fig1}
The graph $(V\cup\{\varrho\},E)$ corresponding to~$D_N$ being the square of $6\times 6$ vertices and all edges emanating from~$D_N$ routed to the boundary vertex~$\varrho$. Note that the graph $(V\cup\{\varrho\},E)$ is planar whenever~$\Z^2\smallsetminus D_N$ is connected.}
\normalsize
\end{quote}
\end{figure}

\noindent
This note is a continuation of earlier work by the first two authors who in~\cite{AB} studied various exceptional  level  sets associated with the local time of random walks in lattice versions~$D_N\subset\Z^2$ of bounded open domains~$D\subset\R^2$, at times proportional to the cover time of~$D_N$. The walks in~\cite{AB} move as the ordinary constant-speed continuous-time simple symmetric random walk on~$D_N$ and, upon exit from~$D_N$, reenter~$D_N$ through a uniformly-chosen boundary edge. The re-entrance mechanism is conveniently realized by addition to~$D_N$ of a boundary vertex~$\varrho$ with all edges emanating out of~$D_N$ on~$\Z^2$ now ending in~$\varrho$. See Fig.~\ref{fig1} for an example. 

In~\cite{AB}, the local time was parametrized by the time spent at~$\varrho$. Through the use of the Second Ray-Knight Theorem (Eisenbaum, Kaspi, Marcus, Rosen and Shi~\cite{EKMRS}) this enabled a connection to the level sets of the Discrete Gaussian Free Field (DGFF) studied earlier by the second author and O.~Louidor~\cite{BL4}.
The goal of the present paper is to extend the results of~\cite{AB} to the more natural setting of a discrete-time random walk parametrized by its actual time. As we shall see, a close connection to the DGFF still persists, albeit now to that conditioned on vanishing arithmetic mean over~$D_N$. As no version of the Second Ray-Knight Theorem seems available for this specific setting, we have to proceed by suitable, and sometimes tedious, approximations. A key point is to control the fluctuations of the total time of the random walk at a given occupation time of the boundary vertex.

In order to give the precise setting of our problem, we first consider a general finite, unoriented, connected graph $G = (V\cup\{\varrho\},E)$, where~$\varrho$ is a distinguished vertex (not belonging to~$V$). Let~$X$ denote a sample path of the simple random walk on~$G$; i.e., a discrete-time Markov chain on~$V\cup\{\varrho\}$ with the transition probabilities
\begin{equation}
\cmss P (u, v) := 
\begin{cases} 
\frac1{\deg(u)}, &~\text{if}~e:=(u, v) \in E, \\
0, &~\text{otherwise},
\end{cases}
\end{equation}
where $\deg(u)$ is the degree of~$u$. As usual, we will write~$P^u$ to denote the law of~$X$ subject to the initial condition~$P^u(X_0=u)=1$.

Given a path~$X$ of the chain, the local time at~$v\in V\cup\{\varrho\}$ at time~$n$ is then given by 
\begin{equation} 
\label{E:local_time}
\ell_n^V(v) := \frac{1}{\deg(v)} \sum_{k=0}^n1_{\{X_k = v \}},\quad n \geq 0.
\end{equation}
Our aim is to observe the Markov chain at times when most, or even all, of the vertices have already been visited. This requires looking at the chain at times (at least) proportional to the total degree~$\deg(V):=\sum_{v\in V\cup\{\varrho\}}\deg(V)$. To simplify our later notations, we thus abbreviate, for any~$t>0$,
\begin{equation}
\label{E:LVt}
 L_t^V(v):=\ell_{\lfloor t\deg(V)\rfloor}^V(v),\quad v\in V.
\end{equation}
In this parametrization, we have $L_t^V(v)= t+o(t)$ with high probability as~$t\to\infty$.

Our derivations will make heavy use of the connection between the above Markov chain and an instance of the Discrete Gaussian Free Field (DGFF). Denoting by
\begin{equation}
H_{v} := \inf \bigl\{n \geq 0 \colon X_n = v \bigr\}
\end{equation}
the first hitting time of vertex~$v$, this DGFF is the centered Gaussian process $\{h_v^{V}\colon v \in V\}$ with covariances given by
\begin{equation}
\label{E:cov}
\E\bigl(h_u^{V} h_v^{V} \bigr) = G^{V} (u, v) :=E^u\bigl(\ell_{H_{\varrho}}^V(v)\bigr),
\end{equation} 
where~$\E$ the expectation with respect to the law of~$h^V$ and~$G^V$ is the Green function. The field naturally extends to~$\varrho$ by~$h^V_\varrho=0$. 

We will apply the above to~$V$ ranging through a sequence of lattice approximations of a well-behaved continuum domain. The following definitions are taken from~\cite{BL2}:

\begin{definition} 
An admissible domain 
is a bounded open subset of $\mathbb{R}^2$ 
that consists of a finite number of connected components
and whose boundary is composed of a finite number of connected sets each of which has 
positive Euclidean diameter.
\end{definition} 

We will write $\mathfrak{D}$ to denote the family of all admissible domains and let $d_{\infty} (\cdot, \cdot)$ denote the $\ell^{\infty}$-distance on $\mathbb{R}^2$. The lattice domains are then assumed to obey:

\begin{definition}
\label{dfn:admissible}
An admissible lattice approximation of $D \in \mathfrak{D}$ is a sequence $\{D_N\}_{N\ge1}$ of sets $D_N\subset\mathbb{Z}^2$ such that the following holds: There is~$N_0\in\N$ such that for all~$N\ge N_0$ we have
\begin{equation}
\label{E:1.8i}
D_N \subseteq \Bigl\{x \in \mathbb{Z}^2 \colon
d_{\infty}\bigl(\ffrac{x}{N}, \mathbb{R}^2 \smallsetminus D\bigr) > \frac{1}{N} \Bigr\}
\end{equation}
and, for any~$\delta>0$ there is also~$N_1\in\N$ such that for all~$N\ge N_1$,
\begin{equation}
\label{E:1.8ii}
D_N \supseteq \bigl\{x \in \mathbb{Z}^2 \colon
d_{\infty} (\ffrac{x}{N}, \mathbb{R}^2 \smallsetminus D) > \delta \bigr\}.
\end{equation}
\end{definition} 

As shown in \cite[Appendix~A]{BL2}, the conditions \twoeqref{E:1.8i}{E:1.8ii} ensure that the discrete harmonic measure on~$D_N$ tends, under scaling of space by~$N$, weakly to the harmonic measure on~$D$. This yields a precise asymptotic expansion of the associated Green function; see \cite[Chapter~1]{B-notes}. In particular, we have $G^{D_N}(x,x)=g\log N+O(1)$ for
\begin{equation}
g:=\frac1{2\pi}
\end{equation}
whenever~$x$ is deep inside~$D_N$. (This is by a factor~$4$ smaller than the corresponding constant in~\cite{B-notes,BL2} due to a different normalization of the Green function.)

\section{Main results}
\noindent
Let us move to discussing our main results. We pick an admissible domain~$D\in\mathfrak D$ and a sequence of admissible lattice approximation~$\{D_N\}_{N\ge1}$ and consider these fixed throughout the rest of the derivations. 

\subsection{Setting the scales}
We begin by setting the scales for the time that the random walk is observed for and determining the range of values taken by the local time:

\begin{theorem}
\label{thm-minmax}
Let~$\{t_N\}_{N\ge1}$ be a positive sequence such that, for some~$\theta>0$,
\begin{equation}
\label{E:1.12}
\lim_{N\to\infty}\frac{t_N}{(\log N)^2}=2g\theta.
\end{equation}
Then for any choices of~$x_N\in D_N$, the following limits hold in $P^{x_N}$-probability: 
\begin{equation}
\label{E:max}
\frac1{(\log N)^2}\,\max_{x\in D_N} L^{D_N}_{t_N}(x)\,\,\,\underset{ N\to\infty}\longrightarrow\,\,\,2 g\bigl(\sqrt\theta+1\bigr)^2
\end{equation}
and
\begin{equation}
\label{E:min}
\frac1{(\log N)^2}\,\min_{x\in D_N} L_{t_N}^{D_N}(x)\,\,\,\underset{ N\to\infty}\longrightarrow\,\,\,2 g\bigl[(\sqrt\theta-1)\vee0\,\bigr]^2.
\end{equation}
\end{theorem}

The conclusion \eqref{E:min} indicates (and our later results on avoided points prove) that the choice $\theta:=1$ identifies the leading order of the \myemph{cover time} of~$D_N$ --- defined as the first time that every vertex of the graph  has been visited. The cover time is random but  it is typically  concentrated  (more precisely, whenever the maximal hitting time is much smaller than the expected cover time; see Aldous~\cite{Aldous}). The scaling \eqref{E:1.12} thus corresponds to the walk run for a~$\theta$-multiple of the cover time.

As it turns out, under \eqref{E:1.12}, the asymptotic $[2g\theta+o(1)](\log N)^2$ marks the value of~$L_{t_N}^{D_N}$ at all but a vanishing fraction of the vertices in~$D_N$. In light of \twoeqref{E:max}{E:min}, this suggests that we call~$x\in D_N$ a \myemph{$\lambda$-thick point} if (for~$\lambda\in[0,1]$)
\begin{equation}
\label{E:2.4iu}
L_{t_N}^{D_N}(x)\ge 2 g\bigl(\sqrt\theta+\lambda\bigr)^2(\log N)^2
\end{equation}
and a \myemph{$\lambda$-thin point} if (for~$\lambda\in[0,\sqrt\theta)$)
\begin{equation}
\label{E:2.5iu}
L_{t_N}^{D_N}(x)\le 2 g\bigl(\sqrt\theta-\lambda\bigr)^2(\log N)^2.
\end{equation}
One of our goals  is to describe the scaling limit of the sets of thick and thin points. This is best done via random measures of the form
\begin{equation}
\label{E:zetaND}
\zeta^D_N:=\frac1{W_N}\sum_{x\in D_N}\delta_{x/N}\otimes\delta_{(L_{t_N}^{D_N}(x)-a_N)/\sqrt{2a_N}}\,,
\end{equation}
where $a_N$ is a sequence with the asymptotic growth as the right-hand side of \twoeqref{E:2.4iu}{E:2.5iu} and~$W_N$ is a normalizing sequence. The specific choice of the normalization by $\sqrt{2a_N}$ reflects on the natural fluctuations of $L_{t_N}^{D_N}(x)$ (which turn out to be order $\log N$ even between nearest neighbors)   and captures best the connection to the corresponding object for the DGFF to be discussed next.

\subsection{Level sets of zero-average DGFF}
Recall that $h^{D_N}$ denotes a sample of the DGFF in~$D_N$.
As shown by Bolthausen, Deuschel and Giacomin~\cite{BDG}, the maximum of~$h^{D_N}$ is asymptotic to~$2\sqrt g\log N$ and so the $\lambda$-thick points are naturally defined as those where the field exceeds~$2\lambda\sqrt g\log N$. Allowing for sub-leading corrections, these are best captured by the random measure
\begin{equation}
\label{E:etaDGFF}
\eta_N^D:=\frac1{K_N}\sum_{x\in D_N}\delta_{x/N}\otimes\delta_{h^{D_N}_x-\wh a_N},
\end{equation}
where~$\{\wh a_N\}$ is a centering sequence with the asymptotic $\wh a_N\sim 2\lambda\sqrt{g}\log N$ and
\begin{equation}
\label{E:1.19e}
K_N:=\frac{N^2}{\sqrt{\log N}}\texte^{-\frac{(\wh a_N)^2}{2g\log N}}.
\end{equation}
 In \cite[Theorem~2.1]{BL4} it was shown that for each~$\lambda\in(0,1)$, there is a constant~$\cspecial(\lambda)>0$  (independent of~$D$ or the approximating sequence $\{D_N\}_{N\ge1}$)  such that, relative to the topology of vague convergence of measures on $\overline D\times(\R\cup\{+\infty\})$,
\begin{equation}
\label{E:1.19}
\eta_N^D\,\,\underset{N\to\infty}\Lawarrow\,\,\cspecial(\lambda)\,Z^D_\lambda(\textd x)\otimes\texte^{-\alpha\lambda h}\textd h,
\end{equation}
where
\begin{equation}
\alpha:=\frac2{\sqrt g}
\end{equation}
 and where~$Z_\lambda^D$ is a random a.s.-finite Borel measure in~$D$ called the \myemph{Liouville Quantum Gravity} (LQG)  at parameter~$\lambda$-times critical. The measure~$Z^D_\lambda$ is normalized so that, for each Borel set~$A\subseteq D$,
\begin{equation}
\label{E:1.19a}
\E Z_\lambda^D(A)=\int_A r^{\,D}(x)^{2\lambda^2}\textd x,
\end{equation} 
where~$r^D$ is an explicit bounded, continuous function supported on~$D$ that, for~$D$ simply connected, is the conformal radius; see~\cite[(2.10)]{BL4}. 

\begin{figure}[t]
\vglue0.2cm
\centerline{\includegraphics[width = 3.8in]{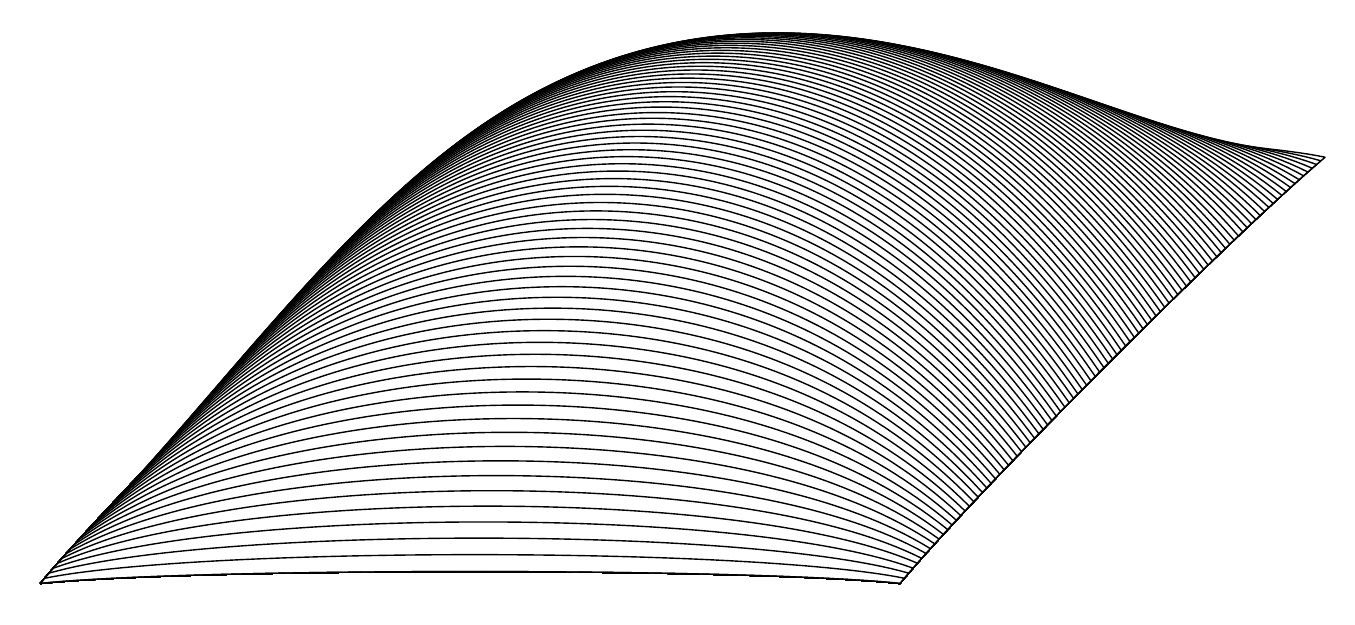}
}
\begin{quote}
\small 
\vglue-0.2cm
\caption{
\label{fig2}
A plot of function~$\frakd$ on~$D:=(0,1)^2$ obtained by solving the differential equation $-\Delta\frakd = \leb(D)/\sigma_D^2$, where~$\Delta$ is the Laplacian, with Dirichlet boundary conditions on~$\partial D$.}
\normalsize
\end{quote}
\end{figure}

As was shown in \cite{AB}, the measures $\{Z^D_\lambda\colon\lambda\in(0,1)\}$ are  quite relevant for the exceptional level sets associated with the continuous-time random walk in the parametrization by the local time spent in the ``boundary vertex.'' Somewhat different measures  will arise for the discrete-time random walk. Let $\Pi^D(x,\cdot)$ be the harmonic measure in~$D$ defined, e.g., as the exit distribution from~$D$ of a Brownian motion started at~$x$. The continuum Green function in~$D$ with Dirichlet boundary condition is then given by
\begin{equation}
\wh G^D(x,y):=-g\log|x-y|+g\int_{\partial D}\Pi^D(x,\textd z)\log|y-z|.
\end{equation}
Writing~$\leb$ for the Lebesgue measure on~$\R^2$, let~$\frakd\colon\R^2\to\R$ be defined by
\begin{equation}
\frakd(x):=\leb(D)\frac{\int_D\textd y\,\wh G^D(x,y)}{\int_{D\times D}\textd z\,\textd y\,
\wh G^D(z,y)}.
\end{equation}
As is readily checked, $\frakd$ is bounded and continuous, vanishes outside~$D$ and integrates to~$\leb(D)$ over~$D$. (We also have~$\frakd\ge0$ because~$\wh G^D\ge0$ and also that the Laplacian of~$\frakd$ is constant on~$D$ but that is of no consequence in the sequel.) See~Fig.~\ref{fig2}. We claim:

\begin{theorem}
\label{thm-DGFF}
For each $\lambda\in(0,1)$ and each~$D\in\mathfrak D$, there is a unique random measure $Z^{D,0}_\lambda$ on~$D$ such that, for any sequence $\{D_N\}_{N\ge1}$ of admissible approximations of~$D$ and any centering sequence  $\{\wh a_N\}_{N\ge1}$ satisfying $\wh a_N\sim 2\lambda\sqrt g\log N$ as~$N\to\infty$, 
\begin{equation}
\label{E:2.16new}
\Bigl(\eta^D_N\,\Big|\,\sum_{x\in D_N}h^{D_N}_x=0\Bigr)
\,\,\,\underset{N\to\infty}\Lawarrow\,\,\,\cspecial(\lambda)\,Z^{D,0}_\lambda(\textd x)\otimes\texte^{-\alpha\lambda h}\textd h,
\end{equation}
where~$\cspecial(\lambda)$ is as in \eqref{E:1.19}. Moreover,
if $Y$ is a normal random variable with mean zero and variance
\begin{equation}
\label{E:2.14new}
\sigma_D^2:=\int_{D\times D}\textd x\textd y\,\wh G^D(x,y),
\end{equation}
then the measure from \twoeqref{E:1.19}{E:1.19a} obeys
\begin{equation}
\label{E:2.15new}
Y\independent Z^{D,0}_\lambda\quad\Rightarrow\quad
Z^D_\lambda(\textd x)\laweq \texte^{\lambda\alpha\frakd(x)Y}
\,Z^{D,0}_\lambda(\textd x).
\end{equation}
The law of $Z^{D,0}_\lambda$ is determined uniquely by \eqref{E:2.15new}.
\end{theorem}

The existence of a random measure $Z^{D,0}_\lambda$ satisfying \eqref{E:2.15new} is part of the proof of \eqref{E:2.16new}. The uniqueness  of the decomposition \eqref{E:2.15new} holds quite generally and constitutes the main technical ingredient of the proof; see Theorem~\ref{lemma-8.1new}  which is of independent interest.  The known properties of~$Z^D_\lambda$ (see \cite[Theorem 2.3]{BL4}) imply that $Z^{D,0}_\lambda$ is a.s.-finite and charges every non-empty open subset of~$D$ a.s.

\subsection{Exceptional local-time sets}
We are now well equipped to state our results concerning the limits of the random measures \eqref{E:zetaND} for a given centering sequence~$\{a_N\}_{N\ge1}$ growing as the right-hand sides of \twoeqref{E:2.4iu}{E:2.5iu} and the normalizing sequence given by
\begin{equation}
\label{E:WN}
W_N:=\frac{N^2}{\sqrt{\log N}}\texte^{-\frac{(\sqrt{2t_N}-\sqrt{2a_N})^2}{2g\log N}}.
\end{equation}
For the thick points we then get:

\begin{theorem}[Thick points]
\label{thm-thick}
Suppose $\{t_N\}_{N\ge1}$ and $\{a_N\}_{N\ge1}$ are positive sequences such that,  for some $\theta>0$ and some $\lambda\in(0,1)$, \eqref{E:1.12} and
\begin{equation}
\label{E:1.20}
\lim_{N\to\infty}\frac{a_N}{(\log N)^2}=2g(\sqrt\theta+\lambda)^2
\end{equation}
hold true.
Then for any~$x_N\in D_N$ and for~$X$ sampled from~$P^{x_N}$, the measures~$\zeta^D_N$ in \eqref{E:zetaND} with~$W_N$ as in \eqref{E:WN} obey
\begin{equation}
\label{E:1.21dis}
\zeta^D_N\,\,\,\underset{N\to\infty}\Lawarrow\,\,\, \sqrt{\frac{\sqrt\theta}{\sqrt\theta+\lambda}}\,\texte^{-\alpha^2\lambda^2/16}\,\cspecial(\lambda) \,\texte^{\alpha\lambda (\frakd(x)-1)Y}Z_\lambda^{D,0}(\textd x)\otimes\texte^{-\alpha\lambda h}\textd h
\end{equation}
in the sense of vague convergence of measures on $\overline D\times(\R\cup\{+\infty\})$,
where $Y=\NN(0,\sigma_D^2)$ and~$Z^{D,0}_\lambda$ are independent 
and $\cspecial(\lambda)$ is as in \eqref{E:1.19}.
\end{theorem}

For the thin points, we similarly obtain:

\begin{theorem}[Thin points]
\label{thm-thin}
Suppose $\{t_N\}_{N\ge1}$ and $\{a_N\}_{N\ge1}$ are positive sequences such that, for some $\theta>0$ and some $\lambda\in(0,\sqrt\theta\wedge1)$, \eqref{E:1.12} and
\begin{equation}
\label{E:1.22}
\lim_{N\to\infty}\frac{a_N}{(\log N)^2}=2g(\sqrt\theta-\lambda)^2
\end{equation}
hold true.
Then for any~$x_N\in D_N$ and for~$X$ sampled from~$P^{x_N}$, the measures~$\zeta^D_N$ in \eqref{E:zetaND} with~$W_N$ as in \eqref{E:WN} obey
\begin{equation}
\label{E:1.23dis}
\zeta^D_N\,\,\,\underset{N\to\infty}\Lawarrow\,\,\, \sqrt{\frac{\sqrt\theta}{\sqrt\theta-\lambda}}\,\texte^{-\alpha^2\lambda^2/16}\,\cspecial(\lambda) \,\texte^{\alpha\lambda (\frakd(x)-1)Y}Z_\lambda^{D,0}(\textd x)\otimes\texte^{+\alpha\lambda h}\textd h
\end{equation}
in the sense of vague convergence of measures on $\overline D\times(\R\cup\{-\infty\})$,
where $Y=\NN(0,\sigma_D^2)$ and~$Z^{D,0}_\lambda$ are independent and $\cspecial(\lambda)$ is as in \eqref{E:1.19}.
\end{theorem}

The limiting spatial distribution of the $\lambda$-thick and $\lambda$-thin points (as well as the distribution of the total number of these points) is governed by the measure 
\begin{equation}
\label{E:2.21iu}
\texte^{\alpha\lambda (\frakd(x)-1)Y}Z_\lambda^{D,0}(\textd x).
\end{equation}
In light of \eqref{E:2.15new}, this is somewhere between the zero-average LQG $Z_\lambda^{D,0}$ and the ``ordinary'' LQG $Z^D_\lambda$, which appeared in the limit for the parametrization by the local time at~$\varrho$. The second component of the measure on the right of \eqref{E:1.21dis} and \eqref{E:1.23dis} is exactly as that for the DGFF \eqref{E:1.19}. This is due to the judicious scaling of the second component of~$\zeta^D_N$ by $\sqrt{2a_N}$ rather than just~$\log N$, as was done in \cite{AB}.

\smallskip
Apart from the thick and thin points, \cite{AB} studied also the sets of points where the local time is order unity, called the \myemph{light} points, and the points where the local time vanishes, called the \myemph{avoided} points. In both cases, the LQG measure that appears is for parameter $\lambda:=\sqrt\theta$ (and $\theta\in(0,1)$). The control extends to the discrete-time problem parametrized by the total time as well. We start with the light points:

\begin{theorem}[Light points]
\label{thm-light}
Suppose $\{t_N\}_{N\ge1}$ is a positive sequence such that \eqref{E:1.12} holds for some $\theta\in(0,1)$. For any~$x_N\in D_N$ and for~$X$ sampled from~$P^{x_N}$,  consider the measure
\begin{equation}
\label{E:varthetaND}
\vartheta^D_N:=\frac1{\wh W_N }\sum_{x\in D_N}\delta_{x/N}\otimes\delta_{L_{t_N}^{D_N}(x)},
\end{equation}
where
\begin{equation}
\label{E:1.31}
\wh W_N :=N^2\texte^{-\frac{t_N}{g\log N}}.
\end{equation}
Then, in the sense of vague convergence of measures on~$\overline D\times[0,\infty)$,
\begin{equation}
\label{E:2.22ii}
\vartheta^D_N\,\,\,\underset{N\to\infty}\Lawarrow\,\,\,  \sqrt{2\pi g}\,\cspecial(\sqrt\theta)\,\,\texte^{\alpha\sqrt\theta (\frakd(x)-1)Y}\,  Z_{\sqrt{\theta}\,}^{D,0}(\textd x)\otimes\mu(\textd h),
\end{equation}
where $\cspecial(\lambda)$ is as in \eqref{E:1.19}, $Y=\NN(0,\sigma_D^2)$ and~$Z^{D,0}_{\sqrt\theta}$ are independent and $\mu:=\sum_{n\ge0}q_n\delta_{n/4}$ for a sequence $\{q_n\colon n\ge0\}$ of non-negative numbers  determined uniquely by
\begin{equation}
\label{E:2.27}
\sum_{n\ge0}q_n(1+s/4)^{-n} = \texte^{\frac{\alpha^2\theta}{2s}},\quad s>0.
\end{equation}
\end{theorem}

That~$\mu$ is supported on $\frac14\N_0:=\{0,\frac14,\frac12,\frac34,1,\dots\}$ arises from the normalization  in~\eqref{E:local_time}.
From \eqref{E:2.22ii} we conclude that the number of the vertices of~$D_N$ visited exactly~$n$ times during the first
\begin{equation}
[8g\theta+o(1)](\log N)^2\deg(D_N)
\end{equation}
steps of the random walk is thus asymptotic to
\begin{equation}
\label{E:2.26}
q_n\,\Bigl[\sqrt{2\pi g}\,\cspecial(\sqrt\theta)\,\,\int_D \texte^{\alpha\sqrt\theta (\frakd(x)-1)Y}\,  Z_{\sqrt{\theta}\,}^{D,0}(\textd x)\Bigr]\,\wh W_N,
\end{equation}
jointly for all~$n\ge0$.
Noting that $q_0=1$, straightforward limit considerations show:

\begin{theorem}[Avoided points]
\label{thm-avoid}
Suppose $\{t_N\}_{N\ge1}$ is a sequence such that \eqref{E:1.12} holds for some~$\theta\in(0,1)$.  For any~$x_N\in D_N$ and for~$X$ sampled from~$P^{x_N}$,  consider the measure
\begin{equation}
\label{E:kappaND}
\kappa^D_N:=\frac1{\wh W_N }\sum_{x\in D_N}1_{\{L_{t_N}^{D_N}(x)=0\}}\,\delta_{x/N},
\end{equation}
where~$\wh W_N $ is as in \eqref{E:1.31}. Then, in the sense of vague convergence of measures on~$\overline D$,
\begin{equation}
\label{E:2.27dis}
\kappa^D_N\,\,\,\underset{N\to\infty}\Lawarrow\,\,\, \sqrt{2\pi g}\,\cspecial(\sqrt\theta)\,\,\texte^{\alpha\sqrt\theta (\frakd(x)-1)Y}\,  Z_{\sqrt{\theta}\,}^{D,0}(\textd x),
\end{equation}
where $Y=\NN(0,\sigma_D^2)$ and~$Z^{D,0}_{\sqrt\theta}$ are independent and $\cspecial(\lambda)$ is as in \eqref{E:1.19}.
\end{theorem}

The above theorems will be deduced from the corresponding statements for a conti\-nu\-ous-time variant of~$X$ observed for a fixed time of order $N^2(\log N)^2$ (see Propositions~\ref{thm-thick-cont},~\ref{thm-thin-cont}, \ref{thm-light-cont} and~\ref{thm-avoid-cont}). These statements are nearly identical to Theorems~\ref{thm-thick}--\ref{thm-avoid} above, respectively, except for the term $\texte^{-\alpha^2\lambda^2/16}$ in \eqref{E:1.21dis} and \eqref{E:1.23dis} that arises from the fluctuations of the (continuous-time) local time at points where the discrete-time local time is large, and the measure~$\mu$ in \eqref{E:2.22ii} which  gets replaced  (in Proposition~\ref{thm-light-cont}) by a continuous, and quite explicit, counterpart. 

The fixed-time results for continuous-time random walk will be inferred from the corresponding results in~\cite{AB} for the parametrization by the local time at~$\varrho$. The main difference is that the measure \eqref{E:2.21iu} gets replaced by the ``pure'' LQG~$Z^D_\lambda$.

\newcommand{\loc}{{\text{\rm loc}}}

\subsection{Local structure}
Similarly as in~\cite{AB}, we are also able to control the local structure of the above exceptional sets. For the thick and thin points, this is achieved by considering the measures on $D\times\R\times\R^{\Z^2}$ of the form
\begin{equation}
\label{E:zetaNDloc}
\zeta^{D,\loc}_N:=\frac1{W_N}\sum_{x\in D_N}\delta_{x/N}\otimes\delta_{(L_{t_N}^{D_N}(x)-a_N)/\sqrt{2a_N}}
\otimes\delta_{\{(L_{t_N}^{D_N}(x)-L_{t_N}^{D_N}(x+z))/\sqrt{2a_N}\colon z\in\Z^2\}},
\end{equation}
where the third coordinate captures the ``shape'' of the local-time configuration near every exceptional point.

In the parametrization by the local time at the boundary vertex, the asymptotic ``law'' of the third component in \eqref{E:zetaNDloc} turned out be that of the pinned DGFF (i.e., the DGFF in~$\Z^2\smallsetminus\{0\}$) reduced by a multiple of the potential kernel~$\fraka$. Here we note that, in our normalization, $\fraka$ is the unique non-negative function on~$\Z^2$ that is discrete harmonic on~$\Z^2\smallsetminus\{0\}$ and obeys $\fraka(0)=0$ and $\fraka(x) \sim g\log|x|+O(1)$ as~$|x|\to\infty$. The pinned DGFF~$\phi$ then has the covariance structure
\begin{equation}
\label{E:2.30}
\Cov(\phi_x,\phi_y)=\fraka(x)+\fraka(y)-\fraka(x-y).
\end{equation}
As it turns out, a different (albeit closely related) Gaussian process arises for the discrete-time walk parametrized by its total time:

\begin{theorem}[Local structure of thick/thin points]
\label{thm-thick-loc}
For the setting and under the conditions of Theorem~\ref{thm-thick}, relative to the vague topology of $\overline D\times(\R\cup\{+\infty\})\times\R^{\Z^2}$,
\begin{equation}
\zeta^{D,\loc}_N\,\,\,\underset{N\to\infty}\Lawarrow\,\,\,\zeta^D\otimes\nu_\lambda,
\end{equation}
where $\zeta^D$ is the measure on the right of \eqref{E:1.21dis} and $\nu_\lambda$ is the law of $\wt\phi+\alpha\lambda\fraka-\frac18\alpha\lambda1_{\{0\}^\cc}$, for~$\wt\phi$ a centered Gaussian process on~$\Z^2$ with covariances
\begin{equation}
\label{E:2.32}
\Cov(\wt\phi_x,\wt\phi_y)=\fraka(x)+\fraka(y)-\fraka(x-y)-\frac18\bigl[1-\delta_{x,0}-\delta_{y,0}+\delta_{x,y}\bigr].
\end{equation}
The same statement (relative to the vague topology on $\overline D\times(\R\cup\{-\infty\})\times\R^{\Z^2}$) holds for the setting of Theorem~\ref{thm-thin} except that~$\nu_\lambda$ is then the law of 
$\wt\phi-\alpha\lambda\fraka+\frac18\alpha\lambda1_{\{0\}^\cc}$.  
\end{theorem}

To demonstrate that $\wt\phi$ is indeed closely related to the pinned DGFF~$\phi$, we note that, for $\{n_z\colon z\in\Z^2\}$ i.i.d\ $\NN(0,\frac18)$ that are independent of~$\wt\phi$,
\begin{equation}
\label{E:2.33}
\{\phi_z\colon z\in\Z^d\}\,\laweq\,\{\wt\phi_z+n_0-n_z\colon z\in\Z^2\}.
\end{equation}
We will verify this relation, along with the fact that \eqref{E:2.32} is positive semidefinite and thus the covariance of a Gaussian process, in Lemma~\ref{lemma-8.4a}. The i.i.d.\ normals appear during a conversion from the continuous-time walk to its discrete-time counterpart. They represent the scaling limit of the fluctuations of the local time due to the random (i.i.d.\ exponential) nature of the jump times.

We will also address the local time structure in the vicinity of the avoided points. This is done by considering the measure on $D\times[0,\infty)^{\Z^2}$ defined by
\begin{equation}
\label{E:kappaNDloc}
\kappa^{D,\loc}_N:=\frac1{\wh W_N }\sum_{x\in D_N}1_{\{L_{t_N}^{D_N}(x)=0\}}\,\delta_{x/N}\otimes\delta_{\{L^{D_N}_{t_N}(x+z)\colon z\in\Z^2\}}.
\end{equation}
For reasons explained earlier, the measure is concentrated on $D\times(\frac14\N_0)^{\Z^2}$. 

Recall from \cite[Theorem~2.8]{AB} that, for the continuous-time random walk parametrized by the local time at the boundary vertex and observed at the time corresponding to~$\theta$-multiple of the cover time, the limit distribution of the local configuration is described by the law $\nu_{\theta}^{\text{\rm RI}}$ of the occupation-time field of random-interlacements at level~$u:=\pi\theta$.  This measure was constructed by Rodriguez~\cite[Theorems~3.3 and~4.2]{R18} (see \cite[Section~2.6]{AB} for a summary of the construction). For the discrete-time random walk parametrized by its total time we get a discrete-time counterpart of $\nu_{\theta}^{\text{\rm RI}}$:

\begin{theorem}[Local structure of avoided points]
\label{thm-avoid-loc}
For each~$u>0$, there is a unique Borel measure $\nu_{u}^{\text{\rm RI,\,dis}}$ on $[0,\infty)^{\Z^2}$ that is supported on $(\frac14\N_0)^{\Z^2}$ and obeys the following: For
\begin{enumerate}
\item[(1)] $\{\ell(z)\colon z\in\Z^2\}$ a sample from $\nu_{u}^{\text{\rm RI,\,dis}}$, and
\item[(2)] $\{\tau_{z,j}\colon z\in\Z^2,\,j\ge1\}$ independent i.i.d.\ Exponential(1),
\end{enumerate}
we have
\begin{equation}
\nu_{u}^{\text{\rm RI}} = \text{\rm law of }\Bigl\{\frac14\sum_{j=1}^{4\ell(z)}\tau_{z,j}\colon z\in\Z^2\Bigr\}.
\end{equation}
For the setting and under the conditions of Theorem~\ref{thm-avoid}, for each~$\theta\in(0,1)$ we then have
\begin{equation}
\kappa^{D,\loc}_N\,\,\,\underset{N\to\infty}\Lawarrow\,\,\,\kappa^D\otimes\nu_{\theta}^{\text{\rm RI,\,dis}}
\end{equation}
where $\kappa^D$ is the measure on the right of \eqref{E:2.27dis}.
\end{theorem}

Similarly as in \cite{AB}, we will not attempt to make statements concerning the local structure of the light points as that would require developing the corresponding extension of the above occupation-time measure to the situation when the local time at the origin does not vanish.

\subsection{Remarks}
We proceed with a couple of remarks. First note that, along with \eqref{E:min} and the fact that~$Z^{D,0}_{\sqrt\theta}$ is supported on all of~$D$ a.s., Theorem~\ref{thm-avoid} implies that the cover time is indeed marked by the choice~$\theta:=1$. Second, note  that 
an explicit formula for~$q_n$ can be extracted from \eqref{E:2.27}. This is achieved using the identity
\begin{equation}
\texte^{x^2/s}=1+\int_0^\infty
\frac{x}{2\sqrt{t}}\,\texte^{t}\,I_1(x\sqrt{t})\,\texte^{-(1+s/4)t}\,\textd t,
\end{equation}
where $I_1(z):=\sum_{n\geq0}\frac{1}{n!(n+1)!}(z/2)^{2n+1}$ is a modified Bessel function. Expanding $\texte^{t}$ and $\frac1{\sqrt{t}}I_1(x\sqrt{t})$ into power series in~$t$ and scaling~$t$ by $(1+s/4)$ then readily shows
\begin{equation}
q_{n+1} = n!\sum_{j=0}^n\frac{(\alpha^2\theta/8)^{j+1}}{j!(j+1)!(n-j)!}
\end{equation}
for each~$n\ge0$. See also \eqref{E:tilde-mu} for the corresponding formulas in continuous time.

Third, as we will see in the proofs, the random variable~$Y$ in the measure \eqref{E:2.21iu} represents the limit of normalized fluctuations of the local time at the boundary vertex for the first $\lfloor t_N\deg(D_N)\rfloor$ steps of the random walk (see Lemma~\ref{lemma-TY}). 
A key point is that this becomes statistically independent of the level-set statistics in the limit. 
Incidentally, through \eqref{E:2.26}, the total mass of the measure \eqref{E:2.21iu} describes the limit law of a normalized total number of uncovered vertices at the time proportional to~$\lambda^2$-multiple of the cover time.

Fourth, the reader may wonder why we had to include the degree of~$\varrho$ into the normalization of the local time \eqref{E:LVt} by~$\deg(V)$. This is because, although $\deg(\varrho)=o(|D_N|)$ under \twoeqref{E:1.8i}{E:1.8ii} (see Lemma~\ref{lemma-5.8}), once the ratio of $\deg(\varrho)/|D_N|$ is larger than $1/\log N$ (which can occur under \twoeqref{E:1.8i}{E:1.8ii}) removing $\deg(\varrho)$ from the normalization changes the scaling of the normalization constants~$W_N$ and~$\wh W_N$ with~$N$.

Fifth, as in \cite{AB}, the above statements deliberately avoid various boundary values of the parameters; i.e., $\lambda=1$ for the thick points, $\lambda=\sqrt\theta\wedge 1$ for the thin points and~$\theta=1$ for the light and avoided points. All of these are closely related to the statistics of nearly-maximal DGFF values, which is different than the regime described in Theorem~\ref{thm-DGFF}. While the nearly-maximal DGFF values are now well understood thanks to the work of the second author with O.~Louidor \cite{BL1,BL2,BL3} and with S.~Guffler and O.~Louidor~\cite{BGL}, the recent work of Cortines, Louidor and Saglietti~\cite{CLS} shows that the connection between the avoided points at~$\theta=1$ (i.e., the time scale of the cover time) and the DGFF extrema is considerably more subtle.

Sixth, a natural setting for the above problem is the random walk on a lattice torus $(\Z/(N\Z))^2$ started from any given vertex~$\varrho$. As our work in progress shows~\cite{ABL}, the scaling of the corresponding measures is then more complicated --- and, in particular, the scaling sequences~$W_N$ and~$\wh W_N$ have to be taken \myemph{random}. This is related to the fact that, for random walks of time-length order $N^2(\log N)^2$, the local time at the starting point of the walk exhibits fluctuations of order $(\log N)^{3/2}$ on the torus while these are only of order $\log N$ at the boundary vertex in our planar domains. 

 Seventh, we note the recent preprints of Jego~\cite{Jego1,Jego2}, where measures of the kind \eqref{E:zetaND} associated with the thick points of planar Brownian motion run until the first exit from a bounded domain are shown to admit a non-trivial scaling limit that is identified with the limit of multiplicative chaos measures associated with the root of the local time. In~\cite{Jego2} the limit measure is shown to obey a list of natural properties that characterize it uniquely.  It remains to be seen whether the limit measure bears  any connection to Gaussian Free Field and/or Liouville Quantum Gravity. 
 
Finally, we note that Dembo, Peres, Rosen and Zeitouni \cite{DPRZ01, DPRZ06}
and Okada \cite{Okada1, Okada2, Okada3} analyzed the fractal nature and clustering of the sets of thick points
and avoided points in the setting of a random walk killed on exit from $D_N$ (for the thick points) and on two-dimensional torus (for the avoided points). 
In particular, for $0 < \beta < 1$, the growth exponents have been obtained for 
\begin{equation}
 \#\Bigl\{(x_1, x_2) \in D_N\times D_N \colon |x_1 - x_2| \le N^{\beta},\,\min\{L^{D_N}_{H_\varrho}(x_1),L^{D_N}_{H_\varrho}(x_2)\}\ge s(\log N)^2\Bigr\}
\end{equation} 
\MB with $s>0$ \eMB and
\begin{equation}
\#\Bigl\{(x_1, x_2) \in D_N\times D_N\colon |x_1 - x_2| \le N^{\beta},\,\MB\max\eMB\{L^{D_N}_{t_N}(x_1),L^{D_N}_{t_N}(x_2)\}=0\Bigr\},
\end{equation} 
\MB as well as the sets where ``$\min$'' and ``$\max$'' are swapped  --- which amounts to changing from the behavior near a typical point in the level set to a typical point in~$D_N$. \eMB
 These conclusions cannot be gleaned from our results because $N^{-1+\beta}$ vanishes as $N \to \infty$. Notwithstanding, the obtained exponents coincide with those for the DGFF thick points computed by Daviaud~\cite{Daviaud} and thus affirm the universality of the DGFF.

\subsection{Outline}
The rest of this paper is organized as follows. In Section~\ref{sec3} we derive the scaling limit for the level sets of zero-average DGFF. Section~\ref{sec4} extends the conclusions of~\cite{AB} on the local time parametrized by the local time at~$\varrho$ to include information on fluctuations of the total time of the walk. This naturally feeds into Section~\ref{sec5} where we establish the scaling limit of exceptional points for the local time of the continuous-time random walk in the parametrization of the total time. Section~\ref{sec6} then controls the effect of starting the walk at an arbitrary point. In Section~\ref{sec7} we then prove our main theorems above concerning the discrete-time walk except for the local behavior, which is deferred to Section~\ref{sec8}.

\section{Zero average DGFF level sets}
\label{sec3}\noindent
We are now ready to commence the proofs. As our first item of business, we will address Theorem~\ref{thm-DGFF} on the level sets of the zero-average DGFF. Our strategy is to derive the statement from the unconditional convergence \eqref{E:1.19}. This leads to a convolution identity whose resolution requires a uniqueness statement that pertains to the whole class of Gaussian Multiplicative Chaos measures:

\begin{theorem}
\label{lemma-8.1new}
Given a bounded open set~$D\subset\R^d$, let $M^D$ and~$\wt M^D$ be two random a.s.-finite Borel measures on~$D$ and let~$\Phi$ be a centered Gaussian field on~$D$ independent of~$M^D$ and~$\wt M^D$ such that, for some bounded measurable functions $\mathfrak h_k\colon D\to\R$, 
\begin{equation}
\label{E:3.1}
\Cov\bigl(\Phi(x),\Phi(y)\bigr)=\sum_{k=0}^\infty \mathfrak h_k(x)\mathfrak h_k(y),\quad \text{\rm locally uniformly in }x,y\in D.
\end{equation}
Then
\begin{equation}
\label{E:8.1new}
\texte^{\Phi(x)}M^D(\textd x)\,\,\laweq\,\,\texte^{\Phi(x)}\wt M^D(\textd x)
\end{equation}
implies $M^D\laweq \wt M^D$.
\end{theorem}

We remark that for the needs of the present paper it would suffice to treat the case when the sum in \eqref{E:3.1} consists of only one non-zero term. However, this still constitutes the bulk of the proof and so we include the more general case as it is interesting in its own right. The result extends (with suitable modifications) even to the case when~$\Phi$ is a generalized Gaussian Field; the statement thus ``reverse engineers'' the base measure from the associated Gaussian Multiplicative Chaos. Our setting goes even somewhat beyond that of, e.g., Shamov~\cite{Shamov} as we make no moment assumptions on~$M^D$ and~$\wt M^D$.

\smallskip
The proof of Theorem~\ref{lemma-8.1new} hinges  on the following technical observation:

\begin{lemma}
\label{lemma-3.2}
Let~$\mathfrak h\colon D\to\R$ and $f\colon D\to[0,\infty)$ be bounded and measurable and let $M^D$ be a random a.s.-finite Borel measure on~$D$. Let~$Y=\NN(0,1)$ be independent of~$M^D$. Define
$\phi\colon \R\times [0,\infty)\to[0,1]$ by
\begin{equation}
\label{E:3.3}
\phi(\lambda,t):=E\bigl(\texte^{-\langle M^D,\,\texte^{\sqrt{t}\,\mathfrak h(\cdot)Y-\lambda\mathfrak h(\cdot)}f\rangle}\bigr).
\end{equation}
Then~$\phi$ is continuous on its domain and smooth on the interior thereof. Moreover, $\phi$ satisfies the heat equation,
\begin{equation}
\label{E:3.4}
\frac{\partial\phi}{\partial t} = \frac12\frac{\partial^2\phi}{\partial \lambda^2},\qquad (\lambda,t)\in\R\times(0,\infty).
\end{equation}
\end{lemma}

\begin{proofsect}{Proof}
The continuity of~$\phi$ on $\R\times[0,\infty)$ follows by the Bounded Convergence Theorem. Using that $\sqrt{t}Y=\NN(0,t)$ and invoking Tonelli's Theorem we get
\begin{equation}
\phi(\lambda,t) = \int\frac{\textd y}{\sqrt{2\pi t}}\,\texte^{-\frac{(y-\lambda)^2}{2t}}\phi(y,0).
\end{equation}
As~$y\mapsto\phi(y,0)$ is bounded, $\phi$ is continuously differentiable on $\R\times(0,t)$. Since the density of $\NN(0,t)$ solves the heat equation \eqref{E:3.4}, the Dominated Convergence Theorem ensures that so does~$\phi$.
\end{proofsect}

We are now ready to give:

\begin{proofsect}{Proof of Theorem~\ref{lemma-8.1new}}
Let us first assume that~$\Phi$ takes the form~$\mathfrak h(x)Y$ for some bounded measurable $\mathfrak h\colon D\to\R$ and~$Y=\NN(0,1)$ independent of~$M^D$ and~$\wt M^D$. Assume that
\begin{equation}
\label{E:3.10}
\texte^{\mathfrak h(x)Y}M^D(\textd x)\,\,\laweq\,\,\texte^{\mathfrak h(x)Y}\wt M^D(\textd x).
\end{equation}
Given any  bounded and measurable~$f\colon D\to[0,\infty)$, let~$\phi(\lambda,t)$, resp., $\wt\phi(\lambda,t)$ denote the functions in \eqref{E:3.3} with the random measure~$M^D$, resp.,~$\wt M^D$. Since also $x\mapsto\texte^{-\lambda\mathfrak h(x)}f(x)$ is non-negative and measurable, from \eqref{E:3.10} we then have
\begin{equation}
\label{E:3.10e}
\phi(\lambda,1)=\wt\phi(\lambda,1),\quad \lambda\in\R.
\end{equation}
In light of Lemma~\ref{lemma-3.2}, the difference $\phi-\wt\phi$ is a bounded solution to the heat equation in $\R\times(0,\infty)$ with a continuous extension to $\R\times[0,\infty)$. A key point is that the heat equation is known to exhibit \myemph{backward uniqueness}. More precisely, Seregin and \v Sver\'ak~\cite[Theorem~4.1]{SS} implies that every bounded solution to \eqref{E:3.4} that vanishes at a given positive time vanishes everywhere. Since \eqref{E:3.10e} implies that~$\phi-\wt\phi$ vanishes at ``time'' $t=1$, we have $\phi=\wt\phi$ on~$\R\times[0,\infty)$. From the equality $\phi(0,0)=\wt\phi(0,0)$ we then infer
\begin{equation}
E\bigl(\texte^{-\langle M^D,f\rangle}\bigr)=E\bigl(\texte^{-\langle \wt M^D,f\rangle}\bigr).
\end{equation}
Since~$f$ was arbitrary, the claim thus holds for any $\Phi$ of the form $\mathfrak h(\cdot)Y$.

To address the general case, we proceed as in Kahane~\cite{Kahane} (see \cite[Section~5.2]{B-notes} for a review). First note that by \eqref{E:3.1} we may write
\begin{equation}
\Phi(x)\,\laweq\,\Phi_n(x)+\sum_{k=0}^n \mathfrak h_k(x)Y_k,
\end{equation}
where $(Y_0,\dots,Y_n)$ are i.i.d.\ standard normal and where $\Phi_n$ is an independent centered Gaussian field with covariance
\begin{equation}
\Cov\bigl(\Phi_n(x),\Phi_n(y)\bigr)=\sum_{k=n+1}^\infty \mathfrak h_k(x)\mathfrak h_k(y).
\end{equation}
The argument for $\Phi$ of the form $\mathfrak h(\cdot)Y$ then shows, inductively, that \eqref{E:8.1new} implies
\begin{equation}
\label{E:3.19}
\texte^{\Phi_n(x)}M^D(\textd x)\,\,\laweq\,\,\texte^{\Phi_n(x)}\wt M^D(\textd x),\quad n\in\N.
\end{equation}
Letting $f\colon D\to[0,\infty)$ be measurable and supported in a compact set $A\subset D$, the assumption of locally-uniform convergence in \eqref{E:3.1} implies that, given~$\epsilon>0$ there is $n\in\N$ such that $\Var(\Phi_n(x))\le\epsilon$ for all~$x\in A$. This also gives $\Cov(\Phi_n(x),\Phi_n(y))\le\epsilon$ for all~$x,y\in A$ and so Kahane's convexity inequality along with Jensen's inequality show, for $Y_\epsilon=\NN(0,\epsilon)$ independent of~$M^D$ and~$\wt M^D$,
\begin{equation}
\begin{aligned}
E\bigl(\texte^{-\texte^{Y_\epsilon}\langle M^D,f\rangle}\bigr)
&=E\bigl(\texte^{-\texte^{\epsilon/2}
\texte^{Y_\epsilon-\epsilon/2}\langle M^D,f\rangle}\bigr)
\\
&\!\!\!\!\overset{\text{Kahane}}\ge
E\bigl(\texte^{-\texte^{\epsilon/2}
\langle M^D,\,\texte^{\Phi_n(\cdot)-\frac12\Var(\Phi_n(\cdot))}f\rangle}\bigr)
\\
&\!\!\overset{\eqref{E:3.19}}=E\bigl(\texte^{-\texte^{\epsilon/2}
\langle \wt M^D,\,\texte^{\Phi_n(\cdot)-\frac12\Var(\Phi_n(\cdot))}f\rangle}\bigr)
\overset{\text{Jensen}}\ge E\bigl(\texte^{-\texte^{\epsilon/2}
\langle \wt M^D,f\rangle}\bigr).
\end{aligned}
\end{equation}
Taking~$\epsilon\downarrow0$ and noting that this implies~$Y_\epsilon\to0$ in probability then shows, with the help of the Bounded Convergence Theorem,
\begin{equation}
E\bigl(\texte^{-\langle M^D,f\rangle}\bigr)\ge E\bigl(\texte^{-\langle\wt M^D,f\rangle}\bigr).
\end{equation}
By symmetry, equality must hold for all~$f$ as above and so $M^D\laweq \wt M^D$, as desired.
\end{proofsect}

Equipped with Theorem~\ref{lemma-8.1new}, we are ready to give:

\begin{proofsect}{Proof of Theorem~\ref{thm-DGFF}}
Abbreviate
\begin{equation}
\label{E:3.21}
Y_N:=\frac1{|D_N|}\sum_{x\in D_N}h^{D_N}_x.
\end{equation}
Then $Y_N$ is normal with mean zero and variance
\begin{equation}
\Var(Y_N)=\frac1{|D_N|^2}\sum_{x,y\in D_N}G^{D_N}(x,y).
\end{equation}
Moreover, denoting
\begin{equation}
\frakd_N(x):=\frac{|D_N|\sum_{y\in D_N}G^{D_N}(\lfloor xN\rfloor,y)}{\sum_{y,z\in D_N}G^{D_N}(z,y)}
\end{equation}
a covariance calculation shows that~$Y_N$ is independent of
\begin{equation}
\wh h^{D_N}_x:=h^{D_N}_x-\frakd_N(x/N)Y_N
\end{equation}
which, we note, has zero average over~$D_N$. Hence, if we define the zero-average variant of~$\eta^D_N$ by
\begin{equation}
\label{E:3.25}
\eta^{D,0}_N:=\frac1{K_N}\sum_{x\in D_N}\delta_{x/N}\otimes\delta_{\,\wh h^{D_N}_x- \wh{a}_N},
\end{equation}
we have
\begin{equation}
\label{E:3.26}
\eta^{D,0}_N\independent Y_N\quad\text{and}\quad \eta^D_N=\eta^{D,0}_N\circ\theta_{\frakd_N(\cdot)Y_N}^{-1},
\end{equation}
where~$\theta_{s(\cdot)}\colon D\times\R\to D\times\R$ is defined by $\theta_{s(\cdot)}(x,h):=(x,h+s(x))$. The stated independence also shows
\begin{equation}
\Bigl(\eta^D_N\,\Big|\,\sum_{x\in D_N}h^{D_N}_x=0\Bigr) \,\laweq\,\,\eta^{D,0}_N
\end{equation}
and so we may and will henceforth focus on the limit of~$\eta^{D,0}_N$.

Using the uniform bound $G^{D_N}(x,y)\le g\log\frac{N}{|x-y|+1}+c$ along with
\begin{equation}
\label{E:3.28}
G^{D_N}\bigl(\lfloor xN\rfloor,\lfloor yN\rfloor\bigr)\underset{N\to\infty}\longrightarrow\,\wh G^D(x,y),\quad x,y\in D,\,x\ne y,
\end{equation}
the Dominated Convergence shows that~$\Var(Y_N)$ converges to $\sigma_D^2$ from \eqref{E:2.14new}. We thus have $Y_N\Lawarrow Y=\NN(0,\sigma_D^2)$. In particular, $\{Y_N\colon N\ge1\}$ is tight and so from the tightness of~$\eta^D_N$, \eqref{E:3.26} and the uniform boundedness of~$\frakd_N$ we get
\begin{equation}
\{\eta^{D,0}_N\colon N\ge1\}\text{ is tight}.
\end{equation}
Similarly we show that $\frakd_N\to\frakd$ uniformly on~$D$. (This implies~$\frakd(x)\ge0$). Writing the equality in \eqref{E:3.26} via Laplace transforms against a test function $f\in C_\cc(D\times\R)$ and invoking \eqref{E:1.19}, any subsequential limit $\eta^{D,0}$ of $\{\eta^{D,0}_N\colon N\ge1\}$ thus obeys
\begin{equation}
\label{E:8.14new}
\eta^{D,0}\circ\theta_{\frakd(\cdot)Y}^{-1} \,\laweq\, \cspecial(\lambda)Z^D_\lambda(\textd x)\otimes\texte^{-\alpha\lambda h}\textd h\,,
\end{equation}
where~$Y=\NN(0,\sigma_D^2)$ is such that $Y\independent \eta^{D,0}$ on the left-hand side. 

Next we note that we may realize \eqref{E:8.14new} as an a.s.\ equality. This is because \eqref{E:8.14new} implies, for any measurable $A\subseteq D$ and $B\subseteq\R$ with $\leb(A)>0$,
\begin{equation}
\frac{\eta^{D,0}\circ\theta_{\frakd(\cdot)Y}^{-1}(A\times B)}{\eta^{D,0}\circ\theta_{\frakd(\cdot)Y}^{-1}(A\times[0,1])}=\alpha\lambda(1-\texte^{-\alpha\lambda})^{-1}\int_B\texte^{-\alpha\lambda h}\,\textd h\quad\text{a.s.}
\end{equation}
due to the fact that equality in law to a constant implies equality a.e. We conclude that the measure
\begin{equation}
A\mapsto\alpha\lambda[\cspecial(\lambda)(1-\texte^{-\alpha\lambda})]^{-1}\,\eta^{D,0}\circ\theta_{\frakd(\cdot)Y}^{-1}(A\times[0,1])
\end{equation}
is equidistributed to~$Z^D_\lambda$.
Replacing $Z^D_\lambda$ by this measure then gives us equality a.s.

Once we have \eqref{E:8.14new} as an a.s.\ equality, and $Z^D_\lambda$ thus as a measurable function of~$\eta^{D,0}$ and~$Y$, we apply a routine change of variables to get
\begin{equation}
\label{E:3.33}
\eta^{D,0}=\cspecial(\lambda)\,\texte^{-\alpha\lambda\frakd(x)Y}\,Z^D_\lambda(\textd x)\otimes\texte^{-\alpha\lambda h}\textd h.
\end{equation}
Setting
\begin{equation}
\label{E:8.19new}
Z^{D,0}_\lambda(\textd x):=\texte^{-  \alpha\lambda\frakd(x)Y}Z^D_\lambda(\textd x)
\end{equation}
 the independence of~$\eta^{D,0}$ of~$Y$ shows $Z^{D,0}_\lambda\independent Y$ and thus proves existence of the decomposition \eqref{E:2.15new}. Since the decomposition is unique by Theorem~\ref{lemma-8.1new} and the fact that~$\frakd$ is bounded and continuous, the law of $Z^{D,0}_\lambda$ does not depend on the subsequential limit~$\eta^{D,0}$. It follows that all subsequential limits of $\{\eta^{D,0}_N\colon N\ge1\}$ are equal in law and so we get the convergence statement \eqref{E:2.16new} as well.
\end{proofsect}

Our use of Theorem~\ref{thm-DGFF} will invariably come through:

\begin{corollary}
\label{cor-3.3}
Under the conditions of Theorem~\ref{thm-DGFF}, and for~$Y_N$ as in \eqref{E:3.21},
\begin{equation}
\eta^D_N\otimes\delta_{Y_N}\,\,\,\underset{N\to\infty}\Lawarrow\,\,\,\cspecial(\lambda)\,\texte^{\alpha\lambda\frakd(x)Y}\,Z^{D,0}_\lambda(\textd x)\otimes\texte^{-\alpha\lambda h}\textd h\otimes\delta_Y,
\end{equation}
where $Y=\NN(0,\sigma_D^2)$, for~$\sigma_D^2$ as in \eqref{E:2.14new}, is such that $Y\independent Z^{D,0}_\lambda$.
\end{corollary}

\begin{proofsect}{Proof}
By \eqref{E:3.26} and the fact that $Y_n\to Y$ in law and $\frakd_N\to\frakd$ uniformly shows
\begin{equation}
\eta_N^D\otimes\delta_{Y_N}\,\,\underset{N\to\infty}\Lawarrow\,\,\,\bigl(\eta^{D,0}\circ\theta^{-1}_{\frakd(\cdot)Y}\bigr)\otimes\delta_Y,
\end{equation}
where $\eta^{D,0}$ is as in \eqref{E:3.33} and obeys $Y\independent\eta^{D,0}$. Invoking \eqref{E:8.19new}, the claim  follows by a routine change of variables.
\end{proofsect}

\section{Augmented boundary vertex measures}
\label{sec4}\noindent
We will now move to the discussion of local time level sets. Our proofs build on the conclusions derived in~\cite{AB} for the local time parametrized by its value at the boundary vertex~$\varrho$. In order to transfer these conclusions to the setting of a fixed total time, we will need to control the fluctuations of the total local time at a fixed local time at~$\varrho$. Our  first step is thus to augment the results of~\cite{AB} by information about these fluctuations.
 
We will again introduce the corresponding quantities on a general finite connected graph with vertex set~$V\cup\{\varrho\}$. Consider a joint law of paths~$X$ of the discrete-time random walk on~$V\cup\{\varrho\}$ and an independent sample $t\mapsto \wt N(t)$ of a rate-1 Poisson process. The continuous-time walk is then defined as
\begin{equation}
\label{E:tilde-X}
\wt X_t:=X_{\wt N(t)},\quad t\ge0.
\end{equation}
The local time naturally associated with~$\wt X$ is given by
\begin{equation}
\label{E:2.32uw}
\wt L_t^V(u):=\frac1{\deg(u)}\int_0^t\textd s\,1_{\{\wt X_s=u\}}.
\end{equation}
Denoting $\hat\tau_{\varrho}(t):=\inf\{s\ge0\colon\wt L_s^V(\varrho)\ge t\}$, the local time parametrized by its value at~$\varrho$ is defined as
\begin{equation}
\label{E:2.27ia}
\wh L_t^V(v):=\wt L_{\hat\tau_{\varrho}(t)}^V(v).
\end{equation}
Note that, in particular, we have $\wh L_t^V(\varrho)=t$ for all~$t\ge0$. The same is true about the expected value at any vertex; i.e., $E^\varrho\wh L^V_t(v)=t$ for all~$v\in V$.

At a given~$t\ge0$, the total (continuous) local time of the walk is computed by adding $\wh L^V_t(v)$ over all~$v\in V\cup\{\varrho\}$. The quantity
\begin{equation}
\label{E:4.4u}
T(t):=\frac1{\sqrt{2t}\,|V|}\sum_{v\in V}\bigl[\,\wh L^V_t(v)-t\bigr]
\end{equation}
then denotes the normalized (empirical) fluctuation of the total local time. (Note that $v=\varrho$ can be freely added to the sum as $\wh L^V_t(\varrho)=t$.) To explain the specific choice of the normalization, we recall the following result from Eisenbaum, Kaspi, Marcus, Rosen and Shi~\cite{EKMRS} (with improvements by Zhai~\cite[Section~5.4]{Zhai}):

\begin{theorem}[Second Ray-Knight Theorem]
\label{lemma-Dynkin}
For each~$t>0$ there exists a coupling of $\wh L^V_t$ (sampled under~$P^\varrho$) and two copies of the DGFF $h^V$ and~$\tilde h^V$ such that
\begin{equation}
\label{E:Dynkin1}
h^V\text{\rm\ and } \wh L^V_t \text{\rm\ are independent}
\end{equation}
and
\begin{equation}
\label{E:Dynkin2}
\wh L_t^V(u)+\frac12(h_u^V)^2 = \frac12\bigl(\tilde h_u^V+\sqrt{2t}\bigr)^2, \quad u\in V.
\end{equation}
\end{theorem}

Using the stated coupling, we readily compute
\begin{equation}
\label{E:4.7}
T(t) = \frac1{|V|}\sum_{u\in V}\tilde h^V_u 
+\frac1{\sqrt{2t}\,|V|}\sum_{u\in V}\frac{(\tilde h^V_u)^2-(h^V_u)^2}2.
\end{equation}
Note that the first term is the average of the field~$\tilde h^V$.

In what follows, the role of~$V$ will be taken by the sets~$D_N$ and~$\varrho$ by the ``boundary vertex.''  We let $h^{D_N}$ be the DGFF on~$D_N$ and, given a sequence $\{t_N\}_{N\ge1}$ and for the continuous-time random walk started at~$\varrho$, let~$\tilde h^{D_N}$ be the DGFF such that \twoeqref{E:Dynkin1}{E:Dynkin2} with $t:=t_N$ holds. We then set
\begin{equation}
\label{E:4.8}
T_N:=\frac1{\sqrt{2t_N}\,|D_N|}
\sum_{x\in D_N}\bigl[\,\wh L^{D_N}_{t_N}(x)-t_N\bigr]
\end{equation}
and denote
\begin{equation}
\label{E:4.9nwt}
Y_N:=\frac1{|D_N|}\sum_{x\in D_N}\tilde h^{D_N}_x.
\end{equation}
We start by noting:

\begin{lemma}
\label{lemma-TY}
For any~$\{t_N\}_{N\ge1}$ with~$t_N\to\infty$ we have 
\begin{equation}
\label{E:4.10}
T_N-Y_N\,\,\underset{N\to\infty}\longrightarrow\,\,0,\quad\text{\rm in probability}.
\end{equation}
In particular,
\begin{equation}
\label{E:4.11}
T_N\,\,\,\underset{N\to\infty}\Lawarrow\,\,\,\NN(0,\sigma_D^2),
\end{equation}
where~$\sigma_D^2$ is as in \eqref{E:2.14new}.
\end{lemma}

\begin{proofsect}{Proof}
The Wick Pairing Theorem gives
\begin{equation}
\begin{aligned}
\Var\Bigl(\,\sum_{x\in D_N} (h^{D_N}_x)^2\Bigr)
&=\sum_{x,y\in D_N}\Cov\bigl((h^{D_N}_x)^2,(h^{D_N}_y)^2\bigr)
\\
&=\sum_{x,y\in D_N} 2\,\bigl[E(h^{D_N}_x h^{D_N}_y)\bigr]^2
=2\sum_{x,y\in D_N}G^{D_N}(x,y)^2.
\end{aligned}
\end{equation}
The uniform bound $G^{D_N}(x,y)\le g\log\frac{N}{|x-y|+1}+c$ shows that the double sum on the right is of order~$|D_N|^2$. From~$t_N\to\infty$ it follows that
\begin{equation}
\frac1{\sqrt{2t_N}\,|D_N|}\sum_{x\in D_N}\bigl[(h^{D_N}_x)^2-\E[(h^{D_N}_x)^2]\bigr]\,\,\underset{N\to\infty}\longrightarrow\,\,0,\quad\text{in probability}.
\end{equation}
Using this along with $\E[(h^{D_N}_x)^2]=\E[(\tilde h^{D_N}_x)^2]$ in \eqref{E:4.7}, we get \eqref{E:4.10}. For \eqref{E:4.11} we invoke the argument after \eqref{E:3.28}.
\end{proofsect}

We are now ready to state and prove convergence theorems for processes associated with exceptional level sets of the boundary vertex local time~$\wh L^{D_N}_{t_N}$ augmented by information about~$T_N$. Starting with the thick and thin points, given positive sequences $\{t_N\}_{N\ge1}$ and $\{a_N\}_{N\ge1}$, define
\begin{equation}
\label{E:zetaND-bv}
\wh\zeta^D_N:=\frac1{W_N}\sum_{x\in D_N}\delta_{x/N}\otimes\delta_{(\wh L_{t_N}^{D_N}(x)-a_N)/\sqrt{2a_N}}\,,
\end{equation}
where~$W_N$ is as in \eqref{E:WN}. For the thick points of~$\wh L^{D_N}_{t_N}$, we then have:

\begin{proposition}[Thick points]
\label{thm-4.3}
Suppose that $\{t_N\}_{N\ge1}$ and $\{a_N\}_{N\ge1}$ are such \eqref{E:1.12} and \eqref{E:1.20} hold for some $\theta>0$ and~$\lambda\in(0,1)$. Then for~$X$ sampled from~$P^\varrho$, relative to the vague convergence of measures on~$\overline D\times(\R\cup\{+\infty\})\times\R$,
\begin{equation}
\label{E:1.21cont}
\wh\zeta^D_N\otimes\delta_{T_N}\,\,\,\underset{N\to\infty}\Lawarrow\,\,\, \sqrt{\frac{\sqrt\theta}{\sqrt\theta+\lambda}}\,\,\cspecial(\lambda) \,\texte^{\alpha\lambda\frakd(x)Y}Z_\lambda^{D,0}(\textd x)\otimes\texte^{-\alpha\lambda h}\textd h\otimes\delta_Y(\textd t)
\end{equation}
where~$Y=\NN(0,\sigma_D^2)$, for~$\sigma_D^2$ as in \eqref{E:2.14new}, is such that~$Y\independent Z^{D,0}_\lambda$.
\end{proposition}

\begin{proofsect}{Proof}
We will rely heavily on the proof of \cite[Theorem~2.2]{AB} but, due to a different normalization of the second coordinate in \eqref{E:zetaND-bv} and also the fact that the limit measure is different than in \cite{AB}, we need to recount the main steps of the proof. Throughout we will assume (for each~$N\ge1$ and each~$t:=t_N$) a coupling of~$\wh L^{D_N}_{t_N}$ and an independent DGFF~$h^{D_N}$ to a DGFF~$\tilde h^{D_N}$ satisfying \eqref{E:Dynkin2}.

First, by~\cite[Corollary~4.2]{AB} the measures $\{\wh\zeta^D_N\colon N\ge1\}$ are tight and, by Lemma~\ref{lemma-TY}, the same applies to the enhanced measures $\{\xi_N\colon N\ge1\}$ where
\begin{equation}
\label{E:4.16}
\xi_N:=\wh\zeta^D_N\otimes\delta_{T_N}.
\end{equation}
Moreover, \cite[Lemma~5.3]{AB} shows that if~$\xi_{N_k}\to\xi$ in law along some increasing sequence~$\{N_k\}_{k\ge1}$, then the extended measures
\begin{equation}
\label{E:xi}
\xi_N^{\text{ext}}:=\frac1{W_N}\sum_{x\in D_N}\delta_{x/N}\otimes\delta_{(\wh L_{t_N}^{D_N}(x)-a_N)/\sqrt{2a_N}}\otimes\delta_{T_N}\otimes\delta_{h^{D_N}_x/(2a_N)^{1/4}}\,,
\end{equation}
where we now normalize the third coordinate differently than in \cite{AB}, obey
\begin{equation}
\label{E:4.18}
\xi_{N_k}^{\text{ext}}\,\,\,\underset{ k \to \infty}\Lawarrow\,\,\,\xi\otimes\mathfrak g
\end{equation}
in which, using \eqref{E:1.20}, $\mathfrak g$ is the law of~$\NN(0,\frac1{\alpha\,(\sqrt\theta+\lambda)})$.

Let~$\eta^D_N$ be the process \eqref{E:etaDGFF} associated with the field~$\tilde h^{D_N}$ and the scale function\begin{equation}
\wh a_N:=\sqrt{2a_N}-\sqrt{2t_N}
\end{equation}
that, by \eqref{E:1.12} and \eqref{E:1.20}, scales as 
 $\wh a_N\sim2\sqrt g\,\lambda\log N$ as~$N\to\infty$.
Let~$Y_N$ be the average of~$\tilde h^{D_N}$ over~$D_N$;  cf~\eqref{E:4.9nwt}. 
Given $f\in C_\cc(D\times\R\times\R)$, in the assumed coupling of~$\wh L^{D_N}_{t_N}$, $h^{D_N}$ and~$\tilde h^{D_N}$, the convergence in Lemma~\ref{lemma-TY} tells us
\begin{equation}
\label{E:4.20}
\langle\eta^D_N\otimes\delta_{Y_N},f\rangle=o(1)+\langle\eta^D_N\otimes\delta_{T_N},f\rangle,
\end{equation}
where~$o(1)\to0$ as~$N\to\infty$ in probability. The calculation in the proof of \cite[Lemma~5.4]{AB} (enabled by the fact that the field~$h^{D_N}$ will be typical at most points contributing to~$\zeta^D_N$, as shown in~\cite[Lemma~5.2]{AB}) then gives
\begin{equation}
\label{E:4.21}
\langle\eta^D_N\otimes\delta_{T_N},f\rangle = o(1)+\langle\xi_N^{\text{ext}},f^{\text{ext}}\rangle,
\end{equation}
where
\begin{equation}
f^{\text{ext}}(x,\ell,t,h):=f\bigl(x,\ell+\tfrac {h^2}2,t\bigr).
\end{equation}
Using Corollary~\ref{cor-3.3} on the left-hand side of \eqref{E:4.20}, from \eqref{E:4.21} and \eqref{E:4.18} and, one more time, \cite[Lemma~5.2]{AB} we conclude that every subsequential limit~$\xi$ of the measures in \eqref{E:4.16}
satisfies the convolution-type identity
\begin{equation}
\label{E:4.23}
\langle\xi,f^{\ast\mathfrak g}\rangle \,\laweq\, \cspecial(\lambda)\int\texte^{\alpha\lambda\frakd(x)Y}\,Z^{D,0}_\lambda(\textd x)\otimes\texte^{-\alpha\lambda \ell}\textd \ell \,\,f(x,\ell,Y),
\end{equation}
where $Y\independent Z^{D,0}_\lambda$ and
\begin{equation}
\label{E:4.24}
f^{\ast\mathfrak g}(x,\ell,t):=\int\mathfrak g(\textd h)f\bigl(x,\ell+\tfrac {h^2}2,t\bigr),
\end{equation}
jointly for all $f\in C_\cc(D\times\R\times\R)$. It remains to ``solve'' \eqref{E:4.23} for~$\xi$.

First we note that the Monotone Convergence Theorem extends \eqref{E:4.23} to all~$f$ of the form $f(x,\ell,t):=1_A(x)\tilde f(\ell)1_{(b,\infty)}(t)$, where~$\tilde f\in C_\cc(\R)$ and where~$A\subseteq D$ is non-empty and open. Denoting $\xi_{A,b}(B):=\xi(A\times B\times(b,\infty))$, a calculation then shows
\begin{equation}
\label{E:4.25}
\langle\xi,f^{\ast\mathfrak g}\rangle = \langle\xi_{A,b},\tilde f\ast\frake\rangle
\end{equation}
where
\begin{equation}
 \frake(z):=\sqrt{\frac\beta\pi}\,\frac{\texte^{\beta z}}{\sqrt{-z}}1_{(-\infty,0)}(z)\quad\text{for}\quad\beta:=\alpha\bigl(\sqrt\theta+\lambda\bigr).
\end{equation}
The identity \eqref{E:4.23} also implies that $\langle\xi_{A,b},1_{[0,\infty)}\rangle<\infty$ a.s.\ and gives
\begin{equation}
\label{E:4.27}
\langle\xi_{A,b},\tilde f\ast\frake\rangle =  \langle\xi_{A,b},1_{[0,\infty)} \ast \frake\rangle \int 
\alpha\lambda\,\texte^{-\alpha\lambda\ell} \tilde f (\ell) \textd \ell,
\end{equation}
where the equality now holds pointwise a.s.\ because once 
$\langle\xi_{A,b},1_{[0,\infty)} \ast \frake \rangle>0$ (which is necessary for the left-hand side to be non-zero), the ratio $\langle\xi_{A,b},\tilde f\ast\frake\rangle/ \langle\xi_{A,b},1_{[0,\infty)} \ast \frake\rangle$ is equal in law, and thus pointwise, to the integral on the right.

Denoting $\mu_\lambda(\textd h):=\texte^{-\alpha\lambda h}\textd h$, a routine change of variables  rewrites \eqref{E:4.27} as
\begin{equation}
\label{E:4.28}
\langle\xi_{A,b},\tilde f\ast\frake\rangle =C\langle\mu_\lambda,\tilde f\rangle
\end{equation}
where~$C$ is a random constant that is finite thanks to~$\beta>\alpha\lambda$. By \cite[Lemma~5.5]{AB}, there is at most one Borel measure~$\xi_{A,b}$ on~$\R$ satisfying \eqref{E:4.28} and, in fact, $\xi_{A,b}(\textd\ell)=C_{A,b}\texte^{-\alpha\lambda\ell}\textd\ell$ for some (random) constant~$C_{A,b}$. It follows that
\begin{equation}
\xi(\textd x\textd\ell\textd t)=M(\textd x\textd t)\otimes \texte^{-\alpha\lambda\ell}\textd\ell,
\end{equation}
where, by plugging this in \eqref{E:4.23},
\begin{equation}
M(\textd x\textd t)\,\laweq\,\Bigl(\int\mathfrak g(\textd h)\texte^{\alpha\lambda\frac{h^2}2}\Bigr)^{-1}
\cspecial(\lambda)\texte^{\alpha\lambda\frakd(x)Y}\,Z^{D,0}_\lambda(\textd x)\otimes\delta_Y(\textd t).
\end{equation}
The integral equals the root of $(\sqrt\theta+\lambda)/\sqrt\theta$. The claim follows.
\end{proofsect}

We proceed with the corresponding result for the thin points:

\begin{proposition}[Thin points]
\label{thm-4.4}
Suppose that $\{t_N\}_{N\ge1}$ and $\{a_N\}_{N\ge1}$ are such \eqref{E:1.12} and \eqref{E:1.22} hold for some $\theta>0$ and~$\lambda\in(0,\sqrt\theta\wedge1)$. Then for~$X$ sampled from~$P^\varrho$, relative to the vague convergence of measures on~$\overline D\times(\R\cup\{-\infty\})\times\R$,
\begin{equation}
\label{E:1.23cont}
\wh\zeta^D_N\otimes\delta_{T_N}\,\,\,\underset{N\to\infty}\Lawarrow\,\,\, \sqrt{\frac{\sqrt\theta}{\sqrt\theta-\lambda}}\,\,\cspecial(\lambda) \,\texte^{-\alpha\lambda\frakd(x)Y}Z_\lambda^{D,0}(\textd x)\otimes\texte^{+\alpha\lambda h}\textd h\otimes\delta_Y(\textd t)
\end{equation}
where~$Y\independent Z^{D,0}_\lambda$ with $Y=\NN(0,\sigma_D^2)$, for~$\sigma_D^2$ as in \eqref{E:2.14new}.
\end{proposition}

\begin{proofsect}{Proof}
The proof is very similar to that of Proposition~\ref{thm-4.3} so we indicate only the needed changes. We will again rely on the coupling of~$\wh L^{D_N}_{t_N}$ and two DGFFs $h^{D_N}$ and~$\tilde h^{D_N}$ such that \twoeqref{E:Dynkin1}{E:Dynkin2} for~$t:=t_N$ hold. Let~$\eta^D_N$ to denote the process associated with~$\tilde h^{D_N}$ and the centering sequence~$-\wh a_N$, where
\begin{equation}
\wh a_N:=\sqrt{2t_N}-\sqrt{2a_N}.
\end{equation}
Note that, under \eqref{E:1.12} and \eqref{E:1.22} we have $\wh a_N\sim2\sqrt g\lambda\log N$. Writing~$Y_N$ for the average of~$\tilde h^{D_N}$ over~$D_N$, 
Corollary~\ref{cor-3.3} along with the symmetry $h^{D_N}\,\laweq\,-h^{D_N}$ ensures
\begin{equation}
\label{E:4.32}
\eta^D_N\otimes\delta_{Y_N}\,\,\,\underset{N\to\infty}\Lawarrow\,\,\,\cspecial(\lambda)\,\texte^{-\lambda\alpha\frakd(x)Y}\,Z_\lambda^{D,0}(\textd x)\otimes\texte^{+\alpha\lambda h}\textd h\otimes\delta_Y(\textd t),
\end{equation}
where~$Y=\NN(0,\sigma_D^2)$ is independent of~$Z^{D,0}_\lambda$.

The argument now proceeds very much like for the thick points. We consider the extended measures \eqref{E:xi}, which are tight by \cite[Corollary~4.8]{AB} and show, with the help of \cite[Lemmas~6.1, 6.2]{AB} and \eqref{E:4.32}, that every subsequential limit~$\xi$ thereof obeys
\begin{equation}
\label{E:4.33}
\langle\xi,f^{\ast\mathfrak g}\rangle \,\laweq\, \cspecial(\lambda)\int
\texte^{- \alpha\lambda\frakd(x)Y}\,Z^{D,0}_\lambda(\textd x)\otimes
\texte^{+\alpha\lambda \ell}\textd \ell \,\,f(x,\ell,Y),
\end{equation}
where~$f^{\ast\mathfrak g}$ is still defined via \eqref{E:4.24} but with
\begin{equation}
\mathfrak g := \text{law of }\NN\bigl(0, \tfrac1{\alpha(\sqrt\theta-\lambda)}\bigr).
\end{equation}
The identity \eqref{E:4.33} readily extends to all~$f$ of the form $f(x,\ell,t):=1_A(x)\tilde f(\ell)1_{(-\infty,b)}(t)$, where~$\tilde f\in C_\cc(\R)$ and where~$A\subseteq D$ is non-empty and open. A calculation then shows \eqref{E:4.25} with~$\frake$ now defined using $\beta:=\alpha(\sqrt\theta-\lambda)$. Proceeding via an analogue of \eqref{E:4.27} (with $1_{[0,\infty)}$ replaced by $1_{(-\infty,0]}$), using \cite[Lemma~6.4]{AB} we then again show
\begin{equation}
\xi(\textd x\textd\ell\textd t)=M(\textd x\textd t)\otimes \texte^{+\alpha\lambda\ell}\textd\ell,
\end{equation}
where, this time,
\begin{equation}
M(\textd x\textd t)\,\laweq\,\Bigl(\int\mathfrak g(\textd h)\texte^{-\alpha\lambda\frac{h^2}2}\Bigr)^{-1}
\cspecial(\lambda)\texte^{-\alpha\lambda\frakd(x)Y}\,Z^{D,0}_\lambda(\textd x)\otimes\delta_Y(\textd t).
\end{equation}
The integral equals the root of $(\sqrt\theta-\lambda)/\sqrt\theta$.
\end{proofsect}

Next we move to the discussion of the light and avoided points. Starting with the light points, we define
\begin{equation}
\label{E:hat-varthetaND}
\wh\vartheta^D_N:=\frac1{\wh W_N }\sum_{x\in D_N}\delta_{x/N}\otimes\delta_{\wh L_{t_N}^{D_N}(x)},
\end{equation}
where $\wh W_N$ is as in \eqref{E:1.31}. We then get:

\begin{proposition}[Light points]
\label{thm-light-bv}
Suppose $\{t_N\}_{N\ge1}$ obeys \eqref{E:1.12} for some $\theta\in(0,1)$. Then, for the random walk sampled from~$P^\varrho$, in the sense of vague convergence of measures on~$\overline D\times[0,\infty)\times\R$,
\begin{equation}
\wh\vartheta^D_N\otimes\delta_{T_N}\,\,\,\underset{N\to\infty}\Lawarrow\,\,\,  \sqrt{2\pi g}\,\cspecial(\sqrt\theta)\,\,\texte^{-\alpha\sqrt\theta\,\frakd(x)Y}\,  Z_{\sqrt{\theta}\,}^{D,0}(\textd x)\otimes\tilde\mu(\textd h)\otimes\delta_Y(\textd t),
\end{equation}
where $Y=\NN(0,\sigma_D^2)$ is independent of~$Z^{D,0}_{\sqrt\theta}$ and
\begin{equation}
\label{E:tilde-mu}
\tilde\mu(\textd h):=\delta_0(\textd h)+\biggl(\,\sum_{n=0}^\infty\frac1{n!(n+1)!}\Bigl(\frac{\alpha^2\theta}2\Bigr)^{n+1} h^n\biggr)1_{(0,\infty)}(h)\,\textd h.
\end{equation}
\end{proposition}

\begin{proofsect}{Proof}
Assuming again the coupling from \twoeqref{E:Dynkin1}{E:Dynkin2}, we set
\begin{equation}
\xi_N:=\wh\vartheta^D_N\otimes\delta_{T_N}.
\end{equation}
The family~$\{\xi_N\colon N\ge1\}$ is tight by \cite[Corollary~4.6]{AB} and so we may consider a subsequential limit~$\xi$ thereof. By \cite[Lemma~7.1]{AB}, the extended measure
\begin{equation}
\xi_N^{\text{ext}}:=\frac{\sqrt{\log N}}{\wh W_N }\sum_{x\in D_N}\delta_{x/N}\otimes\delta_{\wh L_{t_N}^{D_N}(x)}\otimes\delta_{T_N}\otimes\delta_{h^{D_N}_x},
\end{equation}
then converges to~$\xi\otimes\frac1{\sqrt{2\pi g}}\leb$ along the same subsequence. We now pick a test function $f\in C_\cc(D\times[0,\infty)\times\R)$, denote
\begin{equation}
f^{\text{ext}}(x,\ell,t,h):=f\bigl(x,\ell+\tfrac{h^2}2,t\bigr)
\end{equation}
 and observe that \eqref{E:Dynkin2} implies
\begin{equation}
\sum_{x\in D_N}f^{\text{ext}}\Bigl(\ffrac xN,\wh L^{D_N}_{t_N}(x), T_N,h^{D_N}_x\Bigr)
=\sum_{x\in D_N}f\Bigl(\ffrac xN,\tfrac12\bigl(\tilde h^{D_N}_x+\sqrt{2t_N}\bigr)^2,T_N\Bigr).
\end{equation}
Writing this in terms of the above measures, Lemma~\ref{lemma-TY} gives
\begin{equation}
\langle\xi_N^{\text{ext}},f^{\text{ext}}\rangle = o(1)+\bigl\langle\eta^D_N\otimes\delta_{Y_N},f(~\cdot~, \tfrac{1}{2} |\cdot|^2, ~\cdot~)\bigr\rangle,
\end{equation}
where~$\eta^D_N$ is the DGFF process associated with the scale sequence $\wh a_N:=-\sqrt{2t_N}$. As $\wh a_N\sim-2\sqrt g\sqrt\theta\log N$, from \eqref{E:4.32} we  get
\begin{equation}
\langle \xi,f^{\ast\leb}\rangle \,\laweq\,
\cspecial(\sqrt\theta)\int\texte^{-\alpha\sqrt\theta\,\frakd(x)Y}\,Z^{D,0}_{\sqrt\theta}(\textd x)\otimes\texte^{+\alpha\sqrt\theta\, h}\textd h \,\,f\bigl(x,\tfrac12 h^2,Y\bigr),
\end{equation}
where
\begin{equation}
f^{\ast\leb}(x,\ell,t):=\frac1{\sqrt{2\pi g}}\int\textd h\, f\bigl(x,\ell+\tfrac{h^2}2,t\bigr).
\end{equation}
By the Monotone Convergence Theorem, this extends to all~$f$ of the form
\begin{equation}
f(x,\ell,t):=1_A(x)\texte^{-s\ell}1_{[0,\infty)}(\ell)1_{[b,\infty)}(t)
\end{equation}
for~$A\subseteq D$ open,~$b\in\R$ and~$s>0$. For $\xi_{A,b}(B):=\xi(A\times B\times[b,\infty))$, we then get
\begin{equation}
\int_0^\infty\xi_{A,b}(\textd\ell) \texte^{-s\ell}
\,\laweq\, \sqrt{2\pi g}\, \cspecial(\sqrt\theta)\,
\Bigl(\int_{A}
\texte^{-\alpha\sqrt\theta\,\frakd(x)Y}\,Z^{D,0}_{\sqrt\theta}(\textd x)\Bigr)\,\texte^{\frac{\alpha^2\theta}{2s}}\,1_{[b,\infty)}(Y).
\end{equation}
Since the Laplace transform of a measure, if exists, determines the measure uniquely, this proves that $\xi$ takes the product form
\begin{equation}
\xi \,\laweq\,  \sqrt{2\pi g}\, \cspecial(\sqrt\theta)\,\texte^{-\alpha\sqrt\theta\,\frakd(x)Y}\,Z^{D,0}_{\sqrt\theta}(\textd x)\otimes\tilde\mu(\textd \ell)\otimes\delta_Y(\textd t)
\end{equation}
for some deterministic measure~$\tilde\mu$ on $[0,\infty)$ with Laplace transform $s\mapsto\texte^{\frac{\alpha^2\theta}{2s}}$. A calculation shows that the measure \eqref{E:tilde-mu} has this property.
\end{proofsect}

A direct consequence of our control of the light points is:

\begin{proposition}[Avoided points]
\label{thm-avoid-bv}
Suppose $\{t_N\}_{N\ge1}$ is such that \eqref{E:1.12} holds for some $\theta\in(0,1)$ and let
\begin{equation}
\wh\kappa_N^D:=\frac1{\wh W_N}\sum_{x\in D_N}1_{\{\wh L^{D_N}_{t_N}(x)=0\}}\delta_{x/N}.
\end{equation}
Then, for the random walk distributed according to~$P^\varrho$, in the sense of vague convergence of measures on~$\overline D\times\R$,
\begin{equation}
\label{E:2.22cont}
\wh\kappa^D_N\otimes\delta_{T_N}\,\,\,\underset{N\to\infty}\Lawarrow\,\,\,  \sqrt{2\pi g}\,\cspecial(\sqrt\theta)\,\,\texte^{-\alpha\sqrt\theta\,\frakd(x)Y}\,  Z_{\sqrt{\theta}\,}^{D,0}(\textd x)\otimes\delta_Y(\textd t),
\end{equation}
where $Y=\NN(0,\sigma_D^2)$ is independent of~$Z^{D,0}_{\sqrt\theta}$.
\end{proposition}

\begin{proofsect}{Proof}
The proof of \cite[Theorem~2.5]{AB} carries over essentially \myemph{verbatim}.
\end{proofsect}

\newcommand{\ee}{\texte}
\newcommand{\ind}{1}

\newcommand{\cO}{\OO}
\newcommand{\cN}{\NN}
\newcommand{\cE}{\EE}

\section{Fixed total time}
\label{sec5}\noindent
Equipped with the enhanced limit results that include the limit value of suitably-norma\-lized fluctuations of the total local time, we now proceed to derive from these the corresponding conclusions for a fixed total time.
We keep working with the random walk started at the boundary vertex~$\varrho$; general starting points will be dealt with in Section~\ref{sec6}. 

\subsection{Time conversions}
The transition from a fixed local time at~$\varrho$ to a fixed total time is based on a simple inversion formula. Recall that, in our context,
\begin{equation}
\hat\tau_\varrho(t):=\inf\bigl\{s\ge0\colon \wt L_s^{D_N}(\varrho)\ge t\bigr\}
\end{equation}
and $\deg(D_N)=\sum_{x\in D_N\cup\{\varrho\}}\deg(x)$. Given a sequence $\{t_N\}_{N\ge1}$ with $t_N\ge1$, define 
\begin{equation}
\label{E:5.2}
	t_{N}^\star = \inf\bigl\{ t \geq 0 \colon \hat{\tau}_{\varrho}(t) \geq \deg(D_N) t_{N} \bigr\}.
\end{equation}
This is an inverse of~$\hat\tau_\varrho$ evaluated at $\deg(D_N) t_N$ and so we expect $\hat\tau_\varrho(t_N^\star)\approx\deg(D_N) t_N$. By \eqref{E:LVt} and \eqref{E:2.27ia}, we should therefore have $\wt L_{\deg(D_N) t_N}^{D_N}(\cdot)\approx \wh L^{D_N}_{t_N^\star}(\cdot)$. Besides their approximate nature, any use of these identifications are complicated by the appearance of the random time~$t_N^\star$ for which we have no better formula than \eqref{E:5.2}. We will thus base the time conversion on a slightly different (still random) quantity that will turn out to be better adapted to our needs.

Recall the definition of~$T_N$ from \eqref{E:4.8}. We note that this actually coincides with the value of $T_N(t_N)$, where (in accord with \eqref{E:4.4u}) we set
\begin{equation}
	T_{N}(t) := \frac{U_{N}(t)}{\sqrt{2t}}
	\qquad \text{for} \qquad
	U_{N}(t) := \frac{1}{|D_{N}|}\sum_{x\in D_{N}}\bigl[\,\wh L^{D_{N}}_{t}(x) - t\bigr].
\end{equation}
Now let
\begin{equation}
	\label{E:tNcirc}
	t^{\circ}_{N} := t_{N} - \left(1-\frac{\deg(\varrho)}{\deg(D_{N})}\right) \sqrt{2t_{N}} \, T_{N}(t_{N}).
\end{equation}
We then have:

\begin{proposition}[Time conversion]
\label{P:tNbound}
Fix any sequence $(b_{N})_{N \geq 1}$ in $(0, \infty)$ such that $b_{N} \to \infty$ and $b_{N}/t_{N}^{1/4} \to 0$ as $N \to \infty$. Then there exist constants $c_{1} > 0 $ such that
\begin{equation}
		\label{E:tNbound}
		\hat{\tau}_{\varrho}\big(t_{N}^{\circ} - b_{N} t_{N}^{1/4}\big)
		\leq \deg(D_N) t_{N}
		\leq \hat{\tau}_{\varrho}\big(t^{\circ}_{N} + b_{N} t_{N}^{1/4}\big)
\end{equation}
and thus, in particular,
\begin{equation}
\label{E:5.6i}
\wh L^{D_N}_{t_N^\circ-b_Nt_N^{1/4}}(\cdot)\le\wt L^{D_N}_{\deg(D_N) t_N}(\cdot)\le\wh L^{D_N}_{t_N^\circ+b_Nt_N^{1/4}}(\cdot)
\end{equation}
hold true with $P^{\varrho}$-probability at least $1- c_{1} b_{N}^{-1}$.
\end{proposition}

The proof will be split into several intermediate results, some of which will be useful later as well. The first item to note is the ``stability'' (or slow variation) of the fluctuation of the total local time:

\begin{lemma}
	\label{L:UN_max}
	There exists a constant $ c_{2} > 0 $ such that for all $s, t \geq 0$ and all $ r > 0 $,
\begin{equation}
		\label{E:UN_max}
		P^{\varrho}\left(\sup_{0\leq u \leq t} \left| U_{N}(s+u) - U_{N}(s) \right| \geq r \right)
		\leq \frac{c_{2}t}{r^{2}}.
\end{equation}
\end{lemma}

\begin{proofsect}{Proof}
	Note that $ U_{N} $ is a compensated compound Poisson process. In view of stationarity, it suffices to consider the case $ s = 0 $. Moreover, since $ U_{N} $ is a martingale, Doob's maximal inequality is applicable and hence
\begin{equation}
		P^{\varrho}\left(\sup_{0\leq u \leq t} \left| U_{N}(u) \right| \geq r \right)
		\leq \frac{4 \Var_{P^\varrho}(U_{N}(t))}{r^{2}}.
\end{equation}
It suffices to show that $ \Var_{P^\varrho}(U_{N}(t)) $ is bounded by $ Ct $ for some $ C > 0 $. To this end, we note that $t \mapsto (U_{N}(t) + t) $ is a compound Poisson process with rate $ \deg(\varrho) $ and jump size distributed as $ \sum_{x \in D_{N}} \ell(x) / |D_{N}| $, where $ \ell(\cdot) $ is the local time for a single excursion. Hence,
\begin{equation}
		\Var_{P^\varrho}(U_{N}(t))
		= \Var_{P^\varrho}(U_{N}(t)+t)
		= \frac{1}{|D_N|^2} \deg(\varrho) t\, E^{\varrho}\Biggl[\biggl(\sum_{x \in D_{N}} \ell(x) \biggr)^{2}\Biggr].
\end{equation}
	The last expectation can be computed via the Kac moment formula,
\begin{equation}
		\label{E:5.7}
		\Var_{P^\varrho}(U_{N}(t)) = \frac{2t}{|D_{N}|^{2}} \sum_{x,y \in D_{N}} G^{D_{N}}(x, y).
\end{equation}
The uniform bound $G^{D_N}(x,y)\le g\log\frac{N}{|x-y|+1}+c$ shows that the sum is at most a constant times~$|D_N|^2$, uniformly in~$N\ge1$. 
\end{proofsect}

 The next lemma quantifies the difference between $ \hat{\tau}_{\varrho}(t_{N}^\star) $ and $ \deg(D_N) t_{N} $:

\begin{lemma}
	\label{L:tN*asymp1}
	Let $ (b_{N})_{N\geq 1} $ be as in the statement of Proposition \ref{P:tNbound}. Then there exists a constant $ c_{3} > 0 $ such that
\begin{equation}
		\label{E:tN*asymp1}
		\left| \frac{\hat{\tau}_{\varrho}(t_{N}^\star)}{\deg(D_N)} - t_{N} \right| \leq b_{N}
		\qquad \text{\rm and} \qquad
		\left|t_{N}^\star - t_{N}\right| < b_{N}\sqrt{t_{N}}
\end{equation}
	hold with $P^{\varrho}$-probability at least $1 - c_{3} b_{N}^{-2}$.
\end{lemma}

\begin{proofsect}{Proof}
 Note that~$\hat\tau_\varrho(t)=\sum_{x\in D_N\cup\{\varrho\}}\deg(x)\wh L^{D_N}_t(x)$. 
	The proof is a straightforward application of Chebyshev's inequality together with some variance estimates. We begin by noting that $ \hat{\tau}_{\varrho}(t_{N}^\star) - \deg(D_N) t_{N} $ is the first time to hit $ \varrho $ starting from the point $\wt X_{\deg(D_N) t_{N}}$. Writing~$H_\varrho$ for the first hitting time of~$\varrho$, the Markov property tells 
\begin{equation}
		E^{\varrho}\left[ \left(\hat{\tau}_{\varrho}(t_{N}^\star) - \deg(D_N) t_{N} \right)^{2} \right]
		= E^{\varrho}\left[ E^{\wt X_{\deg(D_N) t_{N}}}\big[ H_{\varrho}^{2} \big] \right]
		\leq \max_{x \in D_{N}} E^{x}\big[ H_{\varrho}^{2} \big].
\end{equation}
	As in the proof of the previous lemma, applying the Kac moment formula shows
\begin{equation}
		E^{x}\big[ H_{\varrho}^{2} \big]
		= 2 \sum_{y,z\in D_{N}} \deg(y)\deg(z) G^{D_{N}}(x, y)G^{D_{N}}(y, z)
		\leq c_{4}|D_{N}|^{2}
\end{equation}
	for some absolute constant $ c_{4} > 0 $. (This also conforms to the knowledge that the length of a typical excursion on $D_N$ is comparable to the volume of $D_N$.) Then by the Chebyshev inequality,
\begin{equation}
		\label{E:5.12}
		P^{\varrho}\left( \left| \frac{\hat{\tau}_{\varrho}(t_{N}^\star)}{\deg(D_N)} - t_{N} \right| \geq b_{N} \right)
		\leq \frac{c_{4}|D_{N}|^{2}}{(\deg(D_N) b_{N})^{2}}
		\leq \frac{c_{4}}{16 b_{N}^{2}},
\end{equation}
	where the last step follows from $ \deg(D_N) = \deg(\varrho) + 4|D_{N}| $. Also, by the computation similar to the previous proof, we get
\begin{equation}
		E^{\varrho}\Biggl[\biggl(\frac{\hat{\tau}_{\varrho}(t)}{\deg(D_N)} - t\biggr)^2\Biggr]
		= \frac{2t}{\deg(D_N)^{2}}\sum_{x,y \in D_{N}} \deg(x)\deg(y) G^{D_{N}}(x, y)
		\leq c_{5} t
\end{equation}
	for some constant $ c_{5} > 0 $. So again, by Chebyshev's inequality,
\begin{equation}
		\label{E:5.14}
		P^{\varrho}\bigg( \frac{\hat{\tau}_{\varrho}(t_{N} - b_{N}\sqrt{t_{N}})}{\deg(D_N)} \geq t_{N} - b_{N}\sqrt{t_{N}}/2 \bigg)
		\leq \frac{c_{5}(t_{N} - b_{N}\sqrt{t_{N}})}{(b_{N}\sqrt{t_{N}}/2)^{2}}
		\leq \frac{4c_{5}}{b_{N}^{2}}
\end{equation}
	and likewise
\begin{equation}
		\label{E:5.15}
		P^{\varrho}\bigg( \frac{\hat{\tau}_{\varrho}(t_{N} + b_{N}\sqrt{t_{N}})}{\deg(D_N)} \leq t_{N} + b_{N}\sqrt{t_{N}}/2 \bigg)
		\leq \frac{4c_{5}(1 + b_{N}t_{N}^{-1/2})}{b_{N}^{2}}.
\end{equation}
	Combining \eqref{E:5.12}, \eqref{E:5.14}, and \eqref{E:5.15} we find that there exists a constant $ c_{3} > 0 $, depending only on $ (t_{N})_{N\geq 1} $ and $ (b_{N})_{N\geq 1} $, such that all of
\begin{equation}
		\label{E:5.16}
		\begin{aligned}
			\hat{\tau}_{\varrho}(t_{N} - b_{N}\sqrt{t_{N}}) &< \deg(D_N)\big(t_{N} - b_{N}\sqrt{t_{N}}/2\big), \\
			\hat{\tau}_{\varrho}(t_{N} + b_{N}\sqrt{t_{N}}) &> \deg(D_N)\big(t_{N} + b_{N}\sqrt{t_{N}}/2\big), \\
			\left| \hat{\tau}_{\varrho}(t_{N}^\star) - \deg(D_N) t_{N} \right| &\leq \deg(D_N) b_{N}/2
		\end{aligned}
\end{equation}
	simultaneously hold with $P^{\varrho}$-probability at least $1 - c_{3}b_{N}^{-2}$. But if all of \eqref{E:5.16} hold, then we get
\begin{equation}
		\hat{\tau}_{\varrho}(t_{N} - b_{N}\sqrt{t_{N}})
		< \hat{\tau}_{\varrho}(t_{N}^\star)
		< \hat{\tau}_{\varrho}(t_{N} + b_{N}\sqrt{t_{N}}).
\end{equation}
	By the monotonicity of $ \hat{\tau}_{\varrho} $, these altogether imply \eqref{E:tN*asymp1} as required.
\end{proofsect}

Next we will quantify the difference between $t_N^\star$ and~$t_N^\circ$:

\begin{lemma}
	\label{L:tN*asymp2}
	Assume $t_N\ge1$ and let $ (b_{N})_{N\geq 1} $ be as in the statement of Proposition \ref{P:tNbound}. Then there exists a constant $ c_{6} > 0 $ such that
\begin{equation}
		\label{E:tN*asymp2}
		|t_{N}^\star - t_{N}^{\circ}| \leq b_{N}t_{N}^{1/4}
\end{equation}
	holds with $P^{\varrho}$-probability at least $1 - c_{6} b_{N}^{-1}$.
\end{lemma}

\begin{proofsect}{Proof}
We note that, by \eqref{E:2.27ia} and the fact that $\deg(x)=4$ for~$x\in D_N$, 
\begin{equation}
		\begin{split}
		U_{N}(t)
		&= \frac{1}{|D_{N}|} \sum_{x \in D_{N}} \left( \frac{1}{\deg(x)} \int_{0}^{\hat{\tau}_{\varrho}(t)}   1_{\{\wt X_{s} = x\}} \, \dd s - t \right) \\
		&\quad= \frac{1}{|D_{N}|} \left( \frac{1}{4} \left( \hat{\tau}_{\varrho}(t) - \deg(\varrho)t \right) - |D_{N}|t \right)
		= \frac{\hat{\tau}_{\varrho}(t) - \deg(D_N) t}{4|D_{N}|}.
		\end{split}
\end{equation}
	Rearranging the identity in terms of $ t $, we get
\begin{equation}
		\label{E:5.20}
		t = \frac{\hat{\tau}_{\varrho}(t)}{\deg(D_N)} - \left(1-\frac{\deg(\varrho)}{\deg(D_N)}\right) U_{N}(t).
\end{equation}
	This will be used to prove the desired bound. Plugging $ t := t_{N}^\star $, we notice that the right-hand side of \eqref{E:5.20} 
almost looks like the definition \eqref{E:tNcirc} of $ t_{N}^{\circ} $, except that we need~$t_{N} $ in place of $ \hat{\tau}_{\varrho}(t_{N}^\star) / \deg(D_N)$ and $ U_{N}(t_{N}) $ in place of $ U_{N}(t_{N}^\star) $. This amounts to estimating their respective differences, and this is where the previous lemmas come handy.
	
First, we plug $s := t_{N} - b_{N}\sqrt{t_{N}}$ and $t := 2b_{N}\sqrt{t_{N}}$ in \eqref{E:UN_max} to get
\begin{equation}
		P^{\varrho}\left( \sup_{|u| \leq b_{N}\sqrt{t_{N}}} \left| U_{N}(t_{N}+u) - U_{N}(t_{N}) \right| \geq b_{N}t_{N}^{1/4} \right)
		\leq \frac{ 8 c_{2}b_{N}\sqrt{t_{N}}}{\big(b_{N}t_{N}^{1/4}\big)^{2}}
		= \frac{ 8 c_{2}}{b_{N}}.
\end{equation}
	Combining this with Lemma \ref{L:tN*asymp1}, we can find $ c_{7} > 0 $ such that both \eqref{E:tN*asymp1} and
\begin{equation}
		\label{E:5.22}
		\left| U_{N}(t_{N}+u) - U_{N}(t_{N}) \right| \leq b_{N}t_{N}^{1/4}
		\qquad \text{for all } |u| \leq b_{N}\sqrt{t_{N}}
\end{equation}
	hold with $ P^{\varrho} $-probability at least $ 1-c_{7}b_{N}^{-1} $. Moreover, given \eqref{E:tN*asymp1} and \eqref{E:5.22}, we also get $ \left| U_{N}(t_{N}^\star) - U_{N}(t_{N}) \right| \leq b_{N}t_{N}^{1/4} $. Putting this together, we get
\begin{equation}
	\begin{aligned}
		|t_{N}^\star - t_{N}^{\circ}|
		&\leq \left| \frac{\hat{\tau}_{\varrho}(t_{N}^\star)}{\deg(D_N)} - t_{N} \right| + \left| U_{N}(t_{N}^\star) - U_{N}(t_{N}) \right|
		\\
		&\leq b_{N}\big(1+t_{N}^{1/4}\big)\le 2b_N t_N^{1/4}.
	\end{aligned}
\end{equation}
	Although this bound is slightly larger than that appearing in the statement, we can repeat all the above argument with $\{b_{N}/2\}_{N\geq 1} $ in place of $\{b_{N}\}_{N\geq 1} $, then the desired claim follows with $ c_{6} = 2c_{7} $.
\end{proofsect}

We are now ready to prove the main statement:

\begin{proofsect}{Proof of Proposition \ref{P:tNbound}}
	Let $ (b_{N})_{N\geq 1} $ be as in the statement. Then by the definition of $ t_{N}^\star $ and Lemma \ref{L:tN*asymp2},
\begin{equation}
		\deg(D_N) t_{N} \leq \hat{\tau}_{\varrho}(t_{N}^\star) \leq \hat{\tau}_{\varrho}\big(t_{N}^{\circ} + b_{N}t_{N}^{1/4}\big)
\end{equation}
	holds with $ P^{\varrho} $-probability at least $ 1-\cO(b_{N}^{-1}) $. Next, regarding $ t_{N} \mapsto t_{N}^\star $ and $ t_{N} \mapsto t_{N}^{\circ} $ as functions of $ t_{N} $ for each fixed $ N $, Lemma \ref{L:tN*asymp1} applied to $ (t_{N}-b_{N}/4)_{N\geq 1} $ and $ (b_{N}/4)_{N\geq 1} $ in place of $ (t_{N})_{N\geq 1} $ and $ (b_{N})_{N\geq 1} $, respectively, show that both
\begin{equation}
\deg(D_N) t_{N}
			\geq \hat{\tau}_{\varrho}\big( (t_{N} - b_{N}/4)^\star \big) \geq \deg(D_N)(t_{N} - b_{N}/2)
\end{equation}
and
\begin{equation}
\big\lvert \big( t_{N} - b_{N}/4 \big)^\star - t_{N} \big\rvert \leq b_{N}\sqrt{t_{N}}/2
\end{equation}
are satisfied with $ P^{\varrho} $-probability at least $ 1-\cO(b_{N}^{-1}) $. Then using \eqref{E:5.22} and repeating the argument as in the previous proof, we can bound $ (t_{N} - b_{N}/4)^\star $ from below by $ t_{N}^{\circ} - b_{N}t_{N}^{1/4} $ again with probability at least $ 1-\cO(b_{N}^{-1}) $.
\end{proofsect}

\subsection{Continuous-time exceptional level sets}
We are now ready to adapt the convergence theorems for the exceptional level-set measures for the boundary-vertex local times $\wh L^{D_{N}}$ to those associated with the local time $\wt L^{D_{N}}$ of the continuous-time walk~$\wt X$ run for a fixed time of order~$N^2(\log N)^2$. We begin by the thick points; the arguments will be readily adapted to the other families of exceptional points as well. Given two positive sequences $\{t_{N}\}_{N\geq 1}$ and $\{a_{N}\}_{N\geq 1}$ as before, define
\begin{equation}
\label{E:5.29b}
\wt\zeta^{D}_{N} = \frac{1}{W_{N}} \sum_{x \in D_{N}} \delta_{x/N} \otimes \delta_{(\wt L^{D_{N}}_{\deg(D_N) t_{N}}(x) - a_{N})/\sqrt{2a_{N}}},
\end{equation}
where $W_{N}$ is the same as in the case of $\wh\zeta^{D}_{N}$. Then

\begin{proposition}[Continuous-time thick points]
\label{thm-thick-cont}
Under the setting and notation of Theorem~\ref{thm-thick} and for the walk started at the ``boundary vertex,'' we have
\begin{equation}
\label{E:1.21b}
\wt\zeta^D_N\,\,\,\underset{N\to\infty}\Lawarrow\,\,\, \sqrt{\frac{\sqrt{\theta}}{\sqrt{\theta}+\lambda}}\,\,\cspecial(\lambda)\,\ee^{\alpha \lambda (\mathfrak{d}(x) - 1) T}\,Z_\lambda^{D,0}(\textd x)\otimes\texte^{-\alpha\lambda h}\textd h,
\end{equation}
where $T$ and $Z_{\lambda}^{D,0}$ are independent with $T \sim \cN(0, \sigma_{D}^{2})$.
\end{proposition}

The key point is to carefully track the effects of the random time shift $\sqrt{2t_{N}} \, T_{N}$ in the quantity~$t_{N}^{\circ}$ from \eqref{E:tNcirc}. Let~$\{b_N\}_{N\ge1}$ be a sequence with~$b_N\to\infty$ and~$b_N/t_N^{1/4}\to0$. Consider the event
\begin{multline}
\label{E:evtgood}
\quad
\cE_N:= \left\{ \hat{\tau}_{\varrho}\big(t_{N}^{\circ} - b_{N} t_{N}^{1/4}\big)
			\leq \deg(D_N) t_{N}
			\leq \hat{\tau}_{\varrho}\big(t^{\circ}_{N} + b_{N} t_{N}^{1/4}\big) \right\}
\\
\cap \,\left\{ \max_{|u|\leq b_{N}\sqrt{t_{N}}}\left| U_{N}(t_{N}+u) - U_{N}(t_{N}) \right| \leq b_{N}t_{N}^{1/4}\right\} 
\cap \left\{ \left| T_{N} \right| \leq {b_{N}} \right\}.
\quad
\end{multline}
We then have:

\begin{lemma}
\label{lemma-5.6}
There is a constant~$c_7>0$ such that the following holds for all~$N\ge1$:
\begin{equation}
\label{E:5.32i}
P^{\varrho}(\cE_N) \geq 1 - c_7 b_{N}^{-1}
\end{equation}
and
\begin{equation}
\label{E:5.33i}
\max_{|u|\le b_N\sqrt{t_N}}\,\left|T_{N}(t_{N}+u) - T_{N}\right|
\le c_7 b_N/t_N^{1/4}\quad\text{\rm on }\EE_N.
\end{equation}
\end{lemma}

\begin{proofsect}{Proof}
The bound \eqref{E:5.32i} follows from Proposition \ref{P:tNbound}, Lemma \ref{L:UN_max} and the fact that~$T_N$ has asymptotically a Gaussian tail. To get \eqref{E:5.33i}, note that for $|u|\le b_N\sqrt{t_N}$,
\begin{equation}
		\label{E:5.29}
		\left|T_{N}(t_{N}+u) - T_{N}\right|
		\leq \frac{b_{N}t_{N}^{1/4}}{\sqrt{2(t_{N} - b_{N}\sqrt{t_{N}})}} + \frac{b_{N} \left|T_{N}\right|}{\sqrt{t_{N} - b_{N}\sqrt{t_{N}}}}.
\end{equation}
As $|T_N|\le b_N$ on~$\cE_N$ and $\{b_N/t_N^{1/4}\}_{N\ge1}$ is bounded, this is at most order $b_{N}/t_{N}^{1/4}$.
\end{proofsect}

The argument to follow will be based on dividing the event~$\cE_N$ depending on the values of $T_{N}$. For this we fix an  $\epsilon > 0$, and let $ \{ \rho_{k} \}_{k\in\Z} $ be a family of continuous functions such that
\begin{equation}
		\label{E:evtdiv}
		0 \leq \rho_{k} \leq \ind{}_{[(k-1)\epsilon,(k+1)\epsilon]}
		\qquad \text{and} \qquad
		\sum_{k\in\Z} \rho_{k} = 1.
\end{equation}
We also define two auxilliary time sequences $\{t^{+}_{N,k}\}_{N\geq 1}$ and $\{t^{-}_{N,k}\}_{N\geq 1}$ by
\begin{equation}
		\label{E:tNshift}
		\begin{gathered}
			t^{+}_{N,k}
			= t_{N} - \Big(1 - \tfrac{\deg(\varrho)}{\deg(D_N)} \Big) \epsilon (k-1) \sqrt{2t_{N}} + b_{N}t_{N}^{1/4},\\
			t^{-}_{N,k}
			= t_{N} - \Big(1 - \tfrac{\deg(\varrho)}{\deg(D_N)} \Big) \epsilon (k+1) \sqrt{2t_{N}} - b_{N}t_{N}^{1/4}.
		\end{gathered}
\end{equation}
We then have:

\begin{lemma}
\label{lemma-5.7}
For each~$M>0$ there is~$N_0\in\N$ such that for all~$N\ge N_0$ and all~$k\in\Z$ with~$|k|\le M$, the following holds on $\cE_N\cap\{T_N\in\supp(\rho_k)\}$:
\begin{equation}
\label{E:5.34}
\bigl|T_N(t_{N,k}^\pm)-T_N\bigr|\le c_7 b_N/t_N^{1/4}
\end{equation}
and
\begin{equation}
		\label{E:5.33b}
		\wh L^{D_{N}}_{t^{-}_{N,k}}(\cdot)
		\leq \wt L^{D_{N}}_{\deg(D_N) t_{N}}(\cdot)
		\leq \wh L^{D_{N}}_{t^{+}_{N,k}}(\cdot).
\end{equation}
\end{lemma}

\begin{proofsect}{Proof}
Fix~$M>0$. As~$b_N\to\infty$ and $b_N t_N^{-1/4}\to0$, we can choose~$N_0\in\N$ such that $\epsilon (M+1)\sqrt{2t_N}+b_N t_N^{1/4}\le b_N\sqrt{t_N}$ for all~$N\ge N_0$. Then for all~$N\ge N_0$,
\begin{equation}
\bigl|t_{N,k}^\pm-t_N\bigr|\le b_N\sqrt{t_N},\quad -M\le k\le M.
\end{equation}
The bound \eqref{E:5.34} is then implied by \eqref{E:5.33i}. 

For \eqref{E:5.33b} we note that, on $\{T_N\in\supp(\rho_k)\}$ we have $(k-1)\epsilon\le T_N\le (k+1)\epsilon$ and  thus also
\begin{equation}
\label{E:5.40nw}
t_{N,k}^-\le t_N^\circ-b_N t_N^{1/4}\le t_N^\circ+b_N t_N^{1/4}\le t_{N,k}^+.
\end{equation}
 The bound \eqref{E:5.33b} then follows from the inequalities in \eqref{E:evtgood} and the monotonicity of $t\mapsto\wh L^{D_N}_t(\cdot)$.
\end{proofsect}

The inequalities \eqref{E:5.33b} thus naturally make us consider the level-set measures $\wh\zeta^{D}_{N}$ along different choices of time sequences than the base sequence~$\{t_N\}_{N\ge1}$. We will explicate the dependence on the time sequence by writing $\wh\zeta^{D}_{N}(t'_{N})$ whenever it is along $\{t'_{N}\}_{N\geq 1}$ rather than  $\{t_{N}\}_{N\geq 1}$, and likewise, we will write $W_{N}(t'_{N})$ for the normalizing constants along $\{t'_{N}\}_{N\geq 1}$.
Next we note:

\begin{lemma}
\label{lemma-5.8}
We have $\deg(\varrho)/\deg(D_N)\to0$ as~$N\to\infty$. In particular, for each~$k\in\Z$,
\begin{equation}
\label{E:5.42i}
t^{\pm}_{N,k} \sim 2g\theta (\log N)^{2},\quad N\to\infty.
\end{equation}
Moreover,
\begin{equation}
		\label{E:5.32}
		\begin{gathered}
			W_{N}(t^{+}_{N,k}) = W_{N}(t_N) \, \ee^{-\alpha\lambda \epsilon (k-1) + o(1)}, \\
			W_{N}(t^{-}_{N,k}) = W_{N}(t_N) \, \ee^{-\alpha\lambda \epsilon (k+1) + o(1)},
		\end{gathered}
\end{equation}
where~$o(1)\to0$ uniformly in~$k\in\Z$ with~$|k|\le M$, for any~$M>0$.
\end{lemma}

\begin{proofsect}{Proof}
We start by showing $\deg(\varrho)/\deg(D_N)\to0$. For this we note that~$\deg(D_N)\ge 4|D_N|$ while, for any~$\delta>0$ and~$N$ sufficiently large, $\deg(\varrho)\le 4|D_N\smallsetminus D_N^\delta|$, where $D_N^\delta:=\{x\in D_N\colon d_\infty(x,D_N^\cc) >\delta N\}$. Definition~\ref{dfn:admissible} now ensures
\begin{equation}
\limsup_{N\to\infty}\frac{\deg(\varrho)}{\deg(D_N)}\le
\limsup_{N\to\infty}\frac{|D_N\smallsetminus D_N^\delta|}{|D_N|}\le\frac{\leb(D\smallsetminus D^{2\delta})}{\leb(D)},
\end{equation}
where~$D^\delta:=\{x\in D\colon d_\infty(x,D^\cc)>\delta\}$. As~$D^{2\delta}\uparrow D$ as~$\delta\downarrow0$, we have $\leb(D\smallsetminus D^{2\delta})\to0$ as~$\delta\downarrow0$. 

With $\deg(\varrho)/\deg(D_N)\to0$ settled, the asymptotic \eqref{E:5.42i} is now checked readily from the definition of~$t^\pm_{N,k}$. The bounds in \eqref{E:5.32} follow similarly from the explicit formula for~$W_N$ and some routine estimates.
\end{proofsect}

We are now ready for:
 
\begin{proofsect}{Proof of Proposition~\ref{thm-thick-cont}}
Let $f\colon\overline D\times(\R\cup\{+\infty\})\to[0,\infty)$ be a bounded and continuous function that is non-decreasing in the second coordinate and supported on~$\overline D\times[b,\infty]$ for some~$b\in\R$. Then  \eqref{E:5.32}, \eqref{E:5.33b} and \eqref{E:5.34} show
\begin{equation}
\begin{aligned}
\label{E:5.35}
\ee^{-2\alpha\lambda \epsilon + o(1)} \ee^{-\alpha\lambda T_{N}(t^{-}_{N,k})} \langle \wh\zeta_{N}^{D}(t^{-}_{N,k}), f \rangle 
			& \leq \langle \wt\zeta_{N}^{D}, f \rangle 
			\\
			&\leq \ee^{2\alpha\lambda \epsilon + o(1)} \ee^{-\alpha\lambda T_{N}(t^{+}_{N,k}) } \langle \wh\zeta_{N}^{D}(t^{+}_{N,k}), f \rangle
\end{aligned}
\end{equation}
on $\cE_N\cap \{ T_{N} \in \operatorname{supp}(\rho_{k}) \}$, where~$o(1)$ is a deterministic sequence tending to zero uniformly in~$k\in\Z$ with~$|k|\le M$.

Define the maximal modulus of continuity of~$\{\rho_k\colon|k|\le M\}$ by
\begin{equation}
\text{osc}_{M,\epsilon}(r):=\max_{|k|\le M}\sup_{\begin{subarray}{c}
t,t'\in\R\\|t-t'|\le r
\end{subarray}}
\bigl|\rho_k(t)-\rho_k(t')\bigr|.
\end{equation}
Relying first on the lower bound of \eqref{E:5.35}, we now estimate
\begin{equation}
\label{E:5.46i}
\begin{aligned}
E^{\varrho}\bigl( \ee^{-\langle \wt\zeta_{N}^{D}, f \rangle}& \bigr)-P^{\varrho}(\cE^{\mathrm{c}}_N)
			- P^{\varrho}\bigl(|T_{N}| \geq M/\epsilon\bigr)
			\\
&\leq \sum_{k=- M}^M E^{\varrho} \bigl( \ee^{-\langle \wt\zeta_{N}^{D}, f \rangle} \rho_{k}(T_{N})\ind{}_{\cE_N} \bigr)
		  \\
			&\leq \sum_{k=- M}^M E^{\varrho} \biggl(\ee^{-\ee^{-2\alpha\lambda \epsilon + o(1)} \ee^{-\alpha\lambda T_{N}(t^{-}_{N,k})} \langle \wh\zeta_{N}^{D}(t^{-}_{N,k}), f \rangle}\rho_{k}(T_{N}) \ind_{\cE_N}\biggr)
			\\
			&\leq (2M+1)\text{osc}_{M,\epsilon}\bigl(c_7b_N/t_N^{1/4}\bigr)
			\\&\qquad\quad+
			\sum_{k=- M}^M E^{\varrho} \biggl(\ee^{-\ee^{-2\alpha\lambda \epsilon + o(1)} \ee^{-\alpha\lambda T_{N}(t^{-}_{N,k})} \langle \wh\zeta_{N}^{D}(t^{-}_{N,k}), f \rangle}\rho_{k}\bigl(T_{N}(t_{N,k}^-)\bigr)\ind_{\cE_N} \biggr),
\end{aligned}
\end{equation}
where in the last step we used \eqref{E:5.34}. The key point is that, dropping the indicator of~$\cE_N$, the $k$-th term in the sum is now a continuous function of the process~$\wh\zeta_N^D(t_{N,k}^-)$ and the time~$T_N(t_{N,k}^-)$. In light of \eqref{E:5.42i}, Proposition~\ref{thm-4.3} gives 
\begin{equation}
E^{\varrho} \biggl(\ee^{-\ee^{-2\alpha\lambda \epsilon + o(1)} \ee^{-\alpha\lambda T_{N}(t^{-}_{N,k})} \langle \wh\zeta_{N}^{D}(t^{-}_{N,k}), f \rangle}\rho_{k}\bigl(T_{N}(t_{N,k}^-)\bigr)\biggr)
\underset{N\to\infty}\longrightarrow\, E\Bigl(\ee^{-\ee^{-2\alpha\lambda \epsilon}\texte^{-\alpha\lambda T}\langle\wh\zeta^D,f\rangle}\rho_k(T)\Bigr),
\end{equation}
where
\begin{equation}
\label{E:5.48a}
\wh\zeta^D := \sqrt{\frac{\sqrt{\theta}}{\sqrt{\theta}+\lambda}} \,\cspecial(\lambda)\, \ee^{\alpha \lambda \mathfrak{d}(x)T} Z_{\lambda}^{D,0}(\textd x) \otimes \ee^{-\alpha \lambda h} \textd h
\end{equation}
with~$T=\NN(0,\sigma_D^2)$ independent of~$Z_{\lambda}^{D,0}$. Dropping the restriction to~$|k|\le M$, the $N\to\infty$ \myemph{limes superior} of the sum on the extreme right of \eqref{E:5.46i} is then at most $E(\ee^{- \texte^{- 2 \alpha \lambda \epsilon}\texte^{-\alpha\lambda T}
\langle\wh\zeta^D,f\rangle})$. Since $\text{osc}_{M,\epsilon}(r)\to0$ as~$r\downarrow0$, taking $N\to\infty$ followed by~$M\to\infty$ and~$\epsilon\downarrow0$ shows
\begin{equation}
\limsup_{N\to\infty}
E^{\varrho}\bigl( \ee^{-\langle \wt\zeta_{N}^{D}, f \rangle}\bigr)
\le E\bigl(\ee^{- \texte^{-\alpha\lambda T}\langle\wh\zeta^D,f\rangle}\bigr),
\end{equation}
where the two ``error'' terms on the left-hand side of \eqref{E:5.46i} tend to zero in the stated limits thanks to Lemma~\ref{lemma-5.6} and the Gaussian (asymptotic) tail of~$T_N$.

The argument for  a corresponding lower bound is very similar; we need to work with~$t_{N,k}^+$ instead of~$t_{N,k}^-$ and use explicit estimates to get rid of the indicator~$\ind_{\cE_N}$ and the restriction to the range of~$k$ in the sum. As a conclusion, we get
\begin{equation}
\lim_{N\to\infty}
E^{\varrho}\bigl( \ee^{-\langle \wt\zeta_{N}^{D}, f \rangle}\bigr)
= E\bigl(\ee^{- \texte^{-\alpha\lambda T}\langle\wh\zeta^D,f\rangle}\bigr)
\end{equation}
for any function~$f$ as above. 
This is sufficient to give $\wt\zeta_N^D\Lawarrow\texte^{-\alpha\lambda T}\wh\zeta^D$, as desired.
\end{proofsect}

For the thin points we now get:

\begin{proposition}[Continous-time thin points]
\label{thm-thin-cont}
Under the setting and notation of Theorem~\ref{thm-thin} and for the walk started at the ``boundary vertex,'' we have
\begin{equation}
\label{E:1.21}
\wt\zeta^D_N\,\,\,\underset{N\to\infty}\Lawarrow\,\,\, \sqrt{\frac{\sqrt{\theta}}{\sqrt{\theta}-\lambda}}\,\,\cspecial(\lambda)\,\ee^{-\alpha \lambda (\mathfrak{d}(x) - 1) T}\,Z_\lambda^{D,0}(\textd x)\otimes\texte^{+\alpha\lambda h}\textd h,
\end{equation}
where $T$ and $Z_{\lambda}^{D,0}$ are independent with $T \sim \cN(0, \sigma_{D}^{2})$.
\end{proposition}

\begin{proofsect}{Proof}
The argument is similar to that for the thick points: We need to work with compact\-ly-supported, continuous test functions $f\colon\overline D\times(\R\cup\{-\infty\})\to[0,\infty)$ that are non-increasing in the second coordinate. 
The change in monotonicity effectively swaps the inequalities in \eqref{E:5.35} and, due to a sign change in \eqref{E:5.32}, also that in the exponent of $\ee^{-\alpha\lambda T_N(t^\pm_{N,k})}$. We also need to rely on Proposition~\ref{thm-4.4} instead of Proposition~\ref{thm-4.3}. We leave further details to the reader.
\end{proofsect}

Moving to the light points, we define
\begin{equation}
\wt\vartheta^{D}_{N} := \frac{1}{\wh W_{N}} \sum_{x \in D_{N}} \delta_{x/N} \otimes \delta_{\wt L^{D_{N}}_{\deg(D_N) t_{N}}(x)}
\end{equation}
and state:

\begin{proposition}[Continuous-time light points]
\label{thm-light-cont}
Under the setting and assumptions of Theorem~\ref{thm-light} and for the walk started at the ``boundary vertex,'' we have
\begin{equation}
\label{E:2.31ii}
\wt\vartheta^D_N\,\,\,\underset{N\to\infty}\Lawarrow\,\,\,  \sqrt{2\pi g}\,\cspecial(\sqrt\theta)\,\, \ee^{-\alpha \sqrt{\theta} (\mathfrak{d}(x) - 1) T}\,Z_{\sqrt\theta}^{D,0}(\textd x)\otimes\tilde\mu(\textd h),
\end{equation}
where $T=\NN(0,\sigma_D^2)$ is independent of $Z_{\sqrt\theta}^{D,0}$ and $\tilde\mu$ is the measure in \eqref{E:tilde-mu}.
\end{proposition}

\begin{proofsect}{Proof}
Relying on our convention concerning different time sequences, we start by noting
\begin{equation}
		\label{E:5.32u}
		\begin{gathered}
			\wh W_{N}(t^{+}_{N,k}) = \wh W_{N}(t_N) \, \ee^{\alpha\sqrt\theta \epsilon (k-1) + o(1)}, \\
			\wh W_{N}(t^{-}_{N,k}) = \wh W_{N}(t_N) \, \ee^{\alpha\sqrt\theta \epsilon (k+1) + o(1)}.
		\end{gathered}
\end{equation}
Given a compactly-supported, continuous function $f\colon\overline D\times[0,\infty)\to[0,\infty)$ that is non-increasing in the second coordinate, from \eqref{E:5.32u}, \eqref{E:5.33b} and \eqref{E:5.34} we then have
\begin{equation}
\begin{aligned}
\label{E:5.35u}
 \ee^{-2\alpha\sqrt\theta \epsilon + o(1)} \ee^{\alpha\sqrt\theta T_{N}(t^{+}_{N,k}) } \langle \wh\vartheta_{N}^{D}(t^{+}_{N,k}), f \rangle
			& \leq \langle \wt\vartheta_{N}^{D}, f \rangle 
			\\
			&\leq \ee^{2\alpha\sqrt\theta \epsilon + o(1)} \ee^{\alpha\sqrt\theta T_{N}(t^{-}_{N,k})} \langle \wh\vartheta_{N}^{D}(t^{-}_{N,k}), f \rangle.
\end{aligned}
\end{equation}
The rest of the argument for the thick points (with Proposition~\ref{thm-light-bv} instead of Proposition~\ref{thm-4.3}) can now be applied to get
\begin{equation}
\langle \wt\vartheta_{N}^{D}, f \rangle \,\,\underset{N\to\infty}\Lawarrow\,\, \ee^{+\alpha\sqrt\theta T}\langle\wh\vartheta^D,f\rangle,
\end{equation}
where
\begin{equation}
\wh\vartheta^D:=\sqrt{2\pi g}\,\cspecial(\sqrt\theta)\,\,\texte^{-\alpha\sqrt\theta\,\frakd(x)T}\,  Z_{\sqrt{\theta}\,}^{D,0}(\textd x)\otimes\tilde\mu(\textd h).
\end{equation}
The claim now follows by a density argument.
\end{proofsect}

Finally, for the avoided points we set
\begin{equation}
\wt\kappa^{D}_{N} := \frac{1}{\wh W_{N}} \sum_{x \in D_{N}}  1_{\{\wt L^{D_{N}}_{\deg(D_N) t_{N}}(x)=0\}}\,\delta_{x/N}
\end{equation}
and state:

\begin{proposition}[Continuous-time avoided points]
\label{thm-avoid-cont}
Under the setting and assumptions of Theorem~\ref{thm-light} and for the walk started at the ``boundary vertex,'' we have
\begin{equation}
\label{E:2.27cont}
\wt\kappa^D_N\,\,\,\underset{N\to\infty}\Lawarrow\,\,\, \sqrt{2\pi g}\,\cspecial(\sqrt\theta)\,\, \ee^{-\alpha \sqrt{\theta} (\mathfrak{d}(x) - 1) T}\,Z_{\sqrt\theta}^{D,0}(\textd x),
\end{equation}
where $T=\NN(0,\sigma_D^2)$ is independent of $Z_{\sqrt\theta}^{D,0}$.
\end{proposition}

\begin{proofsect}{Proof}
Given a continuous $f\colon\overline D\to\R$, the identity \eqref{E:5.35u}
applies with $\wt\vartheta^D_N$, resp., $\wh\vartheta^D_N$ replaced by $\wt\kappa^D_N$, resp., $\wh\kappa^D_N$. The argument then proceeds as for Proposition~\ref{thm-light-cont}.
\end{proofsect}

\section{Arbitrary starting points}
\label{sec6}\noindent
As our next item of business, we augment the continuous-time conclusions from the previous section to allow the random walk to start at an arbitrary point of~$D_N$. The formal statement is the content of:

\begin{theorem}[Arbitrary starting points]
\label{thm-6.1}
The statements of Propositions~\ref{thm-thick-cont}, \ref{thm-thin-cont}, \ref{thm-light-cont} and \ref{thm-avoid-cont} apply for random walk starting from an arbitrary point~$x_N\in D_N$.
\end{theorem}

We will start with the thick points as that is the hardest case. Assume that~$\{a_N\}_{N\ge1}$ and~$\{t_N\}_{N\ge1}$ satisfy the conditions of Propositions~\ref{thm-thick-cont}. The integrals of $\{\wt\zeta^D_N\colon N\ge1\}$ from \eqref{E:5.29b} against $f\in C_\cc(\overline D\times(\R\cup\{+\infty\}))$ are tight random variables. 
Our strategy is to use the strong Markov property after the first hitting of the ``boundary vertex.'' For this let us recall that~$H_x$ denotes the first hitting time of vertex~$x$ and let~$\theta_t$ denote the shift on the path space acting as $(\wt X\circ\theta_t)_s=\wt X_{t+s}$. We will write $\{(\wt L^{D_N}\circ\theta_t)_s\colon s\ge0\}$ for the local time process associated with the time-shifted path $\{(\wt X\circ\theta_t)_s\colon s\ge0\}$. Our first observation is then:

\begin{lemma}
\label{lemma-6.2}
On $\{H_\varrho<t\}$, we have
\begin{equation}
\label{E:6.1}
\wt L_t^{D_N}(\cdot) = \wt L^{D_N}_{H_\varrho}(\cdot)+(\wt L^{D_N}\circ\theta_{H_\varrho})_{t-H_\varrho}(\cdot).
\end{equation}
In particular, under the conditions of Proposition~\ref{thm-thick-cont}, for any $f\in C_\cc(\overline D\times(\R\cup\{+\infty\}))$ that is non-decreasing in the second variable
and any $x_N\in D_N$,
\begin{equation}
\label{E:6.2}
\limsup_{N\to\infty} E^{x_N}\bigl(\texte^{-\langle\wt\zeta^D_N,f\rangle}\bigr)
\le \lim_{N\to\infty}\,E^{\varrho}\bigl(\texte^{-\langle\wt\zeta^D_N,f\rangle}\bigr).
\end{equation}
\end{lemma}

\begin{proofsect}{Proof}
The relation \eqref{E:6.1} is a direct consequence of the additivity of the local time. As to \eqref{E:6.2},  for~$f$ as above and any~$m>0$ with $t_N>m$, dropping the term $\wt L^{D_N}_{H_\varrho}$ while noting that $W(t_N-m)\ge \texte^{-c \frac{m}{\log N}} W(t_N)$ for some $c > 0$  shows
\begin{equation}
\bigl\langle\wt\zeta^{D_N}_N(t_N),f\bigr\rangle\ge 
\texte^{- c \frac{m}{\log N}}
\bigl\langle\wt\zeta^{D_N}_N(t_N-m),f\bigr\rangle\circ\theta_{H_\varrho}\quad\text{on }\{H_\varrho< m\deg(D_N)\}.
\end{equation}
The strong Markov property then gives
\begin{equation}
\label{E:6.4}
\begin{aligned}
E^{x_N}\bigl(\texte^{-\langle\wt\zeta^D_N(t_N),f\rangle}\bigr)&
\le P^{x_N}\bigl(H_\varrho\ge m\deg(D_N)\bigr)+E^{x_N}\bigl(\texte^{-\langle\wt\zeta^D_N(t_N),f\rangle}1_{\{H_\varrho<m\deg(D_N)\}}\bigr)
\\
&\le 
P^{x_N}\bigl(H_\varrho\ge m\deg(D_N)\bigr)+E^{\varrho}
\bigl(\texte^{-\texte^{-c \frac{m}{\log N}}\langle\wt\zeta^D_N(t_N-m),f\rangle}\bigr).
\end{aligned}
\end{equation}
Since the random walk on~$D_N$ coincides with the random walk on~$\Z^2$ until time~$H_\varrho$, the Central Limit Theorem shows that the probability tends to zero in the limits~$N\to\infty$ and~$m\to\infty$. The expectation on the right converges by Proposition~\ref{thm-thick-cont}.
\end{proofsect}

Our next goal is to prove a complementary bound to \eqref{E:6.2} for the \myemph{limes inferior}. For this we must control the effect of the first term on the right of \eqref{E:6.1}. Writing $\{(\wh L^{D_N}\circ\theta_t)_s\colon s\ge0\}$ for the local time of the process $\wt X\circ\theta_t$ parametrized at the time spent at the boundary vertex, we then have:

\begin{lemma}
\label{lemma-6.3}
Under the conditions of Proposition~\ref{thm-thick-cont}, for each~$b\in\R$ there is~$c>0$ such that for all $N\ge1$ and all~$x\in D_N$,
\begin{equation}
\label{E:6.5uie}
\sum_{z\in D_N}
P^x\biggl(\wt L^{D_N}_{H_\varrho}(z)+(\wh L^{D_N}\circ\theta_{H_\varrho})_{t_N}(z)\ge a_N+b\log N,\, H_z<H_\varrho\biggr)
\le c\frac {W_N}{\log N}.
\end{equation}
\end{lemma}

\begin{proofsect}{Proof}
Let us for simplicity assume (e.g., by redefining~$a_N$) that~$b=0$. The strong Markov property bounds the probability under the sum by
\begin{equation}
\label{E:6.6}
\sum_{m\ge0}P^x\bigl(H_z<H_\varrho\bigr)
P^z\Bigl(\wt L^{D_N}_{H_\varrho}(z)\ge mG^{D_N}(z,z)\Bigr) P^\varrho\Bigl(\wh L^{D_N}_{t_N}(z)\ge a_N- (m+1)G^{D_N}(z,z)\Bigr).
\end{equation}
 We start by estimating the second term.
Denoting  $p:=P^z(\hat H_z<H_\varrho)$ where $\hat H_z$ is the first return time to~$z$, we have $\wt L^{D_N}_{H_\varrho}(z)\laweq \frac{1}{4}\sum_{i=1}^N\tau_i$ for $N:=\text{Geometric}(p)$ and $\tau_1,\tau_2,\dots$ i.i.d.\ Exponential(1) independent of~$N$. For any $q\in(0,1)$, the Chernoff bound gives
\begin{equation}
P\biggl(\sum_{i=1}^N \tau_i>r\biggr)\underset{0\le s<1-p}\le \texte^{-sr}\frac{1-p}{1-p-s}
\underset{s:=q(1-p)}\le \frac1{1-q}\texte^{-rq(1-p)}.
\end{equation}
As $1-p = P^z(\hat H_z>H_\varrho)=\frac1{4G^{D_N}(z,z)}$, we thus get
\begin{equation}
\label{E:6.9}
P^z\Bigl(\wt L^{D_N}_{H_\varrho}(z)\ge m G^{D_N}(z,z)\Bigr)\le \frac1{1-q}\texte^{-mq},\quad m\ge0.
\end{equation}
for all~$q\in(0,1)$.

Using \eqref{E:6.9} in conjunction with the uniform estimate $G^{D_N}(z,z)\le g\log N+c$, we dominate the part of the sum in \eqref{E:6.6} for~$m$ satisfying
$(m+2) G^{D_N}(z,z)\ge a_N - t_N$ by a quantity of order $N^{-2q[(\sqrt\theta+\lambda)^2-\theta+o(1)]}$. Recalling that~$W_N=N^{2(1-\lambda^2)+o(1)}$, this is $o(W_N N^{-2}/\log N)$ when~$1-q>0$ is so small that $q[(\sqrt\theta+\lambda)^2-\theta]>\lambda^2$.  

In the complementary regime, we have $a_N-(m+2) G^{D_N}(z,z)>t_N$ which permits us to estimate the last term on the right of \eqref{E:6.6} via \cite[Lemma~4.1]{AB} with the choices $a:=a_N$, $t:=t_N$ and~$b:=(m+2) G^{D_N}(z,z)$ to get
\begin{multline}
\label{E:6.10}
P^z\Bigl(\wt L^{D_N}_{H_\varrho}(z)\ge mG^{D_N}(z,z)\Bigr) P^\varrho\Bigl(\wh L^{D_N}_{t_N}(z)\ge a_N- (m+1)G^{D_N}(z,z)\Bigr)
\\
\le 
\frac1{1-q}\,\frac{\sqrt{G^{D_N} (z,z)}}{\sqrt{2a_N - 2(m+1)G^{D_N} (z,z)} - \sqrt{2t_N}} 
\frac{\sqrt{\log N}}{N^2} W_N\,\text{e}^{-qm + (m+1)\frac{\sqrt{2a_N}-\sqrt{2t_N}}{\sqrt{2a_N}}}
\end{multline}
As $\frac{\sqrt{2a_N}-\sqrt{2t_N}}{\sqrt{2a_N}}\to\frac{\lambda}{\sqrt\theta+\lambda}$ as~$N\to\infty$, we  choose~$q\in(\frac{\lambda}{\sqrt\theta+\lambda},1)$ and proceed as follows: For $(m+1) G^{D_N} (z, z) > \frac{1}{2} (a_N - t_N)$, the prefactor is order $\sqrt{\log N}\, W_N/N^2$ but, thanks to the uniform upper bound on~$G^{D_N}(z,z)$, the sum of the exponential terms decays polynomially with~$N$. For~$m$ with $(m+1) G^{D_N} (z, z) \leq \frac{1}{2} (a_N - t_N)$, the prefactor is order $W_N/N^2$ and the sum of the exponentials is bounded.

Combining the above estimates, the sum in \eqref{E:6.5uie} is bounded by a quantity of order
\begin{equation}
\label{E:6.10iu}
o\Bigl(\frac{W_N}{\log N}\Bigr)+
\frac{W_N}{N^2}\sum_{z\in D_N}P^x\bigl(H_z<H_\varrho\bigr).
\end{equation}
Interpreting~$H_\varrho$ as the first exit time of the simple random walk on~$\Z^2$ from~$D_N$, the sum on the right is non-decreasing in~$D_N$. We may thus assume that~$D_N$ is a box of side-length~$2^n$, for~$n=\log_2 N+O(1)$, centered at~$x$. For the probability under the sum we then get, for each~$k=0,\dots,n-1$ and some constant~$c>0$,
\begin{equation}
P^x (H_z < H_{\rho}) = \frac{G^{D_N} (x, z)}{G^{D_N} (z, z)}
\le c\frac {n-k}n,\qquad 2^{k}<|x-z|\le 2^{k+1}.
\end{equation}
The sum in \eqref{E:6.10iu} is thus at most of order $1+\sum_{k=0}^n\frac{n-k}n\,2^{2k}$ which is of order $N^2/\log N$. The claim follows.
\end{proofsect}

We are now ready to give:

\begin{proofsect}{Proof of Theorem~\ref{thm-6.1}, thick points}
Consider a non-negative~$f\in C_\cc(\overline D\times(\R\cup\{+\infty\})$ that is non-decreasing in the second variable and supported in $\overline D\times[b,\infty)$ for some~$b\in\R$.
Note that~$\{H_\varrho<\infty\}$ is a full probability event under~$P^x$. Decomposing the support of~$\zeta^D_N$ according to whether the point was hit before hitting the boundary vertex or not, the monotonicity of~$t\mapsto\wt L^{D_N}_t$ and the assumed monotonicity of~$f$ yield
\begin{multline}
\label{E:6.11}
\quad
\langle\wt\zeta^D_N,f\rangle\le\langle\wt\zeta^D_N,f\rangle\circ \theta_{H_\varrho}
\\
+\frac{\Vert f\Vert_\infty}{W_N}\sum_{z\in D_N}1_{\{H_z<H_\varrho\}}\,1_{\{\wt L^{D_N}_{H_\varrho}(z)+(\wt L^{D_N}\circ\theta_{H_\varrho})_{t_N\deg(D_N)}(z)\ge a_N+ b \sqrt{2 a_N} \}}.
\quad
\end{multline}
Fix a sequence~$b_N\to\infty$ such that~$b_N/t_N^{1/4}\to0$ and let~$\FF_N$ be the event that the inequalities in \eqref{E:5.6i} hold. Fix any~$m>0$ and~$\epsilon>0$. Let~$\GG_N$ be the event that the second term on the right of \eqref{E:6.11} is less than~$\epsilon$. Then
\begin{equation}
\label{E:6.12}
\begin{aligned}
E^x\bigl(\texte^{-\langle\wt\zeta^D_N,f\rangle}\bigr)
&\ge 
E^x\Bigl(\texte^{-\langle\wt\zeta^D_N,f\rangle}1_{\theta_{H_\varrho}^{-1}(\FF_N\cap\{T_N\ge-m\})}\Bigr)
\\
&\ge \texte^{-\epsilon}
E^x\Bigl(\texte^{-\langle\wt\zeta^D_N,f\rangle\circ\theta_{H_\varrho}}1_{\theta_{H_\varrho}^{-1}(\FF_N\cap\{T_N\ge-m\})}
\Bigr)
\\
&\qquad\qquad\qquad-P^x\Bigl(\GG_N^\cc\cap\theta_{H_\varrho}^{-1}(\FF_N\cap\{T_N\ge-m\})\Bigr).
\end{aligned}
\end{equation}
As~$P^x(H_\varrho<\infty)=1$, the strong Markov property gives
\begin{equation}
\begin{aligned}
E^x\Bigl(\texte^{-\langle\wt\zeta^D_N,f\rangle\circ\theta_{H_\varrho}}1_{\theta_{H_\varrho}^{-1}(\FF_N\cap\{T_N\ge-m\})}
\Bigr)
&=E^\varrho\Bigl(\texte^{-\langle\wt\zeta^D_N,f\rangle}1_{\FF_N\cap\{T_N\ge-m\}}
\Bigr)
\\
&\ge E^\varrho\bigl(\texte^{-\langle\wt\zeta^D_N,f\rangle}\bigr)-P^\varrho\bigl((\FF_N\cap\{T_N\ge-m\})^\cc\bigr).
\end{aligned}
\end{equation}
Proposition~\ref{P:tNbound} and the fact that~$\{T_N\colon N\ge1\}$ is tight now ensures that the probability on the right tends to zero in the limits~$N\to\infty$ and~$m\to\infty$.

Concerning the probability on the right of \eqref{E:6.12}, an inspection of \eqref{E:tNcirc} shows that, on $(\FF_N\cap\{T_N\ge-m\})\circ\theta_{H_\varrho}$, we have
\begin{equation}
(\wt L^{D_N}\circ\theta_{H_\varrho})_{t_N\deg(D_N)}(\cdot)
\le(\wh L^{D_N}\circ \theta_{H_\varrho})_{t_N+b_Nt_N^{1/4}+m\sqrt{2t_N}}(\cdot).
\end{equation}
By the Markov inequality, the probability in \eqref{E:6.12} is thus bounded by $\epsilon^{-1}\Vert f\Vert_\infty/W_N(t_N)$ times the sum in Lemma~\ref{lemma-6.3} albeit with~$t_N$ replaced by $t_N':=t_N+b_Nt_N^{1/4}+m\sqrt{2t_N}$. As~$W_N(t_N')/W_N(t_N)$ is bounded by an~$m$-dependent constant uniformly in~$N$, the probability in \eqref{E:6.12} is thus $O(1/\log N)$ uniformly in~$x\in D_N$. 
Taking~$N\to\infty$ followed by~$m\to\infty$ and~$\epsilon\downarrow0$ shows
\begin{equation}
\liminf_{N\to\infty} E^{x_N}\bigl(\texte^{-\langle\wt\zeta^D_N,f\rangle}\bigr)
\ge \lim_{N\to\infty}\,E^{\varrho}\bigl(\texte^{-\langle\wt\zeta^D_N,f\rangle}\bigr).
\end{equation} 
Combining with \eqref{E:6.2}, we then get the desired claim.
\end{proofsect}

The situation for the thin, light and avoided points is similar albeit simpler. Writing $\wt\xi^D_N$ for the corresponding continuous-time point measure (parametrized by the total time), as in Lemma~\ref{lemma-6.2}, the identity \eqref{E:6.1} gives us an easy one-way bound, where the test function~$f$ takes values in~$\overline D\times(\R\cup\{-\infty\})$ for the thin points, $\overline D\times[0,\infty)$ for the light points and~$\overline D$ for the avoided points:

\begin{lemma}
\label{lemma-6.4}
Under the conditions of Proposition~\ref{thm-thin-cont}, \ref{thm-light-cont} and \ref{thm-avoid-cont},  for any any $x_N\in D_N$ and any continuous, compactly-supported, non-negative test function $f$ on the corresponding domain that, for the thin and light points, is non-increasing in the second variable,
\begin{equation}
\liminf_{N\to\infty} E^{x_N}\bigl(\texte^{-\langle\wt\xi^D_N,f\rangle}\bigr)
\ge \lim_{N\to\infty}\,E^{\varrho}\bigl(\texte^{-\langle\wt\xi^D_N,f\rangle}\bigr).
\end{equation}
\end{lemma}

\begin{proofsect}{Proof}
Using \eqref{E:6.1}, on $\{H_\varrho< m\deg(D_N)\}$ we get
\begin{equation}
\bigl\langle\wt\xi^{D_N}_N(t_N),f\bigr\rangle\le
\texte^{c \frac{m}{\log N}}
\bigl\langle\wt\xi^{D_N}_N(t_N-m),f\bigr\rangle\circ\theta_{H_\varrho},
\end{equation}
where we now rely on the fact that~$t\mapsto W_N(t)$, resp., $t\mapsto\wh W_N(t)$ are non-increasing for~$t$ near~$t_N$. The inequalities \eqref{E:6.4} then become
\begin{equation}
\begin{aligned}
E^{x_N}\bigl(\texte^{-\langle\wt\xi^D_N(t_N),f\rangle}\bigr)
&\ge E^{x_N}\bigl(\texte^{-\langle\wt\xi^D_N(t_N),f\rangle}1_{\{H_\varrho<m\deg(D_N)\}}\bigr) \\
&\ge E^{\varrho}\bigl(\texte^{-\texte^{c \frac{m}{\log N}}\langle\wt\xi^D_N(t_N-m),f\rangle}\bigr)-P^{x_N}\bigl(H_\varrho\ge m\deg(D_N)\bigr).
\end{aligned}
\end{equation}
The claim now follows by taking~$N\to\infty$ followed by~$m\to\infty$. 
\end{proofsect}

In replacement of Lemma~\ref{lemma-6.3}, we then need:

\begin{lemma}
\label{lemma-6.5}
Under the conditions of Proposition~\ref{thm-thin-cont}, for each~$b\in\R$ there is~$c>0$ such that for all $N\ge1$ and all~$x\in D_N$,
\begin{equation}
\label{E:6.20uiu}
\sum_{z\in D_N}
P^x\biggl((\wh L^{D_N}\circ\theta_{H_\varrho})_{t_N}(z)\le a_N+b\log N,\, H_z<H_\varrho\biggr)
\le c\frac {W_N}{\log N}.
\end{equation}
Under the conditions of Propositions~\ref{thm-light-cont} and \ref{thm-avoid-cont} the same holds with $a_N+b\log N$ replaced by~$b\ge0$ (including, for the avoided points,~$b=0$) and~$W_N$ replaced by~$\wh W_N$.
\end{lemma}

\begin{proofsect}{Proof}
The Strong Markov property and the estimates from~\cite[Corollary~4.8]{AB} bound the probability in \eqref{E:6.20uiu} by $P^x(H_z<H_\varrho)$ times
\begin{equation}
P^\varrho\bigl(\wh L^{D_N}_{t_N}(z)\le a_N + b \log N \bigr)\le c\frac{W_N}{N^2}
\end{equation}
and so the quantity in \eqref{E:6.20uiu} is at most order $W_N N^{-2}\sum_{z \in D_N}P^x(H_z<H_\varrho)$. The argument then concludes as in the proof of Lemma~\ref{lemma-6.3}.  For the light and avoided points, we instead invoke \cite[Corollary~4.6]{AB} and proceed analogously. 
\end{proofsect}

With this we get:

\begin{proofsect}{Proof of Theorem~\ref{thm-6.1}, thin, light and avoided points}
We proceed similarly as for the thick points. First, writing~$\wt a_N:=a_N+b\log N$ for the thin points and~$\wt a_N:=b$ for the light and (with~$b:=0$) avoided points, given a continuous, compactly-supported~$f$ that is non-increasing in the second variable, in all three cases of interest we have
\begin{equation}
\label{E:6.20}
\langle\wt\xi^D_N,f\rangle\ge\langle\wt\xi^D_N,f\rangle\circ \theta_{H_\varrho}
-\frac{\Vert f\Vert_\infty}{W_N}\sum_{z\in D_N}1_{\{H_z<H_\varrho\}}\,1_{\{(\wt L^{D_N}\circ\theta_{H_\varrho})_{t_N\deg(D_N)-H_\varrho}(z)\le \wt a_N\}}.
\end{equation}
Let $\mathcal{F}_N$ be the event from \eqref{E:5.6i} with $t_N$ replaced by $t_N - m$. 
Abusing our earlier notation, given~$\epsilon>0$, let $\GG_N$ be the event that the second term (without the minus sign) is at most~$\epsilon$. From \eqref{E:6.20}, we then get
\begin{equation}
\label{E:6.21}
\begin{aligned}
E^x\bigl(\texte^{-\langle\wt\xi^D_N,f\rangle}&\bigr)-P^x\bigl(H_\varrho\ge m\deg(D_N)\bigr)
-P^\varrho\bigl((\FF_N\cap  \{T_N (t_N - m) \leq m \} )^\cc\bigr)
\\
&\le
E^x\Bigl(\texte^{-\langle\wt\xi^D_N,f\rangle}1_{\{H_\varrho<m\deg(D_N)\}}1_{\theta_{H_\varrho}^{-1}(\FF_N\cap  \{T_N (t_N - m) \leq m \} )}\Bigr)
\\
&\le\texte^\epsilon E^\varrho\bigl(\texte^{-\langle\wt\xi^D_N,f\rangle}\bigr)
\\
&\qquad\quad+P^x\Bigl(\GG_N^\cc\cap \{H_\varrho<m\deg(D_N)\}\cap\theta_{H_\varrho}^{-1}(\FF_N\cap  \{T_N (t_N - m) \leq m \} )\Bigr).
\end{aligned}
\end{equation}
Thanks to the Central Limit Theorem, the tightness of~$\{T_N\colon N\ge1\}$ and Proposition~\ref{P:tNbound}, the two probabilities on the left-hand side of \eqref{E:6.21} tend to zero in the limits~$N\to\infty$ and~$m\to\infty$, uniformly in~$x\in D_N$. For the probability on the right we observe that, on $\{H_\varrho<m\deg(D_N)\}\cap\theta_{H_\varrho}^{-1}(\FF_N\cap 
 \{T_N (t_N - m) \leq m \} )$, we have
\begin{equation}
(\wt L^{D_N}\circ\theta_{H_\varrho})_{t_N\deg(D_N)-H_\varrho}(\cdot)
\ge(\wh L^{D_N}\circ\theta_{H_\varrho})_{t_N'}(\cdot)
\end{equation}
for~$t_N':=t_N-m-b_N t_N^{1/4}-m\sqrt{2t_N}$. Lemma~\ref{lemma-6.5} and the Markov inequality then bound the probability by an $m$-dependent constant times~$1/\log N$, uniformly in~$x\in D_N$.
Combining these observations we thus get
\begin{equation}
\limsup_{N\to\infty} E^{x_N}\bigl(\texte^{-\langle\wt\xi^D_N,f\rangle}\bigr)
\le \lim_{N\to\infty}\,E^{\varrho}\bigl(\texte^{-\langle\wt\xi^D_N,f\rangle}\bigr).
\end{equation}
In conjunction with Lemma~\ref{lemma-6.4} this proves the claim.
\end{proofsect}

\section{Discrete time conclusions}
\label{sec7}\noindent
We will now move to the proof of our main results except those on the local structure which are deferred to Section~\ref{sec8}. Considering, for a moment, a random walk on a general finite, connected graph on $V\cup\{\varrho\}$, recall that the discrete-time local time $L^V_t$ is parametrized by the total number of steps in units of~$\deg(V)=\sum_{u\in V\cup\{\varrho\}}\deg(u)$ while its continuous-time counterpart $\wt L^V_t$ is parametrized by the total time. Both of these are naturally realized on the same probability space through the definition \eqref{E:tilde-X} of~$\wt X$ via the discrete-time walk~$X$ and an independent (rate-1) Poisson point process~$\wt N(t)$. A key technical tool in what follows is the following lemma:

\begin{lemma}
\label{lemma-7.1}
 There is a family of i.i.d.\ exponentials $\{\tau_j(v)\colon j\ge1, v\in V\}$
with parameter $1$
independent of~$X$ (but not of~$\wt N$) such that
\begin{equation}
\label{E:8.1a}
\wt L_t^V(v) = \frac1{\deg(V)}\sum_{j\ge1}\tau_j(v)1_{\{j\le\deg(V) L^V_{\wt N(t)/\deg(V)}(v)\}},\quad v\ne \wt X_t,
\end{equation}
holds $P^x$-a.s.\ for each $t\ge0$ and each~$x\in V\cup\{\varrho\}$.
\end{lemma}

\begin{proofsect}{Proof}
This is a consequence of the standard representation of the wait times of~$\wt X$ by independent exponentials. (In this representation, the process~$\wt N$ is a function of the exponentials and~$X$, albeit independent of~$X$.) Note that the equality \eqref{E:8.1a} fails at~$\wt X_t$ because the walk is ``in-between'' jumps there.
\end{proofsect}

Moving back to the random walk on~$D_N\cup\{\varrho\}$, this readily yields:

\begin{lemma}
\label{lemma-7.2}
For each~$x\in D_N$, abbreviate
\begin{equation}
\label{E:7.2}
\FF_{N}(x):=\biggl\{\wt L^{D_N}_{(t_N-1)\deg(D_N)}(x)\le\frac14\sum_{j\ge1}\tau_j(x)1_{\{j\le 4L^{D_N}_{t_N}(x)\}}\le\wt L^{D_N}_{(t_N+1)\deg(D_N)}(x)\biggr\}.
\end{equation}
Then for any~$x_N\in D_N$,
\begin{equation}
P^{x_N}\Bigl(\,\sum_{x\in D_N}1_{\FF_N(x)^\cc}>2\Bigr)\,\underset{N\to\infty}\longrightarrow\,0.
\end{equation}
\end{lemma}

\begin{proofsect}{Proof}
The Central Limit Theorem ensures that~$(\wt N(t)-t)/\sqrt{t}$ tends in law to a standard normal as~$t\to\infty$. As~$t_N=o(\deg(D_N))$, the inequalities 
\begin{equation}
\label{E:7.4}
\frac{\wt N((t_N-1)\deg(D_N))}
{\deg(D_N)}\le t_N\le \frac{\wt N((t_N+1)\deg(D_N))}
{\deg(D_N)}
\end{equation}
are satisfied with probability tending to one as~$N\to\infty$. Once \eqref{E:7.4} is in force, the monotonicity of $t\mapsto L^{D_N}_t$ and \eqref{E:8.1a} show that the event~$\FF_N(x)$ occurs at all~$x\in D_N$ except perhaps at the position of~$\wt X$ at times $(t_N\pm 1)\deg(D_N)$.
\end{proofsect}

With these observations in hand, we are now ready to finally present the proofs of our main theorems. The easiest case is that of avoided points:

\begin{proofsect}{Proof of Theorem~\ref{thm-avoid}}
Note that, whenever~$\FF_N(x)$ occurs, $\wt L^{D_N}_{(t_N+1)\deg(D_N)}(x)=0$ forces $L^{D_N}_{t_N}(x)=0$ (a.s.), which in turn forces $\wt L^{D_N}_{(t_N-1)\deg(D_N)}(x)=0$. For any~$f\in C_\cc(\overline D)$ with~$f\ge0$, on the event $\sum_{x\in D_N}1_{\FF_N(x)^\cc}\le 2$ we thus have
\begin{equation}
\label{E:7.5}
\begin{aligned}
\frac{\wh W_N(t_N+1)}{\wh W(t_N)}\bigl\langle\wt\kappa^D_N(t_N+1),f\bigr\rangle-\frac2{\wh W_N} \Vert f\Vert_\infty
&\le \bigl\langle\kappa^D_N,f\bigr\rangle
\\
&\le 
\frac{\wh W_N(t_N-1)}{\wh W(t_N)}\bigl\langle\wt\kappa^D_N(t_N-1),f\bigr\rangle
+\frac2{\wh W_N} \Vert f\Vert_\infty.
\end{aligned}
\end{equation}
As $\{t_N\pm1\}_{N\ge1}$ have the same leading-order asymptotic as~$\{t_N\}_{N\ge1}$, the random variables $\langle\wt\kappa^D_N(t_N\pm1),f\rangle$ have the same weak limit as $\langle\wt\kappa^D_N,f\rangle$. Since~$\wh W_N\to\infty$ and also
\begin{equation}
\frac{\wh W_N(t_N\pm1)}{\wh W_N(t_N)}\,\underset{N\to\infty}\longrightarrow\,1,
\end{equation}
the claim follows from Lemma~\ref{lemma-7.2}, Proposition~\ref{thm-avoid-cont} and Theorem~\ref{thm-6.1}.
\end{proofsect}

Next we tackle the light points:

\begin{proofsect}{Proof of Theorem~\ref{thm-light}}
Denote
\begin{equation}
\label{E:7.7}
\overline L^{D_N}_{t_N}(x):=\frac14\sum_{j\ge1}\tau_j(x)1_{\{j\le 4L^{D_N}_{t_N}(x)\}}
\end{equation}
and consider the auxiliary point measure
\begin{equation}
\overline\vartheta^D_N:=\frac1{\wh W_N}\sum_{x\in D_N}\delta_{x/N}\otimes\delta_{\overline L^{D_N}_{t_N}(x)}.
\end{equation}
Thanks to Lemma~\ref{lemma-7.2}, on the event $\sum_{x\in D_N}1_{\FF_N(x)^\cc}\le 2$, the inequality \eqref{E:7.5} holds for any non-negative $f\in C_\cc(\overline D\times[0,\infty))$ that is non-increasing in the second variable and with~$\wt\kappa^D_N$, resp., $\kappa^D_N$ replaced by $\wt\vartheta^D_N$, resp., $\overline\vartheta^D_N$. As, by Proposition~\ref{thm-light-cont} and Theorem~\ref{thm-6.1}, $\wt\vartheta^D_N$ tends in law to the measure~$\wt\vartheta^D$ on the right of \eqref{E:2.31ii}, we have
\begin{equation}
\langle\overline\vartheta^D_N,f\rangle\,\,\underset{N\to\infty}\Lawarrow\,\,\langle\wt\vartheta^D, f\rangle
\end{equation}
for any non-negative $f\in C_\cc(\overline D\times[0,\infty))$.

Next we observe that, by that fact that for any~$\epsilon>0$ and any random variable~$Y$ taking values in $[0,\epsilon]$,
\begin{equation}
\label{E:7.10}
\exp\{-E(Y)\}\le E(\texte^{-Y})\le \exp\{-\texte^{-\epsilon} E(Y)\},
\end{equation}
the fact that the random variables $\{\tau_j(x)\colon j\ge1,\,x\in D_N\}$ are independent of the random walk   and independent for different~$x\in D_N$ implies
\begin{equation}
E^{x_N}\bigl(\texte^{-E(\langle\overline\vartheta^D_N,f\rangle|\sigma(X))}\bigr)
\le E^{x_N}\bigl(\texte^{-\langle\overline\vartheta^D_N,f\rangle}\bigr)
\le E^{x_N}\bigl(\texte^{-\texte^{-\Vert f\Vert_\infty/\wh W_N}E(\langle\overline\vartheta^D_N,f\rangle|\sigma(X))}\bigr)
\end{equation}
(see~\cite[Lemma 3.12]{BL4}), 
where the conditional expectation is meaningful because $\langle\overline\vartheta^D_N,f\rangle$ is a finite random variable.
Defining $f^{\ast\mathfrak e}\colon\overline D\times[0,\infty)\to\R$ by
\begin{equation}
f^{\ast\mathfrak e}(x,\ell):= E\biggl[\,f\Bigl(x,\frac14\sum_{j=1}^{\lfloor 4\ell\rfloor } \tau_j\Bigr)\biggr],
\end{equation}
where~$\{\tau_j\colon j\ge1\}$ are i.i.d.\ Exponential(1), we have
\begin{equation}
E^{x_N}\bigl(\langle\overline\vartheta^D_N,f\rangle\,\big|\,\sigma(X)\bigr) =\langle\vartheta^D_N,f^{\ast\mathfrak e}\rangle.
\end{equation}
Hence we get (under the laws $\{P^{x_N}\colon N\ge1\}$),
\begin{equation}
\label{E:7.14}
\langle\vartheta^D_N,f^{\ast\mathfrak e}\rangle\,\,\underset{N\to\infty}\Lawarrow\,\,\langle\wt\vartheta^D ,f\rangle
\end{equation}
for any $f\in C_\cc(\overline D\times[0,\infty))$.

We now claim that $\{\vartheta^D_N\colon N\ge1\}$ is tight. For this we pick $M\in\N$, denote $f_M(x,h):=1_{[0,M]}(h)$ and observe that, for all~$n\in \N_0:=\{0,1,2,\dots\}$, we get
\begin{equation}
\label{E:7.15}
f_M^{\ast\mathfrak e}(x, n/4)=P\Bigl(\frac14\sum_{j=1}^{ n}\tau_j\le M\Bigr).
\end{equation}
Markov's inequality then shows $f_{2M}^{\ast\mathfrak e}(x, n/4)\ge\frac121_{[0,M]}( n/4)$ and, therefore,
\begin{equation}
\vartheta^D_N\bigl(\overline D\times[0,M]\bigr)\le2\langle\vartheta^D_N,f_{2M}^{\ast\mathfrak e}\rangle.
\end{equation}
The existence of the limit \eqref{E:7.14} then implies tightness of $\{\vartheta^D_N(\overline D\times[0,M])\colon N\ge1\}$ for all~$M>0$, and thus tightness of $\{\vartheta^D_N\colon N\ge1\}$ as well.

The tightness of $\{\vartheta^D_N\colon N\ge1\}$ permits us to extract a weak subsequential limit~$\vartheta^D$ along a (strictly) increasing sequence~$\{N_k\colon k\ge1\}$ of naturals. This entails the convergence $\langle\vartheta^D_{N_k},f\rangle\Lawarrow\langle\vartheta^D,f\rangle$ for every~$f\in C_\cc(\overline D\times[0,\infty))$. We claim that we even have
\begin{equation}
\label{E:7.17}
\langle\vartheta^D_{N_k},f^{\ast\mathfrak e}\rangle\,\,\underset{k\to\infty}\Lawarrow\,\,\langle\vartheta^D,f^{\ast\mathfrak e}\rangle
\end{equation}
for every $f\in C_\cc(\overline D\times[0,\infty))$. (This is not automatic because $f^{\ast\mathfrak e}$ is not compactly supported in general.) First we note that straightforward comparisons with the Lebesgue measure show, for each~$M>0$,
\begin{equation}
\label{E:7.18}
\lim_{n\to\infty}\frac{P\Bigl(\frac14\sum_{j=1}^{4n}\tau_j\le M\Bigr)}{P\Bigl(\frac14\sum_{j=1}^{4n}\tau_j\le 2M\Bigr)}=0.
\end{equation}
Writing~$\epsilon_n$ for the ratio of the two probabilities, for~$f$ supported in~$\overline D\times[0,M]$ we have $|f^{\ast\mathfrak e}|\le \Vert f\Vert_{\infty} f^{\ast\mathfrak e}_M$ and so, by \eqref{E:7.15},
\begin{equation}
\bigl|f^{\ast\mathfrak e}(x,n/4)\bigr|\le \epsilon_n \Vert f\Vert_{\infty} \,f_{2M}^{\ast\mathfrak e}(x,n),\quad n\in\N_0.
\end{equation}
It follows that the part of the integral $\langle\vartheta^D_N,f^{\ast\mathfrak e}\rangle$ corresponding to the second coordinate in excess of~$n$ is at most~$\epsilon_n \Vert f\Vert_{\infty}$ times $\langle\vartheta^D_N,f^{\ast\mathfrak e}_{2M}\rangle$, which is tight by \eqref{E:7.14}. We can thus approximate $f^{\ast\mathfrak e}$ by a function supported in~$\overline D\times[0,n]$ and pass to the limit $N\to\infty$ followed by~$n\to\infty$. This gives \eqref{E:7.17} as desired.

Combining \eqref{E:7.14} with \eqref{E:7.17} we arrive at the convolution identity
\begin{equation}
\label{E:7.20}
\langle\vartheta^D,f^{\ast\mathfrak e}\rangle\laweq\langle\wt\vartheta^D,f\rangle.
\end{equation}
We have proved this (including the absolute convergence of the integral on the left-hand side) for $f\in C_\cc(\overline D\times[0,\infty))$ but the Monotone Convergence Theorem along with the fact that the second coordinate of~$\wt\vartheta^D$ has subexponentially growing density extends this to all~$f\in C(\overline D\times[0,\infty))$ such that $|f(x,h)|\le c\texte^{-\epsilon h}$ for some~$\epsilon,c>0$. This permits us to consider functions of the form $g_s(x,h):=\tilde f(x)\texte^{-sh}$ for~$s>0$ and~$\tilde f\in C(\overline D)$, for which
\begin{equation}
g_s^{\ast\mathfrak e}(x,n/4)=\tilde f(x)(1+s/4)^{-n},\quad n\in\N_0.
\end{equation}
Since $\vartheta^D$ is supported on $\overline D\times\frac14\N_0$, it makes sense to denote
\begin{equation}
\vartheta^{D,n}(A):=\vartheta^D\bigl(A\times\{n/4\}\bigr).
\end{equation}
The identity \eqref{E:7.20} then becomes
\begin{equation}
\label{E:7.23}
\sum_{n\ge0}\langle\vartheta^{D,n},\tilde f\rangle(1+s/4)^{-n}\laweq\langle\wt\vartheta^D,g_s\rangle.
\end{equation}
Assuming $\tilde f>0$, the explicit form of the right-hand side shows that $\langle\wt\vartheta^D,g_s\rangle/\langle\wt\vartheta^D,g_1\rangle$ is well-defined and equal to a non-random quantity --- namely, the ratio of two Laplace transforms of~$\tilde\mu$. This turns \eqref{E:7.23} into the pointwise identity
\begin{equation}
\label{E:7.24}
\sum_{n\ge0}\langle\vartheta^{D,n},\tilde f\rangle(1+s/4)^{-n}=
\frac{\int\tilde\mu(\textd h)\texte^{-sh}}
{\int\tilde\mu(\textd h)\texte^{-h}}\biggl(\,\sum_{n\ge0}\langle\vartheta^{D,n},\tilde f\rangle(5/4)^{-n}\biggr)
\end{equation}
valid, a.s., for each~$s>0$ and (by elementary extensions) all~$\tilde f\in C(\overline D)$. Thanks to the monotonicity of both sides in~$s$ and almost-sure continuity in~$\tilde f$ of both sides with respect to the supremum norm, the identity actually holds a.s.\ for all~$s>0$ and all~$\tilde f\in C(\overline D)$ simultaneously. 

With \eqref{E:7.24} in hand, we are more or less done. Indeed, as the left-hand side is a generating function of the sequence $\{\langle\vartheta^{D,n},\tilde f\rangle\}_{n\ge0}$, which determines the sequence uniquely, all $\langle\vartheta^{D,n},\tilde f\rangle$ must be the same deterministic multiple of the quantity in the large parentheses on the right-hand side. This shows that $\vartheta^D$ must be as on the right-hand side of \eqref{E:2.22ii} for some~$\mu$ of the form $\mu=\sum_{n\ge0}q_n\delta_{n/4}$ where~$\{q_n\}_{n\ge0}$ is uniquely determined by
\begin{equation}
\label{E:7.18a}
\sum_{n\ge0}q_n(1+s/4)^{-n}=\int_0^\infty\tilde\mu(\textd h)\texte^{-sh},\quad s>0.
\end{equation}
The Laplace transform of~$\tilde\mu$ was calculated in the proof of Proposition~\ref{thm-light-bv}. All subsequential limits of $\{\vartheta^D_N\colon N\ge1\}$ are thus equal in law and so convergence holds.
\end{proofsect}

Moving to the thick points, we first need a version of \eqref{E:7.18}:

\begin{lemma}
\label{lemma-7.3}
For $\{\tau_j\colon j\ge1\}$ be i.i.d.\ Exponential(1), all $k\in\N$ and all reals $s\ge t\ge0$,
\begin{equation}
\frac{P\Bigl(\sum_{j=1}^k(\tau_j-1)\ge s+t\Bigr)}{P\Bigl(\sum_{j=1}^k(\tau_j-1)\ge s\Bigr)}
\le\texte^{-\frac{st}{k+s+t}}.
\end{equation}
\end{lemma}

\begin{proofsect}{Proof}
Since $\sum_{j=1}^k\tau_j$ has density $\frac1{(k-1)!}x^{k-1}\texte^{-x}$, the change of variables $y:=x+t$ gives
\begin{equation}
\begin{aligned}
P\biggl(\,\sum_{j=1}^k(\tau_j-1)\ge s\biggr)
&=\frac1{(k-1)!}\int_{x\ge k+s}\textd x\, x^{k-1}\texte^{-x}
\\
&=\texte^t\,\frac1{(k-1)!}\int_{y\ge k+s+t}\textd y\, (y-t)^{k-1}\texte^{-y}
\\
&
\ge \texte^t\Bigl(1-\frac{t}{k+s+t}\Bigr)^k\,P\biggl(\,\sum_{j=1}^k(\tau_j-1)\ge s+t\biggr).
\end{aligned}
\end{equation}
Using that $s\ge t$, the prefactor can be written as the exponential of
\begin{equation}
\begin{aligned}
t+k\log\Bigl(1-\frac{t}{k+s+t}\Bigr)&=t-k\sum_{n\ge1}\frac1n\frac{t^n}{(k+s+t)^n}
\\
&\ge t-\frac{kt}{k+s+t}-\frac12\frac{k t^2}{(k+s+t)^2}\sum_{n\ge0}2^{-n}.
\end{aligned}
\end{equation}
Noting that right-hand side is no less than $\frac{st}{k+s+t}$, we get the claim.
\end{proofsect}

A convolution identity that inevitably shows up in the proof also requires: 

\begin{lemma}
\label{lemma-7.4}
Suppose~$\nu$ is a Borel measure on~$\R$ such that, for some~$\beta\in\R$ and some~$\sigma^2>0$ and all~$f\in C_\cc(\R)$,
\begin{equation}
\label{E:7.29}
\int_\R \nu(\textd h) \,E\bigl[ \,f(h+\NN(0,\sigma^2))\bigr]  = \int_\R\textd h\,\,\texte^{\beta h}\,f(h)
\end{equation}
Then
\begin{equation}
\nu(\textd h)=\texte^{-\frac12\beta^2\sigma^2+\beta h}\textd h.
\end{equation}
\end{lemma}

\begin{proofsect}{Proof}
Consider the measure $\wt\nu(\textd h):=\texte^{-\beta h+\frac12\beta^2\sigma^2}\nu(\textd h)$.
Absorbing the exponential term on the right of \eqref{E:7.29} into the test function, a calculation shows 
\begin{equation}
\int_{\R\times\R}\wt\nu(\textd h)\otimes\frac{\textd x}{\sqrt{2\pi\sigma^2}}\,\,\texte^{-\frac{(x-h+\beta\sigma^2)^2}{2\sigma^2}}\,f(x)=\int_\R\textd h\,f(h)
\end{equation}
for all~$f\in C_\cc(\R)$. As $C_\cc(\R)$ generates all Borel functions in~$\R$, we get
\begin{equation}
\frac1{\sqrt{2\pi\sigma^2}}\int_\R\wt\nu(\textd h)\,\texte^{-\frac{(x-h+\beta\sigma^2)^2}{2\sigma^2}} = 1,\quad x\in\R.
\end{equation}
This can be interpreted by saying that $\wh\nu(\textd h):=\,\frac1{\sqrt{2\pi\sigma^2}}\texte^{-\frac{(h-\beta\sigma^2)^2}{2\sigma^2}}\wt\nu(\textd h)$ is a measure such that
\begin{equation}
\int_\R\wh\nu(\textd h)\texte^{-xh} = \texte^{-x\beta\sigma^2+x^2\sigma^2/2},\quad x\in\R.
\end{equation}
The right-hand side is the Laplace transform of $\NN(\beta\sigma^2,\sigma^2)$ and so, since the Laplace transform of a measure, if exists, determines the measure uniquely, $\wh\nu$ is the law of $\NN(\beta\sigma^2,\sigma^2)$. Hence~$\wt\nu$ is the Lebesgue measure, thus proving the claim.
\end{proofsect}

\begin{proofsect}{Proof of Theorem~\ref{thm-thick}}
The proof starts by adapting the argument leading to \eqref{E:7.14}. Indeed, working again in the coupling of the random walk~$X$ and the i.i.d.\ exponentials $\{\tau_j(x)\colon x\in D_N,\,j\ge1\}$, let 
\begin{equation}
\overline\zeta^D_N:=\frac1{W_N}\sum_{x\in D_N}\delta_{x/N}\otimes\delta_{(\overline L_{t_N}^{D_N}(x)-a_N)/\sqrt{2a_N}}\,,
\end{equation}
where $\overline L_{t_N}^{D_N}(x)$ is the quantity from \eqref{E:7.7}.
Lemmas~\ref{lemma-7.1}-\ref{lemma-7.2} along with
Proposition~\ref{thm-thick-cont}, Theorem~\ref{thm-6.1} and~\eqref{E:7.10} then show
\begin{equation}
E^{x_N}\bigl(\langle\overline\zeta^D_N,f\rangle\,\big|\,\sigma(X)\bigr)\,\,\underset{N\to\infty}\Lawarrow\,\,\langle\wt\zeta^D, f\rangle
\end{equation}
for every~$f\in C_\cc(\overline D\times(\R\cup\{+\infty\}))$, where $\widetilde{\zeta}^D$ is the measure on the right of (\ref{E:1.21b}).
 Writing $\{\tau_j\colon j\ge1\}$ for generic i.i.d.\ exponentials with parameter 1 and denoting, with some abuse of earlier notation,
\begin{equation}
f^{N,\ast\mathfrak e}(x,h):= E\biggl[\,f\Bigl(x,h+\frac1{4\sqrt{2a_N}}\sum_{j\ge1} (\tau_j-1)1_{\{j\le 4a_N+4h\sqrt{2a_N}\}}\Bigr)\biggr],
\end{equation}
the fact that $L^{D_N}_{t_N}$ takes values in $\frac14\N_0$ then shows
\begin{equation}
E^{x_N}\bigl(\langle\overline\zeta^D_N,f\rangle\,\big|\,\sigma(X)\bigr)\,=\,\langle\zeta^D_N,f^{N,\ast\mathfrak e}\rangle
\end{equation}
thus proving
\begin{equation}
\label{E:7.33}
\langle\zeta^D_N,f^{N,\ast\mathfrak e}\rangle\,\,\underset{N\to\infty}\Lawarrow\,\,\langle\wt\zeta^D, f\rangle
\end{equation}
for every~$f\in C_\cc(\overline D\times(\R\cup\{+\infty\}))$. 

We will now use \eqref{E:7.33} to control the behavior of the measures $\{\zeta^D_N\colon N\ge1\}$. First, writing henceforth $1_{[M,\infty)}$ for the function $(x,h)\mapsto 1_{[M,\infty)}(h)$ we get
\begin{equation}
\bigl(1_{[M,\infty)}\bigr)^{N,\ast\mathfrak e}(x,h)=P\biggl(\,\sum_{j=1}^k (\tau_j-1)\ge (M-h)4\sqrt{2a_N}\biggr),
\end{equation}
where $k:=\lfloor 4a_N+4h\sqrt{2a_N}\rfloor$. Assuming $h\ge 2M$ with~$M>0$ large, Markov's inequality along with $E((\tau_j-1)^2)=1$ then gives
\begin{equation}
\label{E:7.41}
1-\bigl(1_{[M,\infty)}\bigr)^{N,\ast\mathfrak e}(x,h)
\le\frac{4a_N+4h\sqrt{2a_N}}{32a_N(h-M)^2}
\le\frac1{h^2}+\frac1{h\sqrt{2a_N}}.
\end{equation}
For~$M$ large, the right-hand side is at most $1/2$ thus showing
\begin{equation}
1_{[2M,\infty)}(h)
\le2\bigl(1_{[M,\infty)}\bigr)^{N,\ast\mathfrak e}(x,h).
\end{equation}
From \eqref{E:7.33} and the fact that $\wt\zeta^D$ has an exponentially decaying density in the second variable we then get, for each~$\epsilon>0$,
\begin{equation}
\label{E:7.37}
\lim_{M\to\infty}\,\limsup_{N\to\infty}\,P^{x_N}\bigl(\langle\zeta^D_N,1_{[M,\infty)}\rangle>\epsilon\bigr)=0.
\end{equation}
This implies tightness of $\{\zeta^D_N\colon N\ge1\}$ on $\overline D\times(\R\cup\{+\infty\})$ along with their asymptotic concentration on $\overline D\times\R$. In particular, we may extract a weak subsequential limit~$\zeta^D$.

We would like to use the existence of weak subsequential limits to pass to the limit~$N\to\infty$ inside the integral on the left-hand side of \eqref{E:7.33}. For that we need to deal with the fact that the support of $f^{N,\ast\mathfrak e}$ extends to~$-\infty$ in the second variable. Pick any~$b>0$ and, for any $h<-3b$, invoke Lemma~\ref{lemma-7.3} with the choices $s:=4\sqrt{2a_N}(-2b-h)$, $t:=4\sqrt{2a_N}b$
and~$k$ as above to conclude that
\begin{equation}
\bigl(1_{[-b,\infty)}\bigr)^{N,\ast\mathfrak e}(x,h)\le \texte^{-\frac{32a_N b(-2b-h)}
{4a_N- 4\sqrt{2a_N}b}}\,\,\bigl(1_{[-2b,\infty)}\bigr)^{N,\ast\mathfrak e}(x,h),\quad  h < - 3b.
\end{equation}
The prefactor decays to zero as~$h\to-\infty$ uniformly in~$N\ge1$ and so, plugging this into \eqref{E:7.33} and using that $\{\langle\zeta^D_N,(1_{[-2b,\infty)})^{N,\ast\mathfrak e}\rangle\colon N\ge1\}$ is tight we get, for each bounded, continuous~$f$ with $\supp(f)\subseteq\overline D\times[b,\infty]$ and each~$\epsilon>0$,
\begin{equation}
\lim_{M\to\infty}\,\limsup_{N\to\infty}\,P\biggl(\,\Bigl|\,\bigl\langle\zeta^D_N,\,f^{N,\ast\mathfrak e}1_{(-\infty,-M]}\bigr\rangle\Bigr|>\epsilon\biggr)=0.
\end{equation}
Combining this with \eqref{E:7.37}, we may truncate the second variable in the integral on the left of \eqref{E:7.33} to lie in $[-M,M]$ at the cost of errors that tend to zero in probability as $M\to\infty$. The Central Limit Theorem shows
\begin{equation}
\frac1{4\sqrt{2a_N}}\sum_{j\ge1} (\tau_j-1)1_{\{j\le 4a_N\}}\,\,\underset{N\to\infty}\Lawarrow\,\,\NN\bigl(0,\tfrac18\bigr)
\end{equation}
and a simple estimate based, e.g., on Doob's $L^2$-martingale inequality to account for the correction $4\sqrt{2a_N}h$ in the number of terms in the sum then gives
\begin{equation}
\lim_{N\to\infty}\,\,\sup_{h\in[-M,M]}\,\,\sup_{x\in\overline D}\,\,\bigl|f^{N,\ast\mathfrak e}(x,h)-f^{\ast\mathfrak n}(x,h)\bigr|=0,
\end{equation}
where
\begin{equation}
f^{\ast\mathfrak n}(x,h)=E\Bigl[f\bigl(x,h+\NN(0,\tfrac18)\bigr)\Bigr].
\end{equation}
Taking $M\to\infty$ after~$N\to\infty$ we then readily conclude that every subsequential weak limit~$\zeta^D$ of~$\{\zeta^D_N\colon N\ge1\}$ satisfies the distributional identity
\begin{equation}
\label{E:7.49}
\langle\zeta^D,f^{\ast\mathfrak n}\rangle\,\laweq\,\langle\wt\zeta^D,f\rangle
\end{equation}
for all~$f\in C_\cc(\overline D\times(\R\cup\{+\infty\}))$. This includes the fact that the integral on the left-hand side converges absolutely for all such~$f$.

We are now more or less done. Indeed, note that the explicit form of~$\wt\zeta^D$ gives, for $\tilde f\in C_\cc(\R)$ and~$A\subseteq D$ Borel with~$\leb(A)>0$,
\begin{equation}
\frac{\langle\wt\zeta^D,1_A\otimes \tilde f\rangle}{\langle\wt\zeta^D,1_A\otimes 1_{[0,\infty)}\rangle}  =  \alpha\lambda\int \textd h\,\texte^{-\alpha\lambda h}\tilde f(h),\quad \text{a.s.}
\end{equation}
The right-hand side is non-random and so \eqref{E:7.49} becomes the pointwise equality
\begin{equation}
\bigl\langle\zeta^D,(1_A\otimes\tilde f)^{\ast\mathfrak n}\bigr\rangle
=\bigl\langle\zeta^D,(1_A\otimes1_{[0,\infty)})^{\ast\mathfrak n}\bigr\rangle
\,\alpha\lambda\int \textd h\,\texte^{-\alpha\lambda h}\tilde f(h)
\end{equation}
for all~$\tilde f\in C_\cc(\R)$. This shows that, for any~$B\subseteq\R$ Borel,
\begin{equation}
\label{E:7.52}
\zeta^D(A\times B)=\alpha\lambda\bigl\langle\zeta^D,(1_A\otimes1_{[0,\infty)})^{\ast\mathfrak n}\bigr\rangle\otimes\nu(B),
\end{equation}
 where~$\nu$ is a Borel measure on~$\R$ that obeys \eqref{E:7.29} with $\beta:=-\alpha\lambda$ and $\sigma^2:=1/8$. Lemma~\ref{lemma-7.4} then gives $\nu(\textd h)=\texte^{-\alpha^2\lambda^2/16-\alpha\lambda h}\,\textd h$ and, since the first measure on the right of \eqref{E:7.52} has the law of the spatial part of~$\wt\zeta^D$, we get
\begin{equation}
\label{E:7.52b}
\zeta^D\,\,\laweq\,\, \texte^{-\alpha^2\lambda^2/16}\,\wt\zeta^D.
\end{equation}
The claim follows.
\end{proofsect}

Finally, we deal with the changes that are required for the thin points:

\begin{proofsect}{Proof of Theorem~\ref{thm-thin}}
Following the proof  of Theorem~\ref{thm-thick}, the argument is exactly the same up to \eqref{E:7.33}, except that now $f\in C_\cc(\overline D\times(\R\cup\{-\infty\}))$. For the tightness, we then need to consider
\begin{equation}
\bigl(1_{(-\infty,-M]}\bigr)^{N,\ast\mathfrak e}(x,h)=P\biggl(\,\sum_{j=1}^k (\tau_j-1)\le -(M+h)4\sqrt{2a_N}\biggr),
\end{equation}
where $k:=\lfloor 4a_N+4h\sqrt{2a_N}\rfloor$. For $h\le-2M$ the same estimate as \eqref{E:7.41} then shows $1_{(-\infty,-2M]}(h)\le 2(1_{(-\infty,-M]})^{N,\ast\mathfrak e}(x,h)$ and so, for each $\epsilon>0$, we get
\begin{equation}
\lim_{M\to\infty}\,\limsup_{N\to\infty}\,P^{x_N}\bigl(\langle\zeta^D_N,1_{(-\infty,-M]}\rangle>\epsilon\bigr)=0
\end{equation}
from \eqref{E:7.33}. For the upper tail, we need a variation on Lemma~\ref{lemma-7.3}:

\begin{lemma}
\label{lemma-7.5}
For $\{\tau_j\colon j\ge1\}$ i.i.d.\ Exponential(1), all $k\in\N$ and all $s,t\ge0$ with $s+t<k$,
\begin{equation}
\frac{P\Bigl(\sum_{j=1}^k(\tau_j-1)\le -(s+t)\Bigr)}{P\Bigl(\sum_{j=1}^k(\tau_j-1)\le -s\Bigr)}
\le\texte^{-\frac{t(s-1)}{k-s}}.
\end{equation}
\end{lemma}

To use this, let~$b>0$ and invoke the choices $s:=(h-2b)4\sqrt{2a_N}$, $t:=4b\sqrt{2a_N}$ and~$k$ as above while noting that, for~$N$ large and $h>2b$, we have $s+t<k$, to get
\begin{equation}
\bigl(1_{(-\infty,b]}\bigr)^{N,\ast\mathfrak e}(x,h)
\le\exp\Bigl\{-\frac{4b\sqrt{2a_N}[(h-2b)4\sqrt{2a_N}-1]}{4a_N+4h\sqrt{2a_N}-(h-2b)4\sqrt{2a_N}}\Bigr\}
\bigl(1_{(-\infty,2b]}\bigr)^{N,\ast\mathfrak e}(x,h).
\end{equation}
The exponential prefactor tends to zero as~$h\to\infty$ uniformly in~$N$ sufficiently large
and so, for any bounded and continuous~$f$ with $\supp(f)\subseteq\overline D\times(-\infty,b]$ and each~$\epsilon>0$,
\begin{equation}
\lim_{M\to\infty}\,\limsup_{N\to\infty}\,P\biggl(\,\Bigl|\,\bigl\langle\zeta^D_N,\,f^{N,\ast\mathfrak e}1_{[M,\infty)}\bigr\rangle\Bigr|>\epsilon\biggr)=0.
\end{equation}
This again permits us to truncate the tails and derive \eqref{E:7.49} for each $f\in C_\cc(\overline D\times(\R\cup\{-\infty\}))$ and each weak subsequential limit~$\zeta^D$ of $\{\zeta^D_N\colon N\ge1\}$. The rest of the proof of Theorem~\ref{thm-thick} can be followed literally leading to \eqref{E:7.52b}, as before.
\end{proofsect}

It remains to give:

\begin{proofsect}{Proof of Lemma~\ref{lemma-7.5}}
The explicit form of the density along with the substitution $y:=x+t$ again shows
\begin{equation}
\begin{aligned}
P\biggl(\,\sum_{j=1}^k(\tau_j-1)\le-(s+t)\biggr)
&=\frac1{(k-1)!}\int_{0\le x\le k-s-t}\textd x\, x^{k-1}\texte^{-x}
\\
&\le\texte^t\,\frac1{(k-1)!}\int_{t\le y\le k-s}\textd y\, (y-t)^{k-1}\texte^{-y}
\\
&
\le \texte^t\Bigl(1-\frac{t}{k-s}\Bigr)^{k-1}\,P\biggl(\,\sum_{j=1}^k(\tau_j-1)\le -s\biggr)
\end{aligned}
\end{equation}
Using the bound $1-x\le\texte^{-x}$, the prefactor is at most $\texte^{-\frac{t(s-1)}{k-s}}$.
\end{proofsect}

With the help of the above theorems, we can finally settle:

\begin{proofsect}{Proof of Theorem~\ref{thm-minmax}}
For the local time $\wh L^{D_N}_{t_N}$ parametrized by the time at the boundary vertex and the walk started at~$\varrho$, the statement appears as~\cite[Theorem~2.1]{AB}. 
The bounds in Proposition~\ref{P:tNbound} along with the tightness of $\{T_N\colon N\ge1\}$ then extend the conclusion to~$\wh L^{D_N}_{t_N}$ replaced by~$\wt L^{D_N}_{\deg(D_N) t_N}$.
Since the random walk started at~$\varrho$ visits any given $x_N\in D_N$ in time of order $N^2\log N$ while the walk started at~$x_N$ hits~$\varrho$ in time of order~$N^2$ with high probability, shifting~$t_N$ by~$\pm(\log N)^{3/2}$ and invoking the monotonicity of~$t\mapsto \wt L_t^{D_N}$ extends~\cite[Theorem~2.1]{AB} to arbitrary starting points. 
The inequalities \eqref{E:7.4} then extend it to the discrete-time object~$L^{D_N}_{t_N}$ as well.
\end{proofsect}

\section{Local structure}
\label{sec8}\noindent
The last item to be addressed are the proofs of Theorems~\ref{thm-thick-loc} and \ref{thm-avoid-loc} dealing with the local structure of the local time field near thick/thin and avoided points, respectively. We will start with the former setting, as it is technically most demanding.

\subsection{Thick and thin points}
We will again carry the argument primarily for the thick points and only comment on the changes for the thin points. Assuming henceforth the setting and notation of Theorem~\ref{thm-thick}, we start by converting the continuous-time in the boundary-vertex parametrization to that parametrized by the total time. 

\begin{proposition}
\label{prop-8.1}
Let $\overline{\zeta}_N^{D,\loc}$ be given by the same formula as $\zeta_N^{D,\loc}$ in \eqref{E:zetaNDloc} except with $L^{D_N}_{t_N}(x)$ replaced by $\overline L^{D_N}_{t_N}(x)$ from \eqref{E:7.7}. Then, given an~$x_N\in D_N$ for each~$N\ge1$, under $P^{x_N}$,
\begin{equation}
\overline{\zeta}_N^{D,\loc}\,\,\,\underset{N\to\infty}\Lawarrow\,\,\,\wt\zeta^D\otimes\wh\nu_\lambda,
\end{equation}
where $\wt\zeta^D$ is the measure on the right of \eqref{E:1.21b} and $\wh\nu_\lambda$ is the law of $\phi+\alpha\lambda\fraka$, for~$\phi$ the pinned DGFF; i.e., a centered Gaussian process on~$\Z^2$ with covariances \eqref{E:2.30}.
\end{proposition}

The proof will rely heavily on the arguments and notation from Sections~\ref{sec5}--\ref{sec7}. Throughout, we fix a sequence $\{b_N\}_{N\ge1}$ such that $b_N\to\infty$ and $b_N/t_N^{1/4}\to0$. First we condense the ideas underlying Lemmas~\ref{lemma-5.6},~\ref{lemma-5.7} and~\ref{lemma-7.2}  into:

\begin{lemma}
\label{lemma-8.2}
Given~$\epsilon>0$, let $\wt t_{N,k}^{\pm}$ be the quantity from \eqref{E:tNshift} but with $b_N$ replaced by $3b_N$. Abbreviate
\begin{multline}
\label{E:7.2a}
\quad\wt\FF_{N}(x):=\bigcup_{k\in  \Z }
\biggl(\bigl\{(k-1)\epsilon\le  T_N \circ \theta_{H_{\varrho}}  \le(k+1)\epsilon\bigr\}
\\
\cap\Bigl\{(\wh L^{D_N} \circ \theta_{H_{\varrho}})_{\wt t_{N,k}^-} (x)
\le \overline L^{D_N}_{t_N}(x)
\le \wt L^{D_N}_{H_\varrho}(x)+ (\wh L^{D_N} \circ \theta_{H_{\varrho}})_{\wt t_{N,k}^+} (x)\Bigr\}\biggr).
\end{multline}
Then for each~$b\in\R$ and any choice of~$x_N\in D_N$ for each~$N\ge1$,
\begin{equation}
P^{x_N}\Bigl(\,\sum_{x\in D_N}1_{\wt\FF_N(x)^\cc}>2\,\Bigr)\,\underset{N\to\infty}\longrightarrow\,0.
\end{equation}
\end{lemma}

\begin{proofsect}{Proof}
The tightness of~$T_N$ and $H_\varrho/|D_N|$ allows us to effectively truncate the union in \eqref{E:7.2a} to $-M\le k\le M$ and assume $H_{\varrho} \le m \deg(D_N) $.
Recall the event $\FF_N (x)$ from \eqref{E:7.2} and note that on the event
\begin{equation}
\label{E:8.2ev}
\Bigl\{\sum_{x \in D_N} 1_{\FF_N (x)^\cc} \le 2 \Bigr\} \cap \bigl\{H_{\varrho} \le m \deg (D_N) \bigr\},
\end{equation}
we have \begin{multline}
\quad
\wt L^{D_N}_{H_\varrho}+ (\wt L^{D_N} \circ \theta_{H_{\varrho}})_{(t_N +1)\deg (D_N)} (x)
\ge\overline{L}_{t_N}^{D_N} (x)
\\
\ge \wt L_{(t_N-1)\deg(D_N)}^{D_N} (x)
\ge (\wt L^{D_N} \circ \theta_{H_{\varrho}})_{(t_N - m - 1)\deg (D_N)} (x)
\quad
\end{multline}
at all but at most two $x \in D_N$.
Next set $\EE_N^+ := \EE_N (t_N + 1)$ and $\EE_N^- := \EE_N (t_N-m-1)$,
where $\EE_N (t_N')$ is the event~$\EE_N$ from \eqref{E:evtgood} but for $\{t_N\}$ replaced by $\{t_N'\}$.
Recall the notation $(t_N')^\circ$ for the quantity from \eqref{E:tNcirc}. On $\theta_{H_{\varrho}}^{-1} (\EE_N^+ \cap \EE_N^- \cap \{ (k-1)\epsilon  \le T_N \le (k+1)\epsilon\})$ we then get an analogue of \eqref{E:5.40nw} of the form
\begin{equation}
\bigl((t_N+1)^{\circ} + b_N (t_N+1)^{1/4}\bigr) \circ \theta_{H_{\varrho}} \le \wt t_{N, k}^+,
\end{equation}
\begin{equation}
\bigl((t_N-m-1)^{\circ} - b_N (t_N-m-1)^{1/4}\bigr) \circ \theta_{H_{\varrho}} \ge \wt t_{N, k}^-
\end{equation}
once~$N$ is sufficiently large (independent of~$k$).
Consequently, the inequalities
\begin{equation}
\label{E:8.9iui}
\wt L^{D_N}_{H_\varrho}(x)+ (\wh L^{D_N} \circ \theta_{H_{\varrho}})_{\wt t_{N,k}^+} (x)
\ge\overline{L}_{t_N}^{D_N} (x)
\ge (\wh L^{D_N} \circ \theta_{H_{\varrho}})_{\wt t_{N,k}^-} (x)
\end{equation}
apply on the same event as well. 
Lemma~\ref{lemma-7.2} 
shows that \eqref{E:8.9iui} holds at all but two $x\in D_N$ with $P^{x_N}$-probability tending to one as~$N\to\infty$. This proves the claim.
\end{proofsect}

Lemma~\ref{lemma-8.2} eliminates the need to consider other starting points than~$\varrho$. Next comes the main issue to be dealt with in the proof of Proposition~\ref{prop-8.1}: Since we are after differences of the local time, we cannot rely on monotonicity as we did earlier; instead we have to estimate the variation of $t\mapsto\wh L^{D_N}_t$ over time intervals of length of order~$\epsilon\sqrt{2t_N}$. This is the content of:

\begin{lemma}
\label{lemma-8.3}
For all~$\delta>0$, all $b\in\R$ and all~$\{t_N'\}_{N\ge1}$ satisfying $t_N'-t_N=O(\log N)$,
\begin{multline}
\label{E:8.9a}
\quad
\lim_{\epsilon\downarrow0}\,\limsup_{N\to\infty}\frac1{W_N}\sum_{x\in D_N}
P^\varrho\Bigl(\wh L^{D_N}_{t_N'}(x)\ge a_N+b\log N,
\\
\wh L^{D_N}_{t_N'}(x)-\wh L^{D_N}_{t_N'-\epsilon\sqrt{2t_N}}(x)>\delta\sqrt{2t_N}\,\Bigr)
=0.
\quad
\end{multline}
\end{lemma}

\begin{proofsect}{Proof}
The proof is based on tail estimates for the local time which will depend, somewhat sensitively, on a choice of a few parameters. Given~$\delta>0$ let~$\epsilon_0>0$ and~$j_0\in\N$ be such that
\begin{equation}
\label{E:8.10b}
(\sqrt\theta+\lambda)^2-(1+\epsilon_0)\theta>\lambda^2
\end{equation}
and that, for all integers~$j\ge j_0$,
\begin{equation}
\label{E:8.11b}
(j-\delta)\frac{\sqrt{\delta}-\sqrt{\epsilon_0}}{\sqrt{\delta}}>(j+1)\Bigl[\epsilon_0+\frac\lambda{\sqrt\theta+\lambda}\Bigr].
\end{equation}
These choices can be made because $(\theta+\lambda)^2-\theta^2>\lambda^2$ and $\frac\lambda{\sqrt\theta+\lambda}<1$.
Assume~$\epsilon\in(0,\epsilon_0]$ and abbreviate $t_N'':=t_N'-\epsilon\sqrt{2t_N}$ and $\wt a_N:=a_N+b\log N$.
Set $M$ to the least integer such that $(M+1)\sqrt{2t_N}\ge \wt a_N-(1+\epsilon_0)t_N''$.

Using the Markov property of $t\mapsto \wh L^{D_N}_t(x)$, the probability in \eqref{E:8.9a} is bounded by
\begin{multline}
\label{E:8.10a}
P^\varrho\Bigl(\wh L^{D_N}_{t_N''}(x)\ge\wt a_N-j_0\sqrt{2t_N}\Bigr)
P^\varrho\Bigl(\wh L^{D_N}_{\epsilon\sqrt{2t_N}}(x)\ge \delta\sqrt{2t_N}\Bigr)
\\
+
\sum_{j=j_0}^{M} P^\varrho\Bigl(\wh L^{D_N}_{t_N''}(x)\ge\wt a_N-(j+1)\sqrt{2t_N}\Bigr)
P^\varrho\Bigl(\wh L^{D_N}_{\epsilon\sqrt{2t_N}}(x)\ge j\sqrt{2t_N}\Bigr)
\\
+P^\varrho\Bigl(\wh L^{D_N}_{\epsilon\sqrt{2t_N}}(x)\ge (M+1)\sqrt{2t_N}\Bigr).
\end{multline}
We now use \cite[Lemma~4.1]{AB} to bound the individual probabilities on the right-hand side as follows. 
First, noting that by our choice of~$M$,
\begin{equation}
\sqrt{2\bigl(\wt a_N-(M+1)\sqrt{2t_N}\,\bigr)}-\sqrt{2t_N''}
\end{equation}
grows proportionally to~$\log N$ as~$N\to\infty$, \cite[Lemma~4.1]{AB} may be used for the choices $a:=\wt a_N-j_0\sqrt{2t_N}$, $t:=t_N''$ and~$b:=0$. Noting that~$W_N$ defined using $\wt a_N-j_0\sqrt{2t_N}$ and~$t_N''$ instead of~$a_N$ and~$t_N$ is comparable with~$W_N$, the uniform upper bound on~$G^{D_N}(x,x)$ then bounds the very first probability in \eqref{E:8.10a} by a quantity of order~$W_N/N^2$. The Markov inequality shows
\begin{equation}
\label{E:8.49}
P^\varrho\Bigl( \wh L_{\epsilon\sqrt{2t_N}}^{D_N} (x)  >\delta\sqrt{2a_N}\Bigr)\le\frac{\epsilon\sqrt{2t_N}}{\delta\sqrt{2a_N}}
\end{equation}
and so the first term in \eqref{E:8.10a} is order $\epsilon W_N/N^2$ (with a constant that depends on~$j_0$).
 
Next we move to the terms under the sum in \eqref{E:8.10a}. Here we use \cite[Lemma~4.1]{AB} for the choices $a:=\wt a_N$, $t:=t_N''$ and $b:=  -j\sqrt{2t_N} $ to get, for all~$j=j_0,\dots,M+1$,
\begin{equation}
\label{E:8.11a}
P^\varrho\Bigl(\wh L^{D_N}_{t_N''}(x)\ge\wt a_N-j\sqrt{2t_N}\Bigr)
\le c_1\frac{W_N}{N^2}\,\texte^{\,j\frac{\sqrt{2t_N}}{G^{D_N}(x,x)}\frac{\sqrt{2\wt a_N}-\sqrt{2t_N''}}{\sqrt{2\wt a_N}}}
\end{equation}
for some constant~$c_1\in(0,\infty)$ independent of~$N\ge1$,~$j=0,\dots,M+1$ and~$x\in D_N$. 
For the second probability under the sum in \eqref{E:8.10a}, we apply \cite[Lemma~4.1]{AB} with the choices $a:=\delta\sqrt{2t_N}$, $t:=\epsilon\sqrt{2t_N}$ and~$b:=(j-\delta)\sqrt{2t_N}$ to get
\begin{equation}
\label{E:8.12a}
P^\varrho\Bigl(\wh L^{D_N}_{\epsilon\sqrt{2t_N}}(x)\ge j\sqrt{2t_N}\Bigr)
\le c_2\,\texte^{-(j-\delta)\frac{\sqrt{2t_N}}{G^{D_N}(x,x)}\frac{\sqrt{\delta}-\sqrt{\epsilon}}{\sqrt{\delta}}}
\end{equation}
for some constant~$c_2\in(0,\infty)$ independent of~$N\ge1$ and $m\ge1$. Putting \eqref{E:8.11a} and \eqref{E:8.12a} together and invoking \eqref{E:8.11b} along with the uniform upper bound on~$G^{D_N}(x,x)$, the sum over~$j=j_0,\dots,M$ in \eqref{E:8.10a} may be performed with the result of order $ \texte^{-\alpha \sqrt{\theta} j_0\epsilon_0} W_N/N^2$, uniformly in~$x\in D_N$. 

Finally, for the stand-alone probability in \eqref{E:8.10a}, one more use of \cite[Lemma~4.1]{AB} with the choices $a:=(M+1)\sqrt{2t_N}$, $t:=\epsilon\sqrt{2t_N}$ and $b:=0$ yields
\begin{equation}
\label{E:8.16a}
P^\varrho\Bigl(\wh L^{D_N}_{\epsilon\sqrt{2t_N}}(x)\ge (M+1)\sqrt{2t_N}\Bigr)\le
\frac{c_3}{\sqrt{\log N}}\,\texte^{-(1-o(1))\frac{(M+1)\sqrt{2t_N}}{G^{D_N}(x,x)}}
\end{equation}
for a constant~$c_3\in(0,\infty)$ independent of, and $o(1)\to0$ uniformly in,~$N\ge1$ and~$x\in D_N$. Using the definition of~$M$, the right hand side of \eqref{E:8.16a} is order $N^{-2[\sqrt\theta+\lambda)^2-(1+ \epsilon_0 )\theta]+o(1)}$ which is $o(W_N/N^2)$ by $W_N=N^{2(1-\lambda^2)+o(1)}$ and \eqref{E:8.10b}, uniformly in~$x\in D_N$. The claim follows by taking~$N\to\infty$, followed by~$\epsilon\downarrow0$ and~$j_0\to\infty$. 
\end{proofsect}

We are ready to give:

\begin{proofsect}{Proof of Proposition~\ref{prop-8.1}}
Let~$f\in C_\cc( D  \times\R\times\R^{\Z^2})$ be such that~$f(x,h,\phi)$ depends only on coordinates~$\{\phi_z\colon z\in \Lambda_r(0)\}$ for some~$r>0$ and vanishes unless $|h|\le b$ and $\max_{z\in\Lambda_r(0)}|\phi_z|\le b$, for some~$b>0$. Given~$\epsilon>0$, let~$k\in\Z$ be such that $| T_N \circ \theta_{H_{\varrho}}  -k\epsilon|<\epsilon$. Pick~$x\in D_N$ and abbreviate
\begin{equation}
f_{N,r}(x,\ell):=f\biggl(x/N,  \frac{\ell(x)-a_N}{\sqrt{2a_N}} , 
\Bigl\{ \frac{\ell(x)-\ell(x+z)}{\sqrt{2a_N}}  \colon z\in\Lambda_r(0)\Bigr\}\biggr).
\end{equation}
Introducing the oscillation of~$f$ by
\begin{equation}
\text{osc}_{f} (\delta)
:= \sup_{x \in D} \,\sup_{\begin{subarray}{c} u, v \in \R,\\ |u-v| \le \delta \end{subarray}}
\,\,\sup_{\begin{subarray}{c} \phi, \wt \phi \in \R^{\Lambda_r (0)}, \\
\max_{z\in\Lambda_r(0)}|\phi_z - \wt \phi_z| \leq 2 \delta \end{subarray}}
\bigl|f(x, u, \phi) - f(x, v, \wt \phi)\bigr|,
\end{equation}
the difference 
\begin{equation}
\label{E:8.17a}
f_{N,r}\bigl(x, \overline L^{D_N}_{t_N}\bigr)- f_{N,r}\Bigl(x, (\wh L^{D_N}\circ  \theta_{H_\varrho} )_{\wt t_{N,k}^-}\Bigr)
\end{equation}
is bounded in absolute value by the sum over~$z\in\Lambda_r(x)$ of three terms: $2\Vert f\Vert_\infty 1_{\wt\FF_N(z)^\cc}$,
\begin{equation}
\label{E:8.20a}
2\Vert f\Vert_\infty 1_{\wt\FF_N(z)\cap\{H_z<H_\varrho\}}
\Bigl(1_{\{(\wh L^{D_N}\circ  \theta_{H_\varrho} )_{\wt t_{N,k}^-}(z)\ge a_N-2b\sqrt{2a_N}\}}+1_{\{\overline L^{D_N}_{t_N}(z)\ge a_N-2b\sqrt{2a_N}\}}\Bigr)
\end{equation}
and
\begin{multline}
\label{E:8.21a}
\qquad
1_{\wt\FF_N(z)\cap\{H_z>H_\varrho\}}\Bigl(\text{osc}_f(\delta)+\Vert f\Vert_\infty 1_{\{|\overline L^{D_N}_{t_N}(z)-  (\wh L^{D_N} \circ \theta_{H_{\varrho}})_{\wt t_{N,k}^-}  (z)|>\delta\sqrt{2a_N}\}}\Bigr)
\\
\times\Bigl(1_{\{ (\wh L^{D_N} \circ \theta_{H_{\varrho}})_{\wt t_{N,k}^-} (z)\ge a_N-2b\sqrt{2a_N}\}}+1_{\{\overline L^{D_N}_{t_N}(z)\ge a_N-2b\sqrt{2a_N}\}}\Bigr).
\qquad
\end{multline}
To simplify estimates, introduce the events
\begin{equation}
\GG_N(x):=\Bigl\{\wt L^{D_N}_{H_\varrho}(x)+ (\wh L^{D_N} \circ \theta_{H_{\varrho}})_{\wt t_{N,k}^+} (x)\ge a_N-2b\sqrt{2a_N}\Bigr\}\cap \{H_x<H_\varrho\}
\end{equation}
and
\begin{equation}
\HH_N(x):=\Bigl\{\wh L^{D_N}_{\wt t_{N,k}^+}(x)\ge a_N-2b\sqrt{2a_N}\Bigr\}\cap\Bigl\{\wh L^{D_N}_{\wt t_{N,k}^+}(x)-\wh L^{D_N}_{\wt t_{N,k}^-}(x)>\delta\sqrt{2a_N}\Bigr\}.
\end{equation}
Then \eqref{E:8.20a} is bounded by $4\Vert f\Vert_\infty 1_{\GG_N(z)}$ while \eqref{E:8.21a} is bounded by 
\begin{equation}
2\text{osc}_f(\delta)1_{\{(\wh L^{D_N}\circ  \theta_{H_\varrho} )_{\wt t_{N,k}^+}(z)\ge a_N-2b\sqrt{2a_N}\}}+2\Vert f\Vert_\infty 1_{\HH_N(z)}\circ \theta_{H_\varrho}.
\end{equation}
Summarizing these estimates, and writing~$\wh\zeta^{D,\loc}_N(t_N')$ for the measure in \eqref{E:zetaNDloc} except with $L^{D_N}$ replaced by $\wh L^{D_N}$ and~$t_N$ by~$t_N'$, we thus get that, on $\{| T_N \circ \theta_{H_{\varrho}}  -k\epsilon|<\epsilon\}$,
\begin{multline}
\label{E:8.25a}
\biggl|\langle\overline\zeta^{D,\loc}_N,f\rangle-\frac{W_N(\wt t_{N,k}^-)}{W_N}
\bigl\langle\wh\zeta_N^{D,\loc}(\wt t_{N,k}^-),f\bigr\rangle\circ \theta_{H_\varrho}\biggr|
\\
\le 4\Vert f\Vert_\infty|\Lambda_r(0)|\frac1{W_N}\sum_{x\in D_N}\bigl(1_{\wt\FF_N(x)^\cc}+1_{\GG_N(x)}+1_{\HH_N(x)}\circ \theta_{H_\varrho}\bigr)
\\
+2\,\text{osc}_f(\delta)|\Lambda_r(0)|\frac{W_N( \wt t_{N,k}^+ )}{W_N}\bigl\langle\wh\zeta^{D}_N( \wt t_{N,k}^+ ),1_D\otimes 1_{[-2b,\infty)}\bigr\rangle\circ \theta_{H_\varrho}
\end{multline}
Using Lemmas~\ref{lemma-8.2}, \ref{lemma-8.3} and \ref{lemma-6.3}, the first term on the right tends to zero in $P^{x_N}$-probability as~$N\to\infty$  and $\epsilon \downarrow 0$  for each~$\delta>0$. 
The tightness of~$\wh\zeta^D_N$ measures (under~$P^\varrho$) along with the uniform continuity of~$f$ ensure that the second term tends to zero in $P^{x_N}$-probability as~$N\to\infty$ and~$\delta\downarrow0$. 

To finish the proof, note that by \cite[Theorem~2.6]{AB} and the argument underlying Proposition~\ref{thm-4.3} we have, under~$P^\varrho$,
\begin{equation}
\wh\zeta^{D,\loc}_N(t_N')\otimes\delta_{T_N}\,\,\,\underset{N\to\infty}\Lawarrow\,\,\,\wh\zeta^D\otimes\wh\nu_\lambda\otimes\delta_{T}
\end{equation}
for any sequence~$\{t_N'\}_{N\ge1}$ such that $t_N'-t_N=o(t_N)$, where $\wh\zeta^D$ is related to~$T$ as in \eqref{E:5.48a}. 
 Since $W_N(\wt t_{N,k}^-)/W_N=(\texte^{-\alpha\lambda T_N (\wt t_{N, k}^-)}\circ 
\theta_{H_\varrho})\texte^{O(\epsilon)}$ on $\{|T_N \circ \theta_{H_{\varrho}}-k\epsilon|<\epsilon\}\cap\mathcal{E}_N^- \circ \theta_{H_{\rho}}$,
from \eqref{E:8.25a} and the tightness of the random variables~$\{T_N\}_{N\ge1}$ and $\{H_\varrho/|D_N|\}_{N\ge1}$ 
we get, by taking~$N\to\infty$ followed by~$\delta\downarrow0$,~$\epsilon\downarrow0$ and~$m\to\infty$, under~$P^{x_N}$,
\begin{equation}
\overline\zeta^{D,\loc}_N\,\,\,\underset{N\to\infty}\Lawarrow\,\,\,\texte^{-\alpha\lambda T}\wh\zeta^D\otimes\wh\nu_\lambda.
\end{equation}
This is the desired claim.
\end{proofsect}

With Proposition~\ref{prop-8.1} in hand, we are ready to tackle:

\begin{proofsect}{Proof of Theorem~\ref{thm-thick-loc}, thick points}
First observe that the tightness of $\{\zeta_N^{D}\colon N\ge1\}$ 
implies tightness of 
$\{\zeta_N^{D,\loc}\colon N\ge1\}$ and so we may consider subsequential distributional limits $\zeta^{D,\loc}$ of the latter. Using Proposition~\ref{prop-8.1} in the argument from the proof of Theorem~\ref{thm-thick} we conclude that every such subsequential weak limit obeys
\begin{equation}
\langle \zeta^{D,\loc}, f^{\ast\mathfrak n} \rangle 
\,\laweq\,
\langle \wt \zeta^{D}\otimes\wh\nu_\lambda,f\rangle
\end{equation}
for all $f\in C_\cc(D\times\R\times\R^{\Z^2})$, where
\begin{equation}
f^{\ast\mathfrak n}(x,h, \phi):=E\Bigl[f\bigl(x,h+\mathfrak{n}_0, \{\mathfrak{n}_0 
- \mathfrak{n}_z + \phi_z \colon z \in \mathbb{Z}^2 \}\bigr)\Bigr],
\end{equation}
for $\{\mathfrak{n}_z \colon z \in \mathbb{Z}^2\}$ i.i.d. $\NN(0,\tfrac18)$. 

We now proceed similarly as in \twoeqref{E:7.49}{E:7.52}: Given any $\tilde f\in C_\cc(\R \times \R^{\Z^2})$ and any Borel~$A\subseteq D$ with~$\leb(A)>0$, the explicit form of~$\wt\zeta^{D,\loc}$ gives the pointwise equality
\begin{multline}
\bigl\langle\zeta^{D,\loc},(1_A\otimes\tilde f)^{\ast\mathfrak n}\bigr\rangle
\\
=\bigl\langle\zeta^{D,\loc},(1_A\otimes1_{[0,\infty)}\otimes 1_{\R^{\Z^2}})^{\ast\mathfrak n}\bigr\rangle
\,\alpha\lambda\int \textd h\,\texte^{-\alpha\lambda h} \otimes \wh\nu_{\lambda} (\textd \phi) \tilde f(h, \phi).
\end{multline}
Abbreviating $\beta:=-\alpha\lambda$, for each~$A$ as above, the measure $\zeta_A$ on~$\R\times\R^{\Z^2}$ defined by
\begin{equation}
\label{E:8.31a}
\zeta_A(B):=\frac{\zeta^{D,\loc}(A\times B)}
{\alpha\lambda \bigl\langle\zeta^{D,\loc},(1_A\otimes1_{[0,\infty)}\otimes 1_{\R^{\Z^2}})^{\ast\mathfrak n}\bigr\rangle}
\end{equation}
then ``solves'' for~$\mu$ from the convolution equation
\begin{multline}
\label{E:0.24}
\quad
\int_{\R \times \R^{\Z^2}}
\mu(\textd h \textd \phi)
\,E\bigl[ \,f(h+\mathfrak n_0, \{\mathfrak n_0 - \mathfrak n_z + \phi_z \colon z \in \Z^2\})\bigr] 
\\= \int_{\R \times \R^{\Z^2}} \textd h\,\,\texte^{\beta h} \otimes \wh\nu_{\lambda} (\textd \phi)\,f(h, \phi)
\quad
\end{multline}
for all $f\in C_\cc(\R\times\R^{\Z^2})$.  
To solve this equation, we need:

\begin{lemma}
\label{lemma-8.4a}
For each~$x,y\in\Z^2$, let
\begin{equation}
\label{E:8.27}
\wt C(x,y):=\fraka(x)+\fraka(y)-\fraka(x-y)-\frac18\bigl[1-\delta_{x,0}-\delta_{y,0}+\delta_{x,y}\bigr].
\end{equation}
Then $\wt C$ is symmetric and positive semidefinite and so there exists a centered Gaussian process~$\{\wt\phi_x\colon x\in\Z^2\}$ with covariance~$\wt C$. This process then satisfies \eqref{E:2.33}.
\end{lemma}

\begin{proofsect}{Proof}
Recall that (in our normalization) $\fraka$ solves the equation $\Delta\fraka = \delta_0$ and so using Fourier transform techniques we get
\begin{equation}
\fraka(x)=\int_{(-\pi,\pi)^2}\frac{\textd k}{(2\pi)^2}\frac{1-\texte^{-\texti k\cdot x}}{\wh D(k)},
\end{equation}
where 
\begin{equation}
\wh D(k):=4\sin(k_1/2)^2+4\sin(k_2/2)^2.
\end{equation}
Let $v\in\ell^2(\Z^2)$ and denote by $\hat v(k):=\sum_{x\in\Z^2}v(x)\texte^{\texti k\cdot x}$ the Fourier transform of~$v$. A calculation then shows
\begin{equation}
(v,\wt C v) = \int_{(-\pi,\pi)^2}\frac{\textd k}{(2\pi)^2}\biggl(\frac1{\wh D(k)}-\frac18\biggr)\bigl|\hat v(0)-\hat v(k)\bigr|^2
\end{equation}
Noting that $\wt D(k)\le8$, we get that $\wt C$ is indeed positive semidefinite. We now readily check that $x,y\mapsto\frac18[1-\delta_{x,0}-\delta_{y,0}+\delta_{x,y}]$ is the covariance of $\{n_0-n_z\colon z\in\Z^2\}$ for $\{n_z\colon z\in\Z^2\}$ i.i.d.\ $\NN(0,\frac18)$, and so \eqref{E:2.33} holds as well.
\end{proofsect}

The solution of \eqref{E:0.24} will require the following extension of Lemma~\ref{lemma-7.4}:

\begin{lemma}
\label{lemma-8.4}
Let $\wt\phi$ be a centered Gaussian process on~$\Z^2$ such that, for some~$\beta\in\R$ and some~$\sigma^2>0$, the process $\{\wt\phi_x+n_0-n_z\colon z\in\Z^2\}$ with $\{n_z\colon z\in\Z^2\}$ i.i.d.\ $\NN(0,\sigma^2)$ has the law of the pinned DGFF~$\phi$. Denote
\begin{equation}
\nu_{\lambda,\beta}(A):=P\biggl(\wt\phi+\lambda\alpha\fraka+\beta\sigma^21_{\Z^2\smallsetminus\{0\}}\in A\biggr).
\end{equation}
Then \eqref{E:0.24} is solved uniquely by
\begin{equation}
\label{E:0.29}
\mu(\textd h\textd\phi) = \texte^{-\frac1{2}\beta^2\sigma^2+\beta h}\textd h\otimes\nu_{\lambda,\beta}(\textd\phi).
\end{equation}
\end{lemma}

\begin{proofsect}{Proof}
Denote $\wt\mu(\textd h\textd \phi):=\texte^{\frac1{2}\beta^2\sigma^2-\beta h}\mu(\textd h\textd \phi)$.
Pick $\{t_z\colon z\in\Z^2\}$ with finite support and~$t_0=0$ and, writing $\langle\cdot,\cdot\rangle$ for the inner product in $\ell^2(\Z^2)$, apply \eqref{E:0.24} to the test function $h,\phi\mapsto \texte^{-\beta h}\,f(h)\exp\{\langle t,\phi\rangle\}$ with a non-negative~$f\in C_\cc(\R)$. (This is permissible in light of the Monotone Convergence Theorem.) Writing~$x$ for~$h+n_0$ then turns \eqref{E:0.29} into
\begin{multline}
\label{E:0.30}
\int\wt\mu(\textd h\textd \phi)\otimes\textd x\,\texte^{\langle t,\phi\rangle}\,\frac1{\sqrt{2\pi\sigma^2}}\,\texte^{-\frac1{2\sigma^2}(x-h)^2}E\bigl(\texte^{-\langle t,n\rangle})\,\texte^{\bar t(x-h)}\texte^{-\frac1{2}\beta^2\sigma^2+\beta h}\texte^{-\beta x}\,f(x)
\\
=\int \wh \nu_\lambda(\textd\phi)\otimes\textd x\,\texte^{\langle t,\phi\rangle}f(x)
\end{multline}
where $\bar t:=\sum_{z\in\Z^2}t_z$.
By assumption we have 
\begin{equation}
\{\phi_z\colon z\in\Z^d\}\,\laweq\,\{\wt\phi_z+n_0-n_z\colon z\in\Z^2\}
\end{equation}
and so, in light of $t_0=0$,
\begin{equation}
\begin{aligned}
\int \wh \nu_\lambda(\textd\phi)\texte^{\langle t,\phi\rangle}
&=\int P(\textd\phi)\texte^{\langle t,\phi+\alpha\lambda\fraka\rangle}\\
&=\int P(\textd\wt\phi)E\bigl(\texte^{\langle t,\wt\phi+n_0-n+\alpha\lambda\fraka\rangle})\\
&=\int  \nu_{\lambda,\beta}  (\textd\wt\phi)\,\texte^{\langle t,\wt\phi\rangle}\,E\bigl(\texte^{-\langle t,n\rangle}\bigr)E(\texte^{\bar t (n_0-\beta\sigma^2)}),
\end{aligned}
\end{equation}
where the expectation is over $\{n_z\colon z\in\Z^2\}$.
Using this in \eqref{E:0.30} and cancelling $E\bigl(\texte^{-\langle t,n\rangle})$ on both sides, the identity $E(\texte^{\bar t (n_0-\beta\sigma^2)}) = \texte^{\frac12\bar t^2\sigma^2-\beta\bar t\sigma^2}$ along with the fact that functions~$f\in C_\cc(\R)$ separate points yield
\begin{multline}
\label{E:0.33}
\int\wt\mu(\textd h\textd \phi)\,\texte^{\langle t,\phi\rangle}\,\frac1{\sqrt{2\pi\sigma^2}}\,\texte^{-\frac1{2\sigma^2}(x-h)^2}\,\texte^{\bar t(x-h)}\texte^{-\beta x} \texte^{-\frac12\bar t^2\sigma^2+\beta\bar t\sigma^2}\texte^{-\frac1{2}\beta^2\sigma^2+\beta h}
\\
=\int  \nu_{\lambda,\beta}  (\textd\wt\phi)\,\texte^{ \langle t,\wt\phi \rangle }\,
\end{multline}
for all~$x\in\R$. (Continuity is used to get from Lebesgue a.e.~$x\in\R$ to all~$x\in\R$.)
The five exponentials on the left combine into
\begin{equation}
\texte^{-\frac1{2\sigma^2}(x-h-\bar t\sigma^2)^2 -\beta (x-h-\bar t\sigma^2)-\frac1{2}\beta^2\sigma^2}=\texte^{-\frac1{2\sigma^2}(x-h-\bar t\sigma^2+\beta\sigma^2)^2}.
\end{equation}
Shifting~$x$ by $\bar t\sigma^2+\beta\sigma^2$ and scaling it by~$\sigma^2$ shows that $\wh\mu(\textd h\textd\phi):=\frac1{\sqrt{2\pi\sigma^2}}\texte^{-\frac1{2\sigma^2}h^2}\wt\mu(\textd h\textd\phi)$ obeys
\begin{equation}
\label{E:8.45}
\int\wh\mu(\textd h\textd \phi)\,\texte^{\langle t,\phi\rangle-xh}
=\int\nu_{\lambda,\beta}(\textd\wt\phi)\,\texte^{\langle t,\wt\phi\rangle}\texte^{\frac12x^2\sigma^2}
\end{equation}
for all $x\in\R$ and all~$\{t_z\colon z\in\Z^2\}$ with finite support and~$t_0=0$. 

The restriction to~$t_0=0$ is irrelevant in \eqref{E:8.45} since $\nu_{\lambda,\beta}$ is concentrated on $\{\phi\colon \phi_0=0\}$ and, by \eqref{E:0.24} so is~$\mu$ and thus also~$\wh\mu$. The right-hand side of \eqref{E:8.45} is the Laplace transform of the product of the law of~$\NN(0,\sigma^2)$ and~$\nu_{\lambda,\beta}$. Hence
\begin{equation}
\wt\mu(\textd h\textd \phi) = \textd h\otimes \nu_{\lambda,\beta}(\textd\phi)
\end{equation}
and so the claim follows from the definition of~$\wt\mu$.
\end{proofsect}

Returning to the main line of the proof of Theorem~\ref{thm-thick-loc}, it remains to observe that the denominator in \eqref{E:8.31a} has the law of
\begin{equation}
\sqrt{\frac{\sqrt{\theta}}{\sqrt{\theta}+\lambda}}\,\,\cspecial(\lambda)\,\ee^{\alpha \lambda (\mathfrak{d}(x) - 1) Y}\,Z_\lambda^{D,0}(\textd x),
\end{equation}
for~$Y=\NN(0,\sigma_D^2)$ independent of~$Z^{D,0}_\lambda$.
Lemma~\ref{lemma-8.4} with $\beta:=-\alpha\lambda$ and~$\sigma^2:=\frac18$ then yields the claim.
\end{proofsect}

Moving to the thin points, here we go directly for:

\begin{proofsect}{Proof of Theorem~\ref{thm-thick-loc}, thin points}
The proof is considerably simpler because, as a few times earlier, certain key inequalities go in a more favorable direction. Following the argument and the notation from the proof for the thick points, we derive an analogue of \eqref{E:8.25a} with the events $\GG_N(x)$ and $\HH_N(x)$ replaced by
\begin{equation}
\wt\GG_N(x):=\Bigl\{(\wh L^{D_N} \circ \theta_{H_{\varrho}})_{\wt t_{N,k}^-} (x)\le a_N+2b\sqrt{2a_N}\Bigr\}\cap \{H_x<H_\varrho\}
\end{equation}
and
\begin{equation}
\wt\HH_N(x):=\Bigl\{\wh L^{D_N}_{\wt t_{N,k}^-}(x)\le a_N+2b\sqrt{2a_N}\Bigr\}\cap\Bigl\{\wh L^{D_N}_{\wt t_{N,k}^+}(x)-\wh L^{D_N}_{\wt t_{N,k}^-}(x)>\delta\sqrt{2a_N}\Bigr\},
\end{equation}
respectively,
and $1_{[-2b,\infty)}$ replaced by $1_{(-\infty,2b]}$. The $P^{x_N}$-probability of event~$\wt \GG_N(x)$ is controlled using Lemma~\ref{lemma-6.5}. Unlike $\HH_N(x)$ which required a non-trivial decomposition in the proof of Lemma~\ref{lemma-8.3}, the two events constituting~$\wt\HH_N(x)$ can be directly separated using the Markov property of~$t\mapsto\wh L^{D_N}_t$. The expected sum over~$1_{ \wt \HH_N(x) }\circ\theta_{H_\varrho}$ is then shown to be order 
$\epsilon W_N$ by \eqref{E:8.49}
and the fact that $E^\varrho\langle\wh\zeta^D_N(\wt t_{N,k}^-),1_{(-\infty,2b]}\rangle$ is bounded in~$N\ge1$. As a consequence, we get that, under $P^{x_N}$,
\begin{equation}
\overline\zeta^{D,\loc}_N\,\,\,\underset{N\to\infty}\Lawarrow\,\,\,\wt\zeta^D\otimes\wh\nu_\lambda,
\end{equation}
where~$\wt\zeta^D$ is the measure on the right of \eqref{E:1.23dis} without the term $\texte^{-\alpha^2\lambda^2/16}$ and~$\wh\nu_\lambda$ is the law of~$\phi-\alpha\lambda\fraka$.  
The rest of the argument for the thick points may be followed literally.
\end{proofsect}

\subsection{Avoided points}
The proof is a variation on the themes encountered in the proof of convergence of the measure associated with the light and avoided points. In particular, since the local time vanishes at the avoided points, we will be able to use monotonicity arguments. The following observation will be useful:

\begin{lemma}
\label{lemma-8.6}
Let~$\mu$ be a probability measure on~$\N^{\Z^2}$ with samples denoted by $\{\hat n_z\colon z\in\Z^2\}$. Let $\{\tau_j(x)\colon j\ge1,\,x\in\Z^2\}$ be i.i.d.\ Exponential(1), independent of $\{\hat n_z\colon z\in\Z^2\}$. Then for any~$t\in(-1,\infty)^{\Z^2}$ with finite support,
\begin{equation}
E\exp\Bigl\{-\sum_{z\in\Z^2}t(z)\sum_{j=1}^{\hat n_z}\tau_j(z)\Bigr\}=E\exp\Bigl\{-\sum_{z\in\Z^2}t'(z)\hat n_z\Bigr\},
\end{equation}
where $t'(z):=\log(1+t(z))$.
\end{lemma}

\begin{proofsect}{Proof}
This boils down to a calculation of the Laplace transform of Exponential(1).
\end{proofsect}

\begin{proofsect}{Proof of Theorem~\ref{thm-avoid-loc}}
We will establish the existence and uniqueness of the law $\nu_u^{\text{RI},\text{dis}}$ as part of the proof of the convergence. Let~$\tilde f\in C(\overline D)$ be non-negative, pick 
 $t\in(0,\infty)^{\Z^2}$
with finite support and consider the test function
\begin{equation}
f_t(x,\phi):=\tilde f(x)\texte^{-\langle t,\phi\rangle}
\end{equation}
where, abusing notation as before, $\langle\cdot,\cdot\rangle$ denotes the canonical inner product in~$\ell^2(\Z^2)$. The function $x,h,\phi\mapsto \texte^{-hn}f_t(x,\phi)$ is  non-increasing 
in both~$h$ and the coordinates of~$\phi$ and so, thanks to Lemma~\ref{lemma-8.2}, \eqref{E:5.35u} applies to~$f$ replaced by~$\texte^{-hn}f_t$ and $\wt\vartheta^D_N$ by
\begin{equation}
\overline\vartheta^D_N:=\frac1{\wh W_N}\sum_{x\in D_N}\delta_{x/N}\otimes\delta_{\overline L^{D_N}_{t_N}(x)}\otimes\delta_{\{\overline L^{D_N}_{t_N}(x+z)\colon z\in\Z^2\}}.
\end{equation}
Let~$\overline\kappa^D_N$ be the measure tracking the local behavior of $\overline L^{D_N}_{t_N}(x+z)\colon z\in\Z^2$ around every point~$x$ where $\overline L^{D_N}_{t_N}(x)=0$ which, we note, is almost surely equivalent to~$L^{D_N}_{t_N}(x)=0$. Taking the limits~$N\to\infty$ and~$n\to\infty$, from \cite[Theorem~2.8]{AB} we then get, under~$P^{x_N}$,
\begin{equation}
\label{E:8.54}
\langle\overline\kappa_N^{D,\loc},f_t\rangle \,\,\,\underset{N\to\infty}\Lawarrow\,\,\,
\langle\wt\kappa^D\otimes\nu_\theta^{\text{RI}},f_t\rangle,
\end{equation}
where~$\wt\kappa^D$ is the law on the right-hand side of \eqref{E:2.27cont}. 

Next we observe that, by Lemma~\ref{lemma-8.6} and the fact that~$4L^{D_N}_{t_n}(x)$ is a natural,
\begin{equation}
E^\varrho\bigl(\langle\overline\kappa_N^{D,\loc},f_t\rangle\,\big|\,\sigma(X)\bigr) = \langle\kappa_N^{D,\loc},f_{t'}\rangle
\end{equation}
where $t'(z):=4\log(1+t(z)/4)$. From \eqref{E:7.10} and \eqref{E:8.54} we then get that every subsequential weak limit $\kappa^{D,\loc}$ of $\{\kappa_N^{D,\loc}\colon N\ge1\}$ obeys 
\begin{equation}
\langle\kappa^{D,\loc},f_{t'}\rangle \laweq \langle\wt\kappa^D\otimes\nu_\theta^{\text{RI}},f_t\rangle
\end{equation}
jointly for all~$t\in(0,\infty)^{\Z^2}$ with finite support and all~$\tilde f\in C(\overline D)$.
Since~$\nu_\theta^{\text{RI}}$ is non-random, this is readily turned into the a.s.\ identity
\begin{equation}
\int \kappa^{D,\loc}(\textd x\textd\ell) \tilde f(x)\texte^{-\langle t',\ell\rangle} = \Bigl(\int\wt\kappa^D(\textd x)\tilde f(x)\Bigr)\int \nu_\theta^{\text{RI}}(\textd\phi)\texte^{-\langle t,\phi\rangle}.
\end{equation}
This along with the fact that
\begin{equation}
\texte^{-\langle t',\ell\rangle} = E\exp\biggl\{-\sum_{z\in\Z^2}t(z)\frac14\sum_{j=1}^{4\ell(z)} \tau_j(z)\biggr\}
\end{equation}
for $\{\tau_j(z)\colon j\ge1,\,z\in\Z^2\}$ independent i.i.d.\ Exponential(1) implies that
\begin{equation}
\kappa^{D,\loc} = \wt\kappa^D\otimes\nu_\theta^{\text{RI},\text{dis}}
\end{equation}
 where~$\nu_\theta^{\text{RI},\text{dis}}$ is a measure as described in the statement. 

This shows that a measure $\nu_u^{\text{RI},\text{dis}}$ exists with the stated properties for all~$u\in(0,1)$. Since adding independent samples from this measure for parameters~$u\in(0,1)$ and~$v\in(0,1)$ gives us a sample from the measure for parameter~$u+v$, the existence extends to all~$u>0$. The measure is unique by Lemma~\ref{lemma-8.6} and so is thus the distributional limit~$\kappa^{D,\loc}$. This completes the proof.
\end{proofsect}

\section*{Acknowledgments}
\nopagebreak\nopagebreak\noindent
This project has been supported in part by the NSF award DMS-1712632 and GA\v CR project P201/16-15238S.
The first author has been supported in part by JSPS KAKENHI, Grant-in-Aid for Early-Career Scientists 18K13429.

\bibliographystyle{abbrv}

\end{document}